\documentclass[12pt, oneside]{article}
\usepackage{amsmath,amsfonts,amssymb,amscd,theorem}
\date{ }

\newtheorem{thm}{Theorem}
\newtheorem{lem}{Lemma}
\newtheorem{defi}{Definition}
\newtheorem{cor}{Corollary}
\newtheorem{prop}{Proposition}

\begin{document}

\title{The cohomology of a variation of polarized Hodge structures
\\
over a quasi-compact K\"ahler manifold}

\author{J. Jost, Y.-H. Yang\thanks{The author supported partially by NSF of
China (No.10171077) and a direct grant (Title: Cohomology theory
with coefficients in local systems and applications to algebraic
geometry) from the Chinese University of Hong Kong by the third
named author;} ~and K. Zuo\thanks{The author supported partially
by a grant from the Research Grants Council of the Hong Kong
Special Administrative Region, China (Project No. CUHK
4034/02P).}}

\maketitle

\section{Introduction}
Let $(\overline M, \omega_0)$ be a compact K\"ahler manifold of
dimension $n$ and $D$ a divisor of $\overline M$ with at most
normal crossing singularities . Let $D = \cup_{i=1}^lD_i$ , where
the $D_i$ are smooth divisors of $\overline M$. Denote ${\overline
M}\setminus D$ by $M$, called a quasi-compact K\"ahler manifold.
$j:M\to{\overline M}$ is the inclusion mapping. According to
\cite{cg}, one can construct a complete K\"ahler metric $\omega$
on $M$ which is a Poincar\'e-like metric near the divisor and of
finite volume.

Let $\{M, {\bf H}_{\mathbb{C}}={\bf
H}_{\mathbb{Z}}\otimes{\mathbb{C}},
\nabla=\nabla^{1,0}+\nabla^{0,1} , {\bf F}=\{{\bf F}^p\}_{p=1}^m,
{\bf S}\}$ be a rational variation of polarized Hodge structures
with weight $m$ over $M$ such that each ${\bf F}^p$ is a
holomorphic subbundle of the local system ${\bf H_{\mathbb{C}}}$
and $\nabla^{1,0}$ satisfies the Griffiths' infinitesimal period
relation
\[
\nabla^{1,0}{\bf F}^p\subset\Omega^1_M\otimes{\bf F}^{p-1}.
\]
For simplicity, we always assume that all monodromies at the
divisor are {\it unipotent}. Note that by a lemma of Borel (cf.
\cite{sch}, Lemma 4.5), this is not an essential assumption. The
polarization $\bf S$ defines a Hermitian metric
\[
h(\cdot, \cdot)={\bf S}(\mathcal{C}\cdot, \overline{\cdot})
\]
on ${\bf H}_{\mathbb{C}}$ by using the Weil operator
$\mathcal{C}$, called {\it Hodge metric}, where the bar is the
conjugation relative to ${\bf H}_{\mathbb{Z}}$. Take the
successive quotients ${\bf F}^p/{\bf F}^{p-1}$, called {\it Hodge
bundles}, denoted by ${\bf E}^p$. Accordingly, we have the induced
$\cal{O}$-linear map $\theta^p: {\bf E}^p\to{\bf E}^{p-1}$ of
$\nabla^{1,0}$. Set ${\bf E}=\oplus{\bf E}^p$ and
$\theta=\oplus\theta^p$. Since $(\nabla^{1,0})^2=0$,
$\theta\wedge\theta=0$. So $({\bf E}, \theta)$ is a Higgs bundle.

Using the asymptotic behavior of $\theta$ (cf. \S3), one has the
following holomorphic Dolbeault complex of sheaves on $\overline
M$
\[
{\bf E}\stackrel{\theta}{\to}{\bf E}\otimes\Omega^1_{\overline
M}(\log D) \stackrel{\theta}{\to}{\bf E}\otimes\Omega^2_{\overline
M}(\log D)\stackrel{\theta}{\to}\cdots.
\]
Furthermore, using the Poincar\'e-like metric $\omega$, the Hodge
metric $h$ induced on $\bf E$, and the boundedness of $\theta$
(cf. \S3), we can define an $L^2$-subcomplex on $\overline M$
\[
(*)~~~~~~~~~~~~~~~~{\bf E}_{(2)}\stackrel{\theta}{\to}({\bf
E}\otimes\Omega^1_{\overline M}(\log D))_{(2)}
\stackrel{\theta}{\to}({\bf E}\otimes\Omega^2_{\overline M}(\log
D))_{(2)}\stackrel{\theta}{\to}\cdots.
\]
of the above complex by taking the sheaves of germs of local $L^2$
sections. Although the definition of the above $L^2$-subcomplex
depends on the metric on $M$ and the Hodge metric on ${\bf E}$, we
will however see that it is essentially independent of both
metrics and is obtained by taking local sections satisfying a
certain algebraic condition determined by the monodromies of the
variation from the arguments in \S4. One can then consider the
hypercohomology of the $L^2$-subcomplex
$$
\mathbb{H}^*({\overline M}, \{({\bf
E}\otimes\Omega^{\cdot}_{\overline M}(\log D))_{(2)}, \theta\}),
$$
which is independent of both metrics and depends only on the
monodromies of the variation.

On the other hand, by means of a construction of Deligne, we can
also define a complex of fine sheaves as follows (for details, cf.
\S 5). Define $[{\text{Gr}}^*_FA^k({\bf H}_{\mathbb{C}})]_{(2)}$
to be the sheaf obtained by taking local $L^2$ measurable
$k$-forms on $\overline M$ with values in $\bf E$, the
$\overline\partial$ derivatives of which exist in the weak sense
and are $L^2$; and put $D''=\overline\partial+\theta$. Then
$D''([{\text{Gr}}^*_FA^k({\bf
H}_{\mathbb{C}})]_{(2)})\subset[{\text{Gr}}^*_FA^{k+1}({\bf
H}_{\mathbb{C}})]_{(2)}$ and $(D'')^2=0$ due to $\nabla$ being
flat. Thus we have the complex of fine sheaves on $\overline M$
\[
(**)~~~~~~~~~~~~~[{\text{Gr}}^*_FA^0({\bf
H}_{\mathbb{C}})]_{(2)}\stackrel{D''}{\to}
[{\text{Gr}}^*_FA^1({\bf
H}_{\mathbb{C}})]_{(2)}\stackrel{D''}{\to}
[{\text{Gr}}^*_FA^2({\bf
H}_{\mathbb{C}})]_{(2)}\stackrel{D''}{\to}\cdots.
\]
It is not difficult to see that the complex $(*)$ is actually a
subcomplex of $(**)$. Similarly, we can consider the
hypercohomology of the above complex of fine sheaves, which, by a
standard result, is isomorphic to the cohomology of the
corresponding complex of global sections of the sheaves, i.e.
$$
H^*(\{\Gamma ({\overline M}, [{\text Gr}_F^*A^{\cdot}({\bf
H}_{\mathbb{C}})]_{(2)}), D''\}).
$$

One of the purpose of this paper is then to show

\vskip3mm\noindent {\bf Theorem A.} {\it ~~The complexes $(*)$ and
$(**)$ are quasi-isomorphic and hence have the same
hypercohomology, i.e.
\[
{\mathbb{H}}^*({\overline M}, \{({\bf
E}\otimes\Omega^{\cdot}_{\overline M}(\log D))_{(2)}, \theta\})
\simeq H^*(\{\Gamma ({\overline M}, [{\text Gr}_F^*A^{\cdot}({\bf
H}_{\mathbb{C}})]_{(2)}), D''\}).
\]}

\noindent {\bf Remark:}~~When $M$ is compact, Theorem A is a
formal consequence of some standard homological algebra and the
classical $\overline\partial$-Poincar\'e lemma and was stated by
Deligne, but an important point is that he used the Griffiths'
infinitesimal period relation to construct the complex $(**)$ (cf.
\cite{zuc2}, \S 1). If $M$ is quasi-compact K\"ahlerian, the proof
of the theorem becomes highly complicated and was obtained by S.
Zucker for the case of curves (cf. \cite{zuc2}, Theorem 6.3).

\vskip .3cm Let $A^k({\bf H}_{\mathbb{C}})_{(2)}$ be the sheaf of
germs of local $L^2$ measurable ${{\bf H}}_{\mathbb{C}}$-valued
$k$-forms $\phi$ on $\overline M$ for which $D\phi$ exists in the
weak sense as a local $L^2$ form. Here, $D$ is defined as follows:
Let $\phi$ be a smooth $k$-form on $\overline M$ and $v$ a smooth
section of ${\bf H}_{\mathbb{C}}$, then $D(\phi\otimes v) =
d\phi\otimes v + (-1)^k\phi\wedge\nabla v$. It is clear that
$D^2=0$. $A^k({\bf H}_{\mathbb{C}})_{(2)}$ is a fine sheaf and we
then obtain a complex $\{A^{\cdot}({\bf H}_{\mathbb{C}})_{(2)},
D\}$ of fine sheaves again. The corresponding complex of global
sections then computes its hypercohomology, and also the $L^2$- de
Rham cohomology with values in ${\bf H}_{\mathbb{C}}$ on
$\overline M$: $H^*_{(2)}({\overline M}, {\bf H}_{\mathbb{C}})$,
i.e.
$$
H^*_{(2)}({\overline M}, {\bf
H}_{\mathbb{C}})=\mathbb{H}^*({\overline M}, \{A^{\cdot}({\bf
H}_{\mathbb{C}})_{(2)}, D\})=H^*(\{\Gamma({\overline M},
A^{\cdot}({\bf H}_{\mathbb{C}})_{(2)}), D\}).
$$
The standard discussion with the K\"ahler identity and the
harmonic forms of the Laplacians of $D''$ and $D$ for the
situation of variations of Hodge structures (refer to \cite{zuc1},
\S 7) tells us that
\[
H^*(\{\Gamma ({\overline M}, [{\text Gr}_F^*A^{\cdot}({\bf
H}_{\mathbb{C}})]_{(2)}), D''\})\simeq H^*(\{\Gamma({\overline M},
A^{\cdot}({\bf H}_{\mathbb{C}})_{(2)}), D\})
\]
and hence $\simeq H^*_{(2)}({\overline M}, {\bf H}_{\mathbb{C}})$.
On the other hand, Cattani-Kaplan-Schmid's theorem \cite{cks2}
tells us that $H^*_{(2)}({\overline M}, {\bf H}_{\mathbb{C}})$ is
isomorphic to the intersection cohomology $H_{int}^*({\overline
M}, {\bf H}_{\mathbb{C}})$. Therefore, we have

\vskip 3mm \noindent {\bf Theorem B.} {\it There exists a natural
isomorphism
\[
{\mathbb{H}}^*({\overline M}, \{({\bf
E}\otimes\Omega^{\cdot}_{\overline M}(\log D))_{(2)},
\theta\})\simeq H_{int}^*({\overline M}, {\bf H}_{\mathbb{C}}).
\]}

We now give a simple application of the above theorem. More
applications will be presented in a future paper. Let still $\{M,
{\bf H}_{\mathbb{C}}={\bf H}_{\mathbb{Z}}\otimes{\mathbb{C}},
\nabla=\nabla^{1,0}+\nabla^{0,1} , {\bf F}=\{{\bf F}^p\}_{p=1}^m,
{\bf S}\}$ be a rational variation of polarized Hodge structures
with weight $m$ over $M$, $({\bf E}, \theta)$ the corresponding
Higgs bundle. Consider all endomorphisms from ${\bf
H}_{\mathbb{C}}$ to itself, denoted by ${\mathcal{E}nd}({\bf
H}_{\mathbb{C}})$, which obviously has an induced structure of
variation of Hodge structures and to which all results in \S2--5
can be applied; denote by $(\mathcal{E}nd({\bf E}),
\theta^{\mathcal{E}nd})$ the Higgs bundle corresponding to
${\mathcal{E}nd}({\bf H}_{\mathbb{C}})$. By the properties of the
Higgs field $\theta$ (cf. \S3), it can also be considered as a
morphism of sheaves
\[
\theta: {\bf T}_{\overline M}(-{\text{log}}D) \to
{\mathcal{E}nd}({\bf E}),
\]
where ${\bf T}_{\overline M}$ is the holomorphic tangent bundle of
$\overline M$. By an observation due to K. Zuo \cite{zuo}, one has
\begin{prop}(Zuo)~~
$\theta^{\mathcal{E}nd}(\theta({\bf T}_{\overline
M}(-{\text{log}}D)))=0$, i.e. the above morphism of sheaves is a
morphism of Higgs sheaves
\[
\theta: ({\bf T}_{\overline M}(-{\text{log}}D), 0) \to
((\mathcal{E}nd({\bf E}), \theta^{\mathcal{E}nd})).
\]
Here we consider ${\bf T}_{\overline M}(-{\text{log}}D)$ as a
Higgs bundle with $0$ as its Higgs field.
\end{prop}
Therefore, $\theta$ induces a morphism between the
hypercohomologies,
\begin{eqnarray*}
\theta: &&{\mathbb{H}}^*({\bf T}_{\overline M}(-\log D)
\stackrel{0}{\to}{\bf T}_{\overline M}(-\log
D)\otimes\Omega_{\overline M}^1(\log D)\stackrel{0}{\to}\cdots)
\\
&&\to{\mathbb{H}}^*({\mathcal{E}nd}({\bf
E})\stackrel{\theta}{\to}{\mathcal{E}nd}({\bf
E})\otimes\Omega^1_{\overline M}(\log
D)\stackrel{\theta}{\to}\cdots).
\end{eqnarray*}
It is worth noting that the hypercohomology on the left-hand side
is just the $\breve{\text C}$ech cohomology $H^*({\overline M},
{\bf T}_{\overline M}(-\log D))$. Thus, $\theta$ maps
$H^*({\overline M}, {\bf T}_{\overline M}(-\log D))$ into
${\mathbb{H}}^*({\mathcal{E}nd}({\bf
E})\stackrel{\theta}{\to}{\mathcal{E}nd}({\bf
E})\otimes\Omega^1_{\overline M}(\log
D)\stackrel{\theta}{\to}\cdots)$. Actually, by the properties of
$\theta$ (cf. \S3), one has a stronger restriction on the image.
\vskip .3cm \noindent {\bf Theorem C.} The image of $\theta$ lies
in the hypercohomology
$${\mathbb{H}}^*(({\mathcal{E}nd}({\bf
E}))_{(2)}\stackrel{\theta}{\to}({\mathcal{E}nd}({\bf
E})\otimes\Omega^1_{\overline M}(\log
D))_{(2)}\stackrel{\theta}{\to}\cdots).
$$
Applying Theorem B to the variation of Hodge structure
${\mathcal{E}nd}({\bf H}_{\mathbb{C}})$, one has \vskip .3cm
\noindent {\bf Theorem D.} There is a natural map, still denoted
by $\theta$,
$$\theta: H^*({\overline M}, {\bf T}_{\overline M}(-\log D))
\to H^*_{int}({\overline M}, {\mathcal{E}nd}({\bf
H}_{\mathbb{C}})).
$$

\noindent Roughly speaking, the above theorem shows that one can
transform certain geometric invariants on $\overline M$ into
topological invariants. We will give much more details and some
further applications in a future paper.

\vskip .3cm The structure of this paper is as follows. In \S2 we
will review some results about variations of polarized Hodge
structures, most of which are due to Schmid \cite{sch} and
Cattani-Kaplan-Schmid \cite{cks1} and essential for the
development of \S3 and \S4. \S3 is very technical. We first prove
an $L^2$-adapted theorem for the Hodge filtration of the
variation, which, roughly speaking, shows that the meaning of
$L^2$ on ${\bf F}^p$ and that on ${\bf E}^p$ are the same in some
sense; in the process of proving this, we construct an
$L^2$-adapted basis for ${\bf F}^p$ (and hence ${\bf E}^p$), with
the help of which we then prove the $L^2$-boundedness of $\theta$;
finally, we show the relation of $L^2$-adapted bases under
different orderings of coordinates, which is used in the proof of
the $L^2$ $\overline\partial$-Poincar\'e lemma. All the proofs
essentially depend on the nilpotent orbit theorem and the
$SL_2$-orbit theorem. It is worth pointing out that we always use
a $2$-dimensional model to discuss details; the general case is
more complicated but similar. In \S4, we will define the $L^2$
holomorphic Dolbeault complex of sheaves on $\overline M$ and
prove that the complex is actually independent of the Hodge metric
and the Poincar\'e-like metric $\omega$ of the base manifold and
essentially determined by the monodromies of the variation. The
proof heavily depends on the estimates of the Hodge metric near
the divisor. \S5 is devoted to the proof of quasi-isomorphism. In
\S4, 5, we also use the $2$-dimensional model to discuss details.

\vskip .3cm\noindent {\bf Acknowledgements:} This work began when
the second named author was visiting the Institute of Mathematical
Sciences and the Department of Mathematics at the Chinese
University of Hong Kong in the Winter of 2001; a part (\S5) was
completed when he visited the Max-Planck-Institute for Mathematics
in the Sciences during September 2002--March 2003. He thanks all
the above Institutes for hospitality and good working environment.

\section{Variation of Hodge structures and the estimate for
the Hodge norm}
In this section, we will fix some notations and review some
results, due to Schmid\cite{sch} and
Cattani-Kaplan-Schmid\cite{cks1}, which we will use later. The
reader can refer to \cite{sch, cks1} for all these.

Let $\{H^{p,q}\}$ be a Hodge structure with weight $k$ on
$H_{\mathbb{C}} = H_{\mathbb{R}}\otimes_{\mathbb{R}}{\mathbb{C}}$
and $S$ a polarization of the Hodge structure. (Call the complex
dimension of $H^{p,q}$ the {\it Hodge number}, denoted by
$h^{p,q}$.) Namely, $H_{\mathbb{C}}=\sum_{p+q=k}H^{p,q}$,
$H^{p,q}={\overline H}^{q,p}$ (so $h^{p,q}=H^{q,p}$) and $S$ is
symmetric for even $k$, skew for odd $k$, and satisfies
\begin{eqnarray*}
S(H^{p,q}, H^{r,s}) &=& 0 ~~ {\text {unless}} ~ p=s, q=r,  \\
S({\mathcal C}v, {\overline v}) &>& 0 ~~ {\text {if}} ~ v\in
H^{p,q}, v \neq 0,
\end{eqnarray*}
where ${\mathcal C}$ is the Weil operator defined by
\[
{\mathcal C}v = {\sqrt{-1}}^{p-q}v ~~{\text{for}}~ v\in H^{p,q}.
\]
To each Hodge structure $H^{p,q}$ on $H_{\mathbb{C}}$ of weight
$k$ one can assign the {\it Hodge filtration}
\[
0\subset\cdots\subset F^{p+1}\subset F^{p}\subset
F^{p-1}\subset\cdots\subset F^{o} = H_{\mathbb{C}},
\]
with the property
\[
H_{\mathbb{C}} = F^p\oplus {\overline F}^{k-p+1}, ~~ {\text
{for}}~ {\text {each}} ~ p,
\]
by setting
\[
F^p = \oplus_{i\ge p}H^{i,k-i}.
\]
Conversely, every decreasing filtration with the above property
(Hodge filtration) determines a Hodge structure $\{H^{p,q}\}$ of
weight $k$, by setting
\[
H^{p,q} = F^p\cap {\overline F}^q ~~ {\text {with}} ~ p+q=k.
\]
Therefore, weighted Hodge structures and Hodge filtrations
correspond to each other bijectively. In terms of the Hodge
filtration, the conditions on the polarization can be reformulated
as
\begin{eqnarray*}
S(F^p, F^{k-p+1}) &=& 0 ~~ {\text {for}} ~ {\text{all}} ~ p,  \\
S({\mathcal C}v, {\overline v}) &>& 0 ~~ {\text {if}} ~ v\in
H_{\mathbb{C}}, v \neq 0.
\end{eqnarray*}
Using the polarization $S$, one can define a positive Hermitian
form on $H_{\mathbb{C}}$ by $<.,.>=S({\mathcal C}.,
\overline{.})$, where ${\overline .}$ is the conjugation with
respect to the real structure $H_{\mathbb{R}}$. Call the
polarization $S$ {\it defined over} $\mathbb{Q}$ if
$H_{\mathbb{C}}$ has a lattice structure $H_{\mathbb{Z}}$, which
induces the real structure $H_{\mathbb{R}}$, and $S:
H_{\mathbb{Z}}\times H_{\mathbb{Z}}\to {\mathbb{Q}}$.

Now, assume that $H_{\mathbb{C}}$ is a fixed complex vector space
with a lattice structure $H_{\mathbb{Z}}$, which gives rise to the
real structure $H_{\mathbb{R}}$, $k$ is a fixed positive integer,
and the polarization $S$ on $H_{\mathbb{C}}$ is a nondegenerate
bilinear form of $H_{\mathbb{C}}$ defined over $\mathbb{Q}$
relative to the lattice structure $H_{\mathbb{Z}}$, which is
symmetric for even $k$ and skew for odd $k$. Consider the set of
all Hodge structures $\{H^{p,q}\}$ relative to the real structure
$H_{\mathbb{R}}$ with weight $k$, fixed Hodge numbers
$\{h^{p,q}\}$, and polarized by $S$, denoted by $\bf D$, and the
set of all decreasing filtrations $\{F^p\}_{p=0}^{k}$ of
$H_{\mathbb{C}}$ with ${\text{dim}_{\mathbb{C}}}F^p=\sum_{i\ge
k}h^{i,k-i}$ and satisfying
\[
S(F^p, F^{k-p+1}) = 0 ~~ {\text {for}} ~ {\text{all}} ~ p,
\]
denoted by $\tilde{\bf D}$. As seen in Schmid's paper\cite{sch},
$\bf D$ and $\tilde{\bf D}$ are complex homogeneous spaces: Set
\[
G_{\mathbb{C}} = \{ g\in {\text{Gl}}(H_{\mathbb{C}}) ~|~  S(gu,
gv) = S(u, v)~ {\text{for}}~ {\text{all}}~ u, v\in H_{\mathbb{C}}
\}
\]
and
\[
G_{\mathbb{R}} = \{ g\in {\text{Gl}}(H_{\mathbb{R}}) ~|~ S(gu, gv)
= S(u, v)~ {\text{for}}~ {\text{all}}~ u, v\in H_{\mathbb{R}} \},
\]
and fix a point $o\in {\bf D}$ (referred to as the reference point
or base point). Let the isotropy groups of $\tilde{\bf D}$ and
$\bf D$ be $B$ and $V$ respectively, then $\bf D$ and $\tilde{\bf
D}$ can be identified with $G_{\mathbb{R}}/V$ and
$G_{\mathbb{C}}/B$ respectively. Furthermore, $\bf D$ can be
considered as an open subset of $\tilde{\bf D}$, and $\bf D$
($\tilde{\bf D}$) has a $G_{\mathbb{R}}$
($G_{\mathbb{C}}$)-invariant holomorphic tangent subbundle,
denoted by ${\bf T}_h({\bf D})$ (${\bf T}_h({\tilde{\bf D}})$) and
referred to as the holomorphic {\it horizontal tangent
subbundles}. Let $M$ be a complex manifold. A holomorphic map
$\Psi : M\to {\bf D}$ is said to be {\it horizontal} if the image
of the tangent map of $\Psi$ lies in ${\bf T}_h({\bf D})$. This
definition also applies to $\tilde{\bf D}$. Finally, set
\[
G_{\mathbb{Z}} = \{ g\in {\text{Gl}}(H_{\mathbb{R}}) ~|~
gH_{\mathbb{Z}}=H_{\mathbb{Z}} ~{\text{and}}~ S(gu, gv) = S(u, v)~
{\text{for}}~ {\text{all}}~ u, v\in H_{\mathbb{R}} \},
\]
which is an arithmetic subgroup of $G_{\mathbb{C}}$ (cf.
\cite{sch}). Let $\mathfrak g$ and $\mathfrak{g}_0$ be the Lie
algebras of $G_{\mathbb{C}}$ and $G_{\mathbb{R}}$. It is easy to
see that $\mathfrak{g}_0$ is a real form of $\mathfrak g$.

Let $M$ be a complex manifold of dimension $n$, we consider a
family of Hodge structures parameterized by $M$ and with suitable
conditions as follows.
\begin{defi}
A (rational) variation of (polarized) Hodge structure on the base
manifold $M$ consists of the datum
\[
\{M, {\bf H}_{\mathbb{Z}}\subset{\bf H}_{\mathbb{C}}, \{{\bf
F}^p\}_{p=0}^k, \nabla =\nabla^{1,0}+\nabla^{0,1}, {\bf S}\},
\]
where 1) ${\bf H}_{\mathbb{C}}$ is a flat complex vector bundle of
the flat connection $\nabla$ containing a flat lattice ${\bf
H}_{\mathbb{Z}}$ (and hence a flat real structure ${\bf
H}_{\mathbb{R}}$: ${\bf H}_{\mathbb{C}}={\bf
H}_{\mathbb{R}}\otimes{\mathbb{C}}$); 2) $\{{\bf F}^p\}_{p=0}^k$
is a filtration of ${\bf H}_{\mathbb{C}}$, each ${\bf F}^p$ is a
holomorphic subbundle of ${\bf H}_{\mathbb{C}}$ under the
holomorphic structure $\nabla^{0,1}$, and the fibres $\{{\bf
F}^p_t\}_{p=0}^k$ form a Hodge filtration of $({\bf
H}_{\mathbb{C}})_t$ with weight $k$ with respect to the real
structure ${\bf H}_{\mathbb{R}}$ for all $t\in M$; 3) the
Griffiths' infinitesimal period relation
\[
\nabla^{1,0}{\bf F}^p\subset \Omega^1({\bf F}^{p-1})
\]
is satisfied; 4) $\bf S$ is a flat section of ${\bf
H}_{\mathbb{C}}^*\otimes{\bf H}_{\mathbb{C}}^*$ and ${\bf S}(t)$
is a polarization defined over $\mathbb{Q}$, relative to
$(H_{\mathbb{Z}})_t$, of the Hodge filtration $\{{\bf
F}^p_t\}_{p=0}^k$ for all $t\in M$. (Call $\bf S$ is a
polarization of the variation.)
\end{defi}

Similar to the situation for a single Hodge structure, one can
define naturally the Weil operator (a bundle map) $\mathfrak{C}:
{\bf H}_{\mathbb{C}}\to {\bf H}_{\mathbb{C}}$ on the bundle ${\bf
H}_{\mathbb{C}}$ such that, for each $t\in M$, $\mathfrak{C}_t$ is
the usual Weil operator related to the Hodge filtration $\{{\bf
F}^p_t\}_{p=0}^k$. Then one can define a Hermitian metric (usually
called the {\it Hodge metric}) as follows
\[
h(u, v) = {\bf S}({\mathcal C}u, {\overline v}),
\]
where $\overline .$ expresses the conjugation relative to the real
structure ${\bf H}_{\mathbb{R}}$. In the sequel, we denote the
corresponding norm by $\Vert\cdot\Vert$. Note that although $\bf
S$ and the real structure ${\bf H}_{\mathbb{R}}$ (and hence the
conjugation) are flat, the Hodge filtration $\{{\bf
F}^p\}_{p=0}^k$ (and hence the Weil operator $\mathcal C$) is not
necessarily flat. Therefore, $h$ is not necessarily flat. One of
the purpose of this section will be to describe the asymptotic
behavior of the Hodge metric $h$.

Let $\{M, {\bf H}_{\mathbb{Z}}\subset{\bf H}_{\mathbb{C}}, \{{\bf
F}^p\}_{p=0}^k, \nabla =\nabla^{1,0}+\nabla^{0,1}, {\bf S}\}$ be a
variation of polarized Hodge structure.  As usual, one can assign
to the variation a period mapping $\phi :M\to {\bf D}/\Gamma$ or a
$\rho$-equivariant mapping $\tilde{\phi}$ from the universal
covering $\tilde M$ of $M$ to $\bf D$, where $\rho$ is the induced
representation of $\pi_1(M)$ into $G_{\mathbb{Z}}$ by the flat
connection $\nabla$, $\Gamma$ is the image of $\pi_1(M)$ under the
representation $\rho$. Note that in general, ${\bf D}/\Gamma$ is
not necessarily a manifold. The following theorem is due to P.
Griffiths \cite{gri}.
\begin{thm}
The period mapping $\phi$ is holomorphic and comes from the
$\rho$-equivariant mapping $\tilde{\phi}$ which is horizontal.
\end{thm}

In practical problems, the base manifold $M$ of the variation is
usually a quasi-projective variety. By Hironaka's theorem, one can
always consider $M$ as a Zariski open subset of a smooth
projective variety $\overline M$ with $D={\overline M}\setminus M$
being a normal crossing divisor (sometimes called the {\it
singularity} of $M$). Furthermore, since the main concern in the
following is the asymptotic behavior of some objects near the
divisor $D$, one can therefore assume that $M$ is of the form
$(\triangle^*)^{n-k}\times (\triangle )^k$, $\triangle^*$ being
the punctured disk. One can actually assume that $M$ is of the
form $(\triangle^*)^{n}$, since the disk part does not affect the
asymptotic behavior. In the remaining part of this paper, we will
therefore assume that $M$ is $(\triangle^*)^{n}$ and ${\overline
M}=\triangle^n$. Obviously, $\pi_1(M)$ is generated by $n$
elements, denoted by $\sigma_1, \sigma_2, \cdots , \sigma_n$, each
$\sigma_i$ corresponding to the counter-clockwise path around $0$
of the $i$-th component of ${(\triangle}^*)^n$. It is clear that
$\pi_1(M)$ is an Abelian group. The image of $\sigma_i$ in
$\Gamma$ under $\rho$ is denoted by $\gamma_i$, which is possibly
trivial and (if nontrivial) referred to as the $i$-th {\it
Picard-Lefschez} or {\it monodromy transformation} of the
variation. The universal covering $\tilde M$ of $M$ can be
regarded as the product $U^n$, where $U$ is the upper halfplane
$\{z\in {\mathbb{C}}~|~{\text{Im}}z > 0\}$. Take the standard
coordinate systems ${\bf z}=(z_1, z_2, \cdots , z_n)$ and $(t_1,
t_2, \cdots ,t_n)$ on ${\tilde M}=U^n$ and $M=(\triangle^*)^n$
respectively. Then, the covering mapping is given by
\[
\tau : U^n\to (\triangle^*)^n, ~~ \tau (z_1, z_2, \cdots , z_n) =
(e^{2\pi\sqrt{-1}z_1}, e^{2\pi\sqrt{-1}z_2}, \cdots ,
e^{2\pi\sqrt{-1}z_n});
\]
$\sigma_i$, as a deck transformation of $U^n$, corresponds to
$z_i\to z_i+1$; and the $\rho$-equivariant condition of
$\tilde{\phi}$ can be written as
\[
{\tilde{\phi}}(z_1, \cdots ,z_i+1, \cdots ,z_n) = \gamma_i
({\tilde{\phi}}(z_1, \cdots , z_i, \cdots , z_n)),~~~
\]
for all $(z_1, \cdots , z_n)\in U^n$ and $i=1, \cdots , n$. By
Borel's lemma \cite{sch}, all $\gamma_i$ are quasi-unipotent. For
the sake of simplicity, {\it throughout this paper we always
assume that the $\gamma_i$ are unipotent, unless stated
otherwise}. Denote the logarithm of $\gamma_i$ by $N_i$, $i=1, 2
\cdots , n$. It is clear that the $\{N_i\}_{i=1}^n$ are nilpotent
and commutative. By the final remark of the Section, {\it the
nilpotency index of each linear combination of $\{N_i\}$ is not
greater than $k$}, the weight of the variation, i.e.,
$(N_i)^{k+1}=0$. Define a map by
\[
\tilde{\psi}({\bf z}) = \exp (-\sum_{i=1}^nz_iN_i)\cdot \tilde{\phi}({\bf z}).
\]
It is easy to check that $\tilde{\psi}$ remains invariant under
the translations $z_i\to z_i+1, 1\le i\le n$. So, $\tilde{\psi}$
drops to a mapping
\[
\psi : ({\triangle}^*)^n\to {\tilde{\bf D}}.
\]
\begin{thm}({\bf Nilpotent Orbit Theorem}) The map $\psi$ extends
holomorphically to $(\triangle)^n$. Denote the point ${\psi}(0, 0,
\cdots , 0)\in{\tilde{\bf D}}$ by $F$, then the map
\[
{\bf z}(\in {\mathbb{C}}^n)\to{\exp}(\sum_{i=1}^nz_iN_i)\cdot F\in
{\tilde{\bf D}}
\]
is horizontal (equivalently, $N_iF^p\subset F^{p-1}$). Finally,
there exist constants $\alpha , \beta, K\ge 0$, such that under
the restrictions ${\text{Im}}z_i\ge\alpha,~~1\le i\le n$, the
point $\exp(\sum_{i=1}^nz_iN_i)\cdot F$ lies in $\bf D$ and
satisfies the following inequality
\[
{\text d}(\exp (\sum_{i=1}^nz_iN_i)\cdot F, {\tilde\phi}({\bf z}))
\le K\sum_{i=i}^n({\text{Im}}(z_i))^{\beta}\exp(-2\pi{\text{Im}}(z_i)),
\]
where $\text d$ again denotes a $G_{\mathbb{R}}$-invariant
Riemannian distance function on $\bf D$; moreover, the constants
$\alpha , \beta , K$ depend only on the choice of $\text d$ and
the weight and Hodge numbers associated to $\bf D$.
\end{thm}

\begin{defi}
A {\bf nilpotent orbit} is a mapping $\theta: {\mathbb{C}}^n\to
{\tilde{\bf D}}$ of form
\[
\theta ({\bf z}) = \exp (\sum_{i=1}^{n}z_iN_i)\cdot F
\]
where 1) $F=\{F^p\}_{p=0}^k\in {\tilde{\bf D}}$; 2)
$\{N_i\}_{i=1}^n$ is a commutative set of nilpotent elements of
$\mathfrak{g}_0$, horizontal at $F$, i.e., $N_iF^p\subset F^{p-1}$
(hence $\theta$ is a horizontal mapping); 3) There exists an
$\alpha\in {\mathbb{R}}$, such that $\theta ({\bf z})\in {\bf D}$
for ${\text{Im}}(z_i)>\alpha$.
\end{defi}

\noindent {\bf Remark.} By the definition, the map
${\exp}(\sum_{i=1}^nz_iN_i)\cdot F$ in the above theorem is a
nilpotent orbit.

Concerning nilpotent orbits, one has the {\it $SL_2$-orbit
theorem}. We now recall it. We first introduce a notion, the
so-called (polarized) {\it mixed Hodge structure} \cite{d}, which
is very important in Hodge theory. As before, we still let
$H_{\mathbb{C}}$ be a complex vector space with a lattice
structure $H_{\mathbb{Z}}$, $k$ a positive integer, and $S$ a
non-degenerate bilinear form on $H_{\mathbb{C}}$, defined over
$\mathbb{Q}$ and such that $S(u, v)=(-1)^kS(v, u)$ for $u, v\in
H_{\mathbb{C}}$. We also have corresponding other notations, e.g.
$\tilde{\bf D}$, $\bf D$, $G_{\mathbb{C}}$, $\mathfrak g$,
$G_{\mathbb{R}}$, $\mathfrak{g}_0$ etc..
\begin{defi}
A mixed Hodge structure defined over $\mathbb{R}$ on
$H_{\mathbb{C}}$ consists of a pair of finite filtrations of
$H_{\mathbb{C}}$
\begin{eqnarray*}
&&W:~~\cdots\subset W_l\subset W_{l+1}\subset\cdots ~~~
({\text{the}}~{\text{weight}}~{\text{filtration}})  \\
&&F(\in{\tilde D}):~~\cdots\subset F^p\subset F^{p-1}\subset\cdots
~~~({\text{the}}~{\text{decreasing}}~{\text{filtration}})
\end{eqnarray*}
such that i) $W$ is defined over $\mathbb{R}$, ii) the filtration
$F({\text{Gr}}_l(W_*))$ induced on
${\text{Gr}}_l(W_*)=W_l/W_{l-1}$ by $F$ is a Hodge filtration of
weight $l$.
\end{defi}
\noindent {\bf Remark.} This notion, like the notion of Hodge
structure, is compatible with the natural operations, e.g. duality
and tensor product etc..

A {\it splitting} of a mixed Hodge structure $(W, F)$ is a
bigrading $H_{\mathbb{C}}=\oplus J^{p,q}$ such that
\[
W_l = \oplus_{p+q\le l}J^{p,q}, ~~F^p = \oplus_{r\ge p}J^{r,s}.
\]
An $(r, r)$-{\it morphism} $(W, F)$ is an element $X$ of
$\mathfrak{gl}(H_{\mathbb{C}})$ satisfying $X(W_l)\subset
W_{l+2r}$ and $X(F^p)\subset F^{p+r}$. So the nilpotent element
$N$ of a polarized mixed Hodge structure $(W, F, N)$ (see the
following definition) is a $(-1,-1)$-morphism of the mixed Hodge
structure $(W, F)$. Call an $(r, r)$-morphism $X$ of $(W, F)$ {\it
compatible} with the splitting $\{J^{p,q}\}$ if $X(J^{p,q})\subset
J^{p+r,q+r}$. A construction of Deligne \cite{d, cks1} tells us
that a mixed Hodge structure always admits some splittings which
are compatible with all its morphisms. A mixed Hodge structure
$(W, F)$ is said to split over $\mathbb{R}$ if it admits a
splitting $\{J^{p,q}\}$ such that $J^{p,q}=\overline{J^{q,p}}$,
here $\overline{.}$ represents the conjugation relative to the
real structure. Such a splitting is then called a {\it real}
splitting of $(W, F)$. The following proposition will be important
for stating the $SL_2$-orbit theorem for several variables.
\begin{prop}
Given a mixed Hodge structure $(W, F)$ on
$H_{\mathbb{C}}=H_{\mathbb{R}}\otimes_{\mathbb{R}}{\mathbb{C}}$,
there exists a unique $\delta\in{\mathfrak{gl}}(H_{\mathbb{R}})$
such that $(W, \exp(-{\sqrt{-1}}\delta)F)$ is a mixed Hodge
structure which splits over $\mathbb{R}$. Furthermore, every
morphism of $(W, F)$ commutes with $\delta$.
\end{prop}

The following definition is also important for our exposition.
\begin{defi}
A polarized mixed Hodge structure consists of a mixed Hodge
structure $(W, F)$ defined over $\mathbb{R}$ on $H_{\mathbb{C}}$
and a nilpotent element $N\in {\mathfrak g}_0$ such that 1)
$N^{k+1}=0$; 2) $W=W(N)[-k]$, where $W(N)$ is the weight
filtration induced canonically by $N$ and
$(W(N)[-k])_l=(W(N))_{l-k}$; 3) $S(F^p, F^{k-p+1})=0$; 4)
$NF^p\subset F^{p-1}$; 5) The Hodge structure (of weight $s+l$) on
the primitive part $P_{s+l}={\text{Ker}}(N^{l+1}:
{\text{Gr}}_{s+l}(W_*)\to{\text{Gr}}_{s+l}(W_*))$ is polarized by
the form $S(., N^l.)$.
\end{defi}
Then, for polarized mixed Hodge structures, one has the following
proposition, which is a consequence of the $SL_2$-orbit theorem in
a single variable \cite{sch}.
\begin{prop}
Let $\theta (z)=\exp(zN)\cdot F$ be a nilpotent orbit for
$z\in{\mathbb{C}}$, a nilpotent element $N\in{\mathfrak g}_0$ with
nilpotent index $k$, and a filtration $F\in {\tilde{\bf D}}$,
i.e., satisfying that the mapping $\theta$ is horizontal and for
some $\alpha >0$ and ${\text{Im}}z>\alpha$, $\theta (z)\in{\bf
D}$. Then $(W(N)[-k], F, N)$ is a polarized mixed Hodge structure.
\end{prop}

Using the above notions and results about (polarized) mixed Hodge
structures, one can reformulate the $SL_2$-orbit theorem in a
single variable in the following manner, which will be helpful for
the reformulation of $SL_2$-orbit theorem in several variables.
From the proposition above, $(W=W(N)[-k], F)$ is a mixed Hodge
structure polarized by $N$ (so, $N$ is a $(-1, -1)$-morphism of
$(W, F)$); on the other hand, by Proposition 1, there exists a
$\delta\in{\mathfrak{gl}}(H_{\mathbb{R}})$ such that $(W, {\tilde
F}={\text{exp}}(-{\sqrt{-1}}\delta)F)$ is a mixed Hodge structure
splitting over $\mathbb{R}$. Let $\{{\tilde J}^{p,q}\}$ be a real
splitting of $(W, {\tilde F})$, i.e.,
\[
{\tilde F}^p = \oplus_{r\ge p}{\tilde J}^{r,s},~ W_l =
\oplus_{r+s\le l}{\tilde J}^{r,s}, ~ {\tilde J}^{r,s} =
{\overline{{\tilde J}^{s,r}}}.
\]
Define a semisimple transformation ${\tilde Y}$ as follows,
\[
{\tilde Y}v = (p+q-k)v  ~~{\text{for}}~v\in {\tilde J}^{p,q}.
\]
Then, some additional discussion (see $\S 3$ of \cite{cks1}) shows
that $\{{\tilde Y}, N\}$ can be completed to an
$\mathfrak{sl}_2$-triple $\{{\tilde N}^+, {\tilde Y}, N\}$:
Actually, the proposition 1 tells us that $N$ is a
$(-1,-1)$-morphism of $(W, {\tilde F})$ and compatible with the
above splitting $\{{\tilde J}^{p,q}\}$ (this needs some additional
argument), while $\tilde Y$ is a $(0,0)$-morphism; and hence
${\tilde N}^+$ is a $(1,1)$-morphism of $(W, {\tilde F})$. Like
$N$, the transformations $\delta , {\tilde Y}, {\tilde N}^+$ are
infinitesimal isometries of the polarization $S$ and belong to
${\mathfrak g}_0$. Thus, we have constructed a representation
${\tilde{\rho}}_*$ of $\mathfrak{sl}_2({\mathbb{C}})$ into
${\mathfrak g}$ with
\begin{eqnarray*}
&&{\tilde{\rho}}_*(Z) = {\sqrt{-1}}({\tilde N}^+ - N), \\
&&{\tilde{\rho}}_*(X_+) = {\frac 1 2}({\tilde N}^+ + N + {\sqrt{-1}}{\tilde Y}),   \\
&&{\tilde{\rho}}_*(X_-) = {\frac 1 2}({\tilde N}^+ + N - {\sqrt{-1}}{\tilde Y})
\end{eqnarray*}
(or ${\tilde{\rho}}_*({\bf y})={\tilde Y}, {\tilde{\rho}}_*({\bf
n}_+)={\tilde N}^+, {\tilde{\rho}}_*({\bf n}_-)=N$), which lifts
to a homomorphism
\[
{\tilde{\rho}}: SL_2({\mathbb{C}})\to G_{\mathbb{C}}.
\]
It is clear that ${\tilde{\rho}}$ is defined over $\mathbb{R}$.
Again since ${\tilde N}^+$ and ${\tilde Y}$ are $(1,1)$ and
$(0,0)$-morphisms of $(W, {\tilde F})$ respectively, so they fix
the filtration $\tilde F$. Let $B$ be the isotropy subgroup in
$G_{\mathbb{C}}$ of $\tilde{\bf D}$ fixing ${\tilde F}$, the Lie
algebra $\mathfrak b$ of which clearly contains ${\tilde N}^+$ and
${\tilde Y}$; let $L$ be the isotropy subgroup in
$SL_2({\mathbb{C}})$ of ${\bf P}^1$ fixing $0$, the Lie algebra
$\mathfrak l$ of which clearly contains ${\bf n}_+$ and $\bf y$.
Thus, the homomorphism $\tilde{\rho}$ induces an equivariant
embedding
\[
{\bf P}^1\simeq SL_2({\bf{C}})/L\to {\tilde{\bf D}}
\]
by $g\to {\tilde{\rho}}(g){\tilde F}$, which can be written as,
for $z\in {\mathbb{C}}$,
\[
z\to {\text{exp}}(zN)\cdot{\tilde F}.
\]
A further argument shows that $\exp(zN)\cdot{\tilde F}\in{\bf D}$
as ${\text{Im}}z>0$ (in particular, $e^{\sqrt{-1}N}{\tilde
F}\in{\bf D}$), namely the upper half plane is mapped into $\bf D$
\cite{cks1}.

Write $\exp(zN)\cdot{F} = {\tilde g}(z)\exp(zN)\cdot{\tilde F}$.
Since for some $\alpha >0$, as ${\text{Im}}z>0$,
$\exp(zN)\cdot{F}\in{\bf D}$, so ${\tilde g}(z)\in G_{\mathbb{R}}$
for ${\text{Im}}z>0$. In particular, ${\tilde g}({\sqrt{-1}}y)\in
G_{\mathbb{R}}$. In the following, write ${\tilde g}(y)$ instead
of ${\tilde g}({\sqrt{-1}}y)$. Then, the $SL_2$-orbit theorem in a
single variable applies to ${\tilde g}$ here. Especially, ${\tilde
g}(\infty )\in G_{\mathbb{R}}$. Setting $\mu ={\text{log}}{\tilde
g}(\infty )$ and ${\tilde F}_0={\tilde g}(\infty ){\tilde F}$, we
obtain a polarized mixed Hodge structure split over $\mathbb{R}$
$(W, {\tilde F}_0, N)$ \cite{cks1}. The most important is that
\[
(W, {\tilde F}_0) ~{\text{is}}~ {\text{canonically}}~
{\text{attached}}~ {\text{to}}~ (W, F),
\]
although defined in terms of the nilpotent orbit
$\exp(zN)\cdot{F}$, and ${\tilde g}(\infty )$ depends on $\delta$
(equivalently, on the mixed Hodge structure $(W, F)$, but not on
$N$). This point will be very useful in  the following
reformulation for the $SL_2$-orbit theorem for several variables.

We can now begin with reformulating the $SL_2$-orbit theorem for
several variables. Let $F\in\tilde{\bf D}$, $\alpha$ a fixed
positive constant, and $N_1, N_2, \cdots , N_n\in\mathfrak{g}_0$
some commuting nilpotent elements of nilpotent indices $k$, such
that for ${\bf z}=(z_1, \cdots , z_n)\in{\mathbb{C}}^n$
\begin{eqnarray*}
&&{\text{1)}}~~{\text{the}}~{\text{mapping}}~\theta({\bf z}) =
{\text{exp}}(\sum_{i=1}^nz_iN_i)\cdot F
~{\text{is}}~{\text{horizontal}};  \\
&&{\text{2)}}~~{\text{exp}}(\sum_{i=1}^nz_iN_i)\cdot F\in{\bf
D}~{\text{for}}~{\text{Im}}z_1>\alpha, 1\le i\le n,
\end{eqnarray*}
i.e., $\theta({\bf z}) = {\text{exp}}(\sum_{i=1}^nz_iN_i)\cdot F$
is a nilpotent orbit. In the remaining part of this section and
the next sections, we will have to fix an ordering of variables
$z_1, z_2, \cdots , z_n$, and correspondingly the commuting
nilpotent elements $N_1, N_2, \cdots , N_n$. First of all, we
state a purely algebraic result, which is due to Cattani and
Kaplan \cite{ck}, as follows:
\begin{lem}
Let $N_1, N_2, \cdots ,N_m$ be some commutative nilpotent linear
transformations in $\mathfrak{g}_0$ with nilpotent indices $k$.
Set $C(N_1, N_2, \cdots , N_m) = \{\sum_{j=1}^m\lambda_jN_j~|  \\
~\lambda_j\in {\mathbb{R}}, \lambda > 0\}$. Then the weight
filtration on $H_{\mathbb{C}}$ canonically attached to each
nilpotent element $N\in C(N_1, N_2, \cdots , N_m)$ is the same.
\end{lem}
Thus, for our present situation, $C_{\bf r}=C(N_1, \cdots ,
N_r)=\{\sum_{j=1}^r \lambda_jN_j~|~\lambda_j\in {\mathbb{R}},
\lambda
> 0\}$, $1\le r\le n$, determines a unique weight filtration
defined over $\mathbb{R}$ on $H_{\mathbb{C}}$, denoted by $W^{\bf
r}$, as follows
\[
0\subset W_{-k}^{\bf r}\subset W^{\bf
r}_{-k+1}\subset\cdots\subset W_{k-1}^{\bf r}\subset W_{k}^{\bf r}
= H_{\mathbb{C}}.
\]
So, by the proposition 2, we have that $(W^{\bf n}[-k], F)$ is a
mixed Hodge structure polarized by each $N\in C_{\bf n}$. The
previous discussion then shows that one can attach canonically to
$(W^{\bf n}[-k], F)$ another mixed Hodge structure $(W^{\bf
n}[-k], {\tilde F}_{\bf n})$ also polarized by each $N\in C_{\bf
n}$ and splitting over $\mathbb{R}$ for some ${\tilde F}_{\bf
n}\in {\tilde{\bf D}}$. Since $(W^{\bf n}[-k], {\tilde F}_{\bf
n})$ is polarized by every $N\in C_{\bf n}$, the mapping
\[
(z_1, z_2, \cdots , z_{n-1})\to
{\text{exp}}(\sum_{i=1}^{n-1}z_iN_i)\cdot(e^{{\sqrt{-1}}N_n}{\tilde
F}_{\bf n})
\]
is a nilpotent orbit, namely satisfies the properties mentioned
above. For this nilpotent orbit, we can do the same argument as
before:
$(W^{\bf{n-1}}[-k], \\
e^{{\sqrt{-1}}N_n}{\tilde F}_{\bf n})$ is a mixed Hodge structure
polarized by each $N\in C_{\bf{n-1}}$ and one can attach
canonically to it another $\mathbb{R}$-split mixed Hodge structure
$(W^{\bf{n-1}}[-k], {\tilde F}_{\bf{n-1}})$ polarized by each
$N\in C_{\bf{n-1}}$ for some ${\tilde F}_{\bf{n-1}}\in{\tilde{\bf
D}}$. {\it Inductively, one has that, for $1\le r\le n$, there
exists some ${\tilde F}_{\bf{r}}\in{\tilde{\bf D}}$ such that
$(W^{\bf{r}}[-k], {\tilde F}_{\bf{r}})$ is a $\mathbb{R}$-split
mixed Hodge structure polarized by each $N\in C_{\bf r}$ and
canonically attached to $(W^{\bf{r}}[-k],
e^{{\sqrt{-1}}N_{r+1}}{\tilde F}_{\bf{r+1}})$; and
\[
(z_1, z_2, \cdots , z_{r})\to
{\text{exp}}(\sum_{i=1}^{r}z_iN_i)\cdot(e^{{\sqrt{-1}}N_{r+1}}{\tilde
F}_{\bf{r+1}})
\]
is a nilpotent orbit.
}

Letting $\{J^{p,q}_{\bf r}\}$ be the real splitting of
$(W^{\bf{r}}[-k], {\tilde F}_{\bf{r}})$, for $1\le r\le n$, we
define some semisimple transformations ${\tilde Y}_{\bf r}$
\[
{\tilde Y}_{\bf r}v = (p+q-k)v,~~{\text{for}} ~ v\in J^{p,q}_{\bf r}.
\]
One can show that $\{{\tilde Y}_{\bf r}\}_{r=1}^n$ is a commuting
set of semisimple endomorphisms in $\mathfrak{g}_0$. We also
define ${\tilde N}^-_{r}$ as the component of $N_r$ in the
subspace $\cap_{j=1}^{r-1}{\text{ker}}({\text{ad}}{\tilde Y}_{\bf
j})$ relative to the decomposition of $\mathfrak{g}_0$ in
eigenspaces of the commuting set of semisimple endomorphisms
$\{{\text{ad}}({\tilde Y}_{\bf j})\}_{j=1}^{r-1}$. Set ${\tilde
N}^-_{\bf r}=\sum_{j=1}^r{\tilde N}^-_{j}$. Then, one can show
that ${\tilde Y}_{\bf r}, {\tilde N}^-_{\bf r}$ can be expanded to
an $\mathfrak{sl}_2$-triple $({\tilde Y}_{\bf r}, {\tilde
N}^-_{\bf r}, {\tilde N}^+_{\bf r})$ and $\{({\tilde Y}_{\bf r},
{\tilde N}^-_{\bf r}, {\tilde N}^+_{\bf r})\}_{r=1}^n$ induces a
representation $\rho_*$ of $({\mathfrak{sl}}_2({\mathbb{C}}))^n$
into $\mathfrak{g}$, which has a particular property. This can be
explained as follows. The Lie algebra
$\mathfrak{sl}_2({\mathbb{C}})$ admits a natural Hodge structure
of weight $0$, with
\begin{eqnarray*}
&&(\mathfrak{sl}_2({\mathbb{C}}))^{-1,1}={\overline{(\mathfrak{sl}_2({\bf
C}))^{1,-1}}}
={\mathbb{C}}({\sqrt{-1}}{\bf y}+{\bf n}_++{\bf n}_-)          \\
&&(\mathfrak{sl}_2({\mathbb{C}}))^{0,0}={\mathbb{C}}({\bf
n}_+-{\bf n}_-),
\end{eqnarray*}
which gives rise to a Hodge structure of weight $0$ on
$(\mathfrak{sl}_2({\mathbb{C}}))^{n}$. Here,
\[
{\bf y}=\left(\begin{matrix} -1&0\\0&1\end{matrix}\right), {\bf
n}_+=\left(\begin{matrix} 0&0\\1&0\end{matrix}\right), {\bf
n}_-=\left(\begin{matrix} 0&1\\0&0\end{matrix}\right).
\]
Fix a point $o\in{\bf D}$, called the {\it reference point} or
{\it base point}, the corresponding Hodge structure of which is
denoted by $\{H_0^{p,q}\}$. One can then define a Hodge structure
${\mathfrak{g}^{p,-p}}$ of weight $0$ on $\mathfrak{g}$:
${\mathfrak{g}^{p,-p}}=\{X\in \mathfrak{g}~|~XH^{r,s}_0\subset
H^{r+p,s-p}_0\}$. We say that a Lie algebra homomorphism $\xi:
(\mathfrak{sl}_2({\mathbb{C}}))^{n}\to \mathfrak{g}$ is {\it
Hodge} at $o$ if $\xi$ is a $(0,0)$-morphism from the Hodge
structure on $(\mathfrak{sl}_2({\mathbb{C}}))^{n}$ to the induced
Hodge structure above on $\mathfrak g$, i.e.,
$\xi(((\mathfrak{sl}_2({\bf
C}))^n)^{r,-r})\subset{\mathfrak{g}^{r,-r}}$. Then the following
theorem shows that $\rho_*$ above is Hodge. We now are in the
position of stating the $SL_2$-orbit theorem for several
variables.

\begin{thm}({\bf $SL_2$-orbit theorem for several variables})
For the above fixed ordering of the variables ${\bf z}=(z_1, z_2,
\cdots ,z_n)$ in ${\mathbb{C}}^n$ and the nilpotent orbit
$\theta({\bf z})$, there exists a unique Lie group homomorphism
\[
\rho: (SL_2({\mathbb{C}}))^n\to G_{\mathbb{C}}
\]
with the following properties:

i) The induced homomorphism $\rho_*$ is Hodge at the point
$\exp({\sqrt{-1}}N_1)\cdot{\tilde{F_{\bf 1}}}\in {\bf D}$ and
$\tilde{F_{\bf r}}=\exp(-{\sqrt{-1}}{\tilde N}_{\bf
r}^-)(\exp({\sqrt{-1}}N_1)\cdot{\tilde{F_{\bf 1}}})$.

ii) The image of the $r$-th ${\bf n}_-$ (${\bf y}, {\bf n}_+$)
factor of $(\mathfrak{sl}_2({\mathbb{C}}))^n$ under $\rho_*$ is
${\tilde N}^-_r$ (${\tilde Y}_r, {\tilde N}^+_r$ respectively) and
$W({\tilde N}^-_{\bf r})=W^{\bf r}$.

\noindent In addition, setting $z_i=x_i+{\sqrt{-1}}y_i$, there
exist $G_{\mathbb{R}}$-valued functions $g_{\bf r}(y_1, \cdots,
y_r)$ defined for $y_i>0$ if $1\le r\le n-1$ and for $y_i>\alpha$,
$\alpha\in{\mathbb{R}}$, if $r=n$, such that

iii) For $j<r\le n$, $g_{\bf j}(y_1, \cdots, y_j)$ commutes with
${\tilde Y}_{\bf r}$ and leaves the point ${\tilde F}_{\bf r}$
fixed. In particular if $j<r$, $g_{\bf j}(y_1, \cdots, y_j)$ is a
$(0,0)$-morphism of the mixed Hodge structure $(W^{\bf r},
\tilde{F}_{\bf r})$.

iv) $\sum_{s=1}^ry_sN_s =
{\text{Ad}}(\Pi_{j=r-1}^1g_{\bf j}
({\frac{y_1}{y_{j+1}}}, \cdots, {\frac{y_j}{y_{j+1}}}))
\sum_{s=1}^ry_s{\tilde N}^-_s$.

v) $\exp\sqrt{-1}(\sum_{j=1}^ry_jN_j)
\cdot(\exp(\sqrt{-1}N_{r+1}){\tilde F}_{{\bf r}+1}) =  \\
(\Pi_{j=r}^1g_{\bf j}({\frac{y_1}{y_{j+1}}}, \cdots, {\frac{y_j}{y_{j+1}}}))\cdot
((\exp\sqrt{-1}\sum_{j=1}^ry_j{\tilde N}^-_j)
{\tilde F}_{\bf r})$, $1\le r\le n$ and $y_{r+1}=1$.

vi) The functions $g_{\bf r}(y_1, \cdots, y_r)$ and $g_{\bf
r}(y_1, \cdots, y_r)^{-1}$ have power series expansions in
nonpositive powers of ${\frac{y_1}{y_2}}, {\frac{y_2}{y_3}},
\cdots , {\frac{y_{r-1}}{y_r}}, y_r$ with constant term $1$ and
convergent in any region of the form ${\frac{y_1}{y_2}}>\beta,
\cdots, {\frac{y_{r-1}}{y_r}}>\beta, y_r>\beta$, $\beta >0$.
\end{thm}

Using the nilpotent orbit theorem and the $SL_2$-orbit theorem,
one can describe the asymptotic behavior of the Hodge metric. Let
$\{M, {\bf H}_{\mathbb{Z}}\subset{\bf H}_{\mathbb{C}}, \{{\bf
F}^p\}_{p=0}^k, \nabla =\nabla^{1,0}+\nabla^{0,1}, {\bf S}\}$ be a
variation of Hodge structures, where $M$ is $(\triangle^*)^n$. Let
$h$ be the corresponding Hodge metric of the variation. We still
use the previous notations. Fix a point $s\in (\triangle^*)^n$.
$N_i$ can be considered as a nilpotent linear transformation of
the complex vector space $({\bf H}_{\mathbb{C}})_s$. So, by Lemma
1, for $C(N_1, N_2, \cdots , N_j)$, $1\le j\le n$, there exists a
unique weight filtration of $({\bf H}_{\mathbb{C}})_s$, denoted by
$W^{\bf j}$, with some standard properties (see Section 6 of
\cite{sch}). Namely,
\[
0\subset W_{-k}^{\bf j}\subset W_{-k+1}^{\bf j}\subset\cdots
\subset W_{k-1}^{\bf j}\subset W_k^{\bf j} = ({\bf
H}_{\mathbb{C}})_s
\]
such that, for all $N\in C(N_1, N_2, \cdots , N_j)$, $N(W_l^{\bf
j})\subset W_{l-2}^{\bf j}$ and, for $l\ge 0$
\[
N^l: {\text Gr}_l(W_*^{\bf j})\to{\text Gr}_{-l}(W_*^{\bf j})
\]
is an isomorphism, where ${\text{Gr}}_l^{W_*^{\bf j}}=W_l^{\bf
j}/W_{l-1}^{\bf j}$. This filtration is called {\it monodromy
weight filtration}, according to Deligne. On the other hand,
relative to the flat structure of the bundle ${\bf
H}_{\mathbb{C}}\to (\triangle^*)^n$, the Picard-Lefschetz
transformations $\gamma_i$, and hence the logarithms $N_i$, can be
regarded as flat sections of the bundle ${\bf
H}^*_{\mathbb{C}}\otimes {\bf H}_{\mathbb{C}}$, so one can
construct, from the above filtrations, the corresponding
filtrations of ${\bf H}_{\mathbb{C}}$ by some locally constant
sheaves (i.e., some sheaves of locally flat sections), denoted by
$\{{\bf W}_*^{\bf j}\}$, which are $\{\gamma_1, \cdots ,
\gamma_j\}$-invariant. Note that since ${\bf H}_{\mathbb{C}}$ has
a flat lattice ${\bf H}_{\mathbb{Z}}\subset{\bf H}_{\mathbb{C}}$,
so the above filtrations of locally constant sheaves are defined
over $\mathbb{Q}$; in general, if ${\bf H}_{\mathbb{C}}$ is
defined only over $\mathbb{R}$, the filtrations in question are
defined only over $\mathbb{R}$.

Now, we consider a canonical sheaf extension $\overline{\bf
H}_{\mathbb{C}}$ of ${\bf H}_{\mathbb{C}}$ across the singularity,
when considering ${\bf H}_{\mathbb{C}}$ as a holomorphic vector
bundle under the holomorphic structure ${\nabla}^{0,1}$, as
follows. As a sheaf, the germs of the sections of $\overline{\bf
H}_{\mathbb{C}}$ at the singularity are generated by the elements
of $\Gamma (M, {\bf H_{\mathbb{C}}})$ which are of the form
\[
{\tilde v} = {\text{exp}}({\frac{1}{2\pi\sqrt{-1}}}\sum_{i=1}^nN_i\log t_i)v,
\]
where $v$ is a multivalued parallel section of ${\bf
H}_{\mathbb{C}}$. It is clear that $\tilde v$ is holomorphic.
Similarly, we have the canonical extensions $\{\overline{\bf
W}_l^{\bf j}\}_{l=-k}^k$ of $\{{\bf W}_l^{\bf j}\}_{l=-k}^k$,
which induce a filtration $\overline{\bf H}_{\mathbb{C}}$, for
$1\le j\le n$. We can now state the asymptotic behavior near the
singularity of the Hodge metric \cite{cks1} as follows
\begin{thm}
All the notations as above. If $v\in \cap_j{\bf W}_{l_j}^{\bf j}$
is a flat section and the projection of $v$ in each
${\text{Gr}}_{l_j}^{{{\bf W}^{\bf j}_*}}$ is nontrivial, then the
Hodge norm of $\tilde v$ satisfies
\[
\Vert{\tilde v}\Vert \sim ({\frac{\log|t_1|}{\log|t_2|}})^{\frac{l_1}{2}}
({\frac{\log|t_2|}{\log|t_3|}})^{\frac{l_2}{2}}\cdots (-\log|t_n|)^{\frac{l_n}{2}}
\]
on any region of the form
\[
D_{\epsilon} = \{(t_1, t_2, \cdots , t_n)\in
(\triangle^*)^n~|~{\frac{\log|t_1|}{\log|t_2|}}>\epsilon,
{\frac{\log|t_2|}{\log|t_3|}}>\epsilon, \cdots , -\log|t_n|>\epsilon \}
\]
for any $\epsilon >0$. The same estimate holds for $\Vert{v}\Vert$
when restricted to any sector $|{\arg} t_j|<\delta$, $j=1, 2,
\cdots , n$ for any $\delta > 0$.
\end{thm}

To conclude this section, we establish some notation, that will be
used in the next two sections. Let $\bf H$ be a Hermitian vector
bundle over a Riemannian manifold $M$. Let ${\bf v}=\{v_1, v_2,
\cdots , v_q\}$ be a global frame field of $\bf H$. Then $\bf v$
is said to be {\it $L^2$-adapted} if that $\sum_{i+1}^qf_iv_i$ is
square integrable implies that each $f_iv_i$ is square integrable
with respect to the Hermitian and Riemannian metrics, where the
$f_i$ are smooth functions on $M$. Obviously, the
$L^2$-adaptedness condition is invariant under constant matrix
transformations. Similarly, it is invariant under scaling, so one
can normalize the frame without changing the $L^2$-adaptedness of
$\bf v$. The following lemma will be useful in the sequel (for the
proof, see \cite{zuc2}, Lemma 4.5).
\begin{lem}
Let ${\bf v}$ be a frame for $\bf H$ with $\sup\Vert
v_i\Vert<\infty$ for all $i$. A sufficient condition that $\bf v$
be $L^2$-adapted is that the matrix of the inner products $(<v_i,
v_j>)$ has a bounded inverse.
\end{lem}

\section{$L^2$-adaptedness of the filtration on
${\Omega}^.({\bf H}_{\mathbb{C}})_{(2)}$}

Let $(M, {\bf H}_{\mathbb{C}}={\bf
H}_{\mathbb{Z}}\otimes{\mathbb{C}}, \nabla
=\nabla^{1,0}+\nabla^{0,1}, {\bf F}=\{{\bf F}^p\}, {\bf S})$ be a
rational variation of polarized Hodge structure with weight $k$
defined over $\mathbb{Z}$ on $M$. For the sake of simplicity, we
always assume that $M$ is complex $2$-dimensional, unless stated
otherwise, because the case of higher dimension is similar. Thus,
as in \S2, $M$ can be seen as $({\triangle^*})^2$ and ${\overline
M}=\triangle^2$. Let $\omega$ be the Poincar\'e-like metric on
$M$, which is quasi-isometric to, near the singularity,
\[
\eta = {\frac{\sqrt{-1}}{2}}[{\frac{dt_1\wedge d{\overline
t}_1}{|t_1|^2(-{\text{log}}|t_1|)^2}} + {\frac{dt_2\wedge
d{\overline t}_2}{|t_2|^2(-{\text{log}}|t_2|)^2}}],
\]
and of finite volume, where $(t_1, t_2)$, as in \S2, are the
standard coordinates of $(\triangle^*)^2$. Let $N_1$ and $N_2$ be
the logarithmic monodromies of the variation, which are nilpotent.
We fix an ordering $(N_1, N_2)$ from now on as in \S2. (Here, we
again assume that the monodromies $\gamma_1$ and $\gamma_2$ are
unipotent as in \S 2.)

Using the Poincar\'e-like metric $\omega$ on $M$ and the Hodge
metric $h$ on ${\bf H}_{\mathbb{C}}$, we define ${\Omega}^r({\bf
H}_{\mathbb{C}})_{(2)}$ to be the sheaf of germs of local $L^2$
holomorphic $r$-forms valued in $j_*{\bf H}_{\mathbb{C}}$ on
$\overline M$, where $j$ is the inclusion map of $M$ into
$\overline M$. In the next section, we will study in detail the
asymptotic behavior near the divisor of its sections and show that
it is independent of both metrics, but only depends on the
logarithmic monodromies of the variation. According to the Hodge
filtration $\{{\bf F}^p\}$ of the variation, one can then
construct a filtration of $\Omega^r({\bf H}_{\mathbb{C}})_{(2)}$,
denoted by $F^{p}\Omega^r({\bf H}_{\mathbb{C}})_{(2)}$, which is
the sheaf of germs of local $L^2$-holomorphic $r$-forms on
$\overline M$ with values in $j_*{\bf F}^{p-r}$. Then, using the
projection of ${\bf F}^{p-r}$ to ${\bf E}^{p-r}={\bf F}^{p-r}/{\bf
F}^{p-r+1}$, one can also do a projection from $F^{p}\Omega^r({\bf
H}_{\mathbb{C}})_{(2)}$ to $(\Omega^r\otimes {\bf
E}^{p-r})_{(2)}$, the sheaves of germs of local $L^2$-holomorphic
forms on $\overline M$ with values in $j_*{\bf E}^{p-r}$, which,
as seen in the next section, is actually $(j_*{\bf
E}^{p-r}\otimes\Omega^r_{\overline M}(\log D))_{(2)}$, the sheaf
of germs of local $L^2$-sections in $j_*{\bf
E}^{p-r}\otimes\Omega^r_{\overline M}(\log D)$. (Note that here we
use the induced Hodge metric on ${\bf E}^p$.) One of the main
purposes of this section is then to show the following
\begin{thm}
The sequence
\begin{eqnarray*}
0\to F^{p+1}\Omega^r({\bf H}_{\mathbb{C}})_{(2)}\to
F^{p}\Omega^r({\bf H}_{\mathbb{C}})_{(2)}\to (\Omega^r\otimes {\bf
E}^{p-r})_{(2)}\to 0
\end{eqnarray*}
is exact.
\end{thm}

\noindent {\bf Proof.} We will freely use the notations
established in \S2. In the following discussion, we sometimes
return to the general dimension case, i.e. assuming
$M=({\triangle^*})^n$, since this will not complicate the
discussion.

By the Nilpotent Orbit Theorem, corresponding to the variation,
there exists a nilpotent orbit
$$
{\bf z}(\in {\mathbb{C}}^n)\to
{\text{exp}}(\sum_{i=1}^nz_iN_i)\cdot F\in{\tilde{\bf
D}},~~{\text{for some}}~F\in{\tilde{\bf D}},
$$
which sufficiently approximates the period map near the divisor of
the variation. From \S2, we can produce a series of
$\mathbb{R}$-split mixed Hodge structures $(W^{\bf r}[-k], {\tilde
F}_{\bf r})$ which are canonically attached to $(W^{\bf r}[-k],
e^{\sqrt{-1}N_{r+1}}{\tilde F}_{\bf{r+1}})$ satisfying, for some
$G_{\mathbb{C}}$-value functions $g_{\bf r}({\frac{y_1}{y_{r+1}}},
\cdots , {\frac{y_r}{y_{r+1}}})$ (which map into $G_{\mathbb{R}}$
as ${\frac{y_r}{y_{r+1}}}$ sufficiently large),
\begin{eqnarray*}
&&{\text{exp}}(\sqrt{-1}{\frac{y_r}{y_{r+1}}}(N_r +
\sum_{j=1}^{r-1}
{\frac{y_j}{y_{r}}}N_j))e^{\sqrt{-1}N_{r+1}}{\tilde{F}}_{\bf{r}+1}\\
&&= g_{\bf r}({\frac{y_1}{y_{r+1}}}, \cdots ,
{\frac{y_r}{y_{r+1}}})
{\text{exp}}(\sqrt{-1}{\frac{y_r}{y_{r+1}}}(N_r +
\sum_{j=1}^{r-1}{\frac{y_j}{y_{r}}}N_j)){\tilde{F}}_{\bf{r}},
\end{eqnarray*}
where $W^{\bf r}[-k] = W(N_r +
\sum_{j=1}^{r-1}{\frac{y_j}{y_{r}}}N_j)[-k]$,
$y_j={\text{Im}z_j}$, and $0\le r\le n$; here by $k$ we mean the
nilpotent index of $\{N_j\}$, we will see at the end of this
section that it does not actually exceed the weight of the
variation, i.e. $k$, that is why we here just use $k$.
(Essentially, $g_{\bf r}$ is determined by the following two
1-dimensional nilpotent orbits $\theta_{\bf r}({\bf
z})={\text{exp}}(z(N_r + \sum_{j=1}^{r-1}{\frac{y_j}{y_{r}}}N_j))
e^{\sqrt{-1}N_{r+1}}{\tilde{F}}_{\bf{r}+1}$ and
${\tilde\theta}_{\bf r}({\bf z})={\text{exp}}(z(N_r +
\sum_{j=1}^{r-1}{\frac{y_j}{y_{r}}}N_j)){\tilde{F}}_{\bf{r}}$ with
$$
\theta_{\bf r}({\bf z})=g_{({\frac{y_1}{y_{r}}}, \cdots,
{\frac{y_{r-1}}{y_{r}}}, 1)}({\bf z}){\tilde\theta}_{\bf r}({\bf
z})
$$
so that $g_{\bf r}({\frac{y_1}{y_{r+1}}}, \cdots ,
{\frac{y_r}{y_{r+1}}})=g_{({\frac{y_1}{y_{r}}}, \cdots,
{\frac{y_{r-1}}{y_{r}}}, 1)}({\frac{y_r}{y_{r+1}}})$. For details,
cf. \S2.) By the splitting of $(W^{\bf r}[-k], {\tilde F}_{\bf
r})$, as in \S2, one has a $(0, 0)$-morphism $\tilde{Y}_{\bf r}$,
which together with $N_r + \sum_{j=1}^{r-1}{\frac{y_j}{y_{r}}}N_j$
determines an $sl_2$-triple; in addition, by the $SL_2$-orbit
theorem for several variables, $\tilde{Y}_{\bf r}$ and
$\tilde{N}^{-}_{\bf r}$ together determine another $sl_2$-triple,
where $\tilde{N}^{-}_{\bf r} = \sum_{j=1}^{r}\tilde{N}^{-}_j$ and
$\tilde{N}^{-}_j$ is the component of $N_j$ in the subspace
$\cap_{k=1}^{j-1}{\text{Ker(ad}}{\tilde{Y}}_{\bf k})$.

We now consider
$$
{\text{Ad}}{\text{exp}}(-\sqrt{-1}{\frac{y_r}{y_{r+1}}} (N_r +
\sum_{j=1}^{r-1}{\frac{y_j}{y_{r}}}N_j))(g_{\bf
r}({\frac{y_1}{y_{r+1}}}, \cdots , {\frac{y_r}{y_{r+1}}})), ~~
1\le r\le n,
$$
denoted by $A_{\bf r}$, where $y_{n+1} = 1$ and $y_r > 0$. As in
the proof of 6.20 of \cite{sch}, fixing ${\frac{y_1}{y_{r}}},
\cdots, {\frac{y_{r-1}}{y_{r}}}$, the limit
$\lim_{{\frac{y_r}{y_{r+1}}}\to\infty}A_{\bf r}$ exists, which is
denoted by $g_{\bf r}(\infty)$, and since
$e^{\sqrt{-1}N_{r+1}}{\tilde{F}}_{\bf{r}+1} = A_{\bf r}{\tilde
F}_{\bf r}$, so

\[
e^{\sqrt{-1}N_{r+1}}{\tilde{F}}_{\bf{r}+1} = g_{\bf
r}(\infty){\tilde F}_{\bf r}.
\]
It should be pointed out that although the limit may depend on the
variables ${\frac{y_1}{y_{r}}}, \cdots, {\frac{y_{r-1}}{y_{r}}}$,
this will not affect the following arguments. Important is that
$g_{\bf r}(\infty)$ preserves $W^{\bf r}_l$ and acts as the
identity on $Gr^{W^{\bf r}}_l$ for all $({\frac{y_1}{y_{r}}},
\cdots, {\frac{y_{r-1}}{y_{r}}})$ (cf. \cite{sch}, Lemma 6.20), as
will be used later in this section, so one can fix some
${\frac{y_1}{y_{r}}}, \cdots, {\frac{y_{r-1}}{y_{r}}}$ after
taking the limit each time. Thus we have
\begin{eqnarray*}
&&F = g_{\bf n}(\infty){\tilde F}_{\bf n}
= g_{\bf n}(\infty)e^{-\sqrt{-1}N_n}g_{\bf{n}-1}(\infty){\tilde F}_{\bf{n}-1}\\
&&= \cdots = \prod_{j=n}^1(g_{\bf
j}(\infty)e^{-\sqrt{-1}N_j})(e^{\sqrt{-1}N_1}{\tilde F}_{\bf 1})
:= g(\infty)(e^{\sqrt{-1}N_1}{\tilde F}_{\bf 1}).
\end{eqnarray*}
Denoting $\prod_{j=n}^1(g_{\bf j}(\infty)e^{-\sqrt{-1}N_j})$ by
$g(\infty)$, we have
$$
F=g(\infty)(e^{\sqrt{-1}N_1}{\tilde F}_{\bf 1});
$$
In addition, we also know that the point $e^{\sqrt{-1}N_1}{\tilde
F}_{\bf 1}$ of ${\tilde{\bf D}}$ actually lies in $\bf D$, i.e.
corresponds to a pure Hodge structure on
$H_{\mathbb{C}}=H_{\mathbb{R}}\otimes\mathbb{C}$.

Let now $G_{\mathbb{C}}$ be the complex group corresponding to the
classifying space ${\tilde{\bf D}}$ and $\mathfrak{g}$ its Lie
algebra. Then the $SL_2$-orbit theorem for several variables in
\S2 tells us that there exists a Lie algebra homomorphism $\rho_*:
(\mathfrak{sl}(2, \mathbb{C}))^n\to \mathfrak{g}$, which is {\it
Hodge} at $e^{\sqrt{-1}N_1}{\tilde F}_{\bf 1}$ (called {\it
horizontal} in \cite{sch}). This can be explained as follows.
Using the pure Hodge structure $e^{iN_1}{\tilde F}_{\bf 1}$, one
can define a pure Hodge structure $\{\mathfrak{g}^{p,-p}\}$ of
weight $0$ on $\mathfrak{g}$ relative to $\mathfrak{g}_0$; on the
other hand, let ${\bf n}^+_j, {\bf n}^-_j, {\bf y}_j$ be the
generators of the $j$-th factor of $(\mathfrak{sl}(2,
\mathbb{C}))^n$ as in \S2. Then their images under $\rho_*$ are
$\tilde{N}^+_j, \tilde{N}^-_j, \tilde{Y}_j$ respectively, and
\[
\tilde{N}^+_j + \tilde{N}^-_j + \sqrt{-1}\tilde{Y}_j \in
\mathfrak{g}^{-1,1}, \tilde{N}^+_j - \tilde{N}^-_j \in
\mathfrak{g}^{0,0}.
\]

Now, we will temporarily digress from the proof of the theorem and
establish a slightly abstract setting, which can be just applied
to the above situation. Let $\{H^{p,q}\}$ be a pure Hodge
structure of weight $k$ on $H_{\mathbb{C}} =
H_{\mathbb{R}}\otimes_{\mathbb{R}} \mathbb{C}$, which has an
$(\mathfrak{sl}(2, \mathbb{C}))^n$-action defined over
$\mathbb{R}$ satisfying
\[
X^+_jH^{p,q}\subset H^{p-1,q+1}, X^-_jH^{p,q}
\subset H^{p+1,q-1}, ZH^{p,q}\subset H^{p,q},
\]
where $X^+_j = {\bf n}^+_j+{\bf n}^-_j+\sqrt{-1}{\bf y}_j, X^-_j =
{\bf n}^+_j+{\bf n}^-_j+\sqrt{-1}{\bf y}_j, Z_j = {\bf n}^+_j-{\bf
n}^-_j$, and the action of the $i$-th factor and action of the
$j$-th factor commute, $i\neq j$. We call the
$(\mathfrak{sl}(2,\mathbb{C}))^n$-action {\it Hodge} (or {\it
horizontal}) at $\{H^{p,q}\}$. When the Hodge structure happens to
be polarized by a bilinear form $S$ on $H_{\mathbb{C}}$, one will
say that the $(\mathfrak{sl}(2,{\mathbb{C}}))^n$-action is {\it
compatible with the polarization} if
$(\mathfrak{sl}(2,\mathbb{C})^n$ acts as a Lie algebra of
infinitesimal isometries of $S$. Similar to what Schmid did in
\cite{sch}, p. 258, one can define the invariance and
irreducibility with respect to the given Hodge structure and
horizontal $(\mathfrak{sl}(2,\mathbb{C}))^n$-action.

In the following, we will give an analogous result to Lemma 6.24
of \cite{sch}. First of all, we give some basic Hodge structures
which are irreducible with respect to a certain horizontal
$(\mathfrak{sl}(2,\mathbb{C})^n$-action. For simplicity of the
statement, {\it we just return to the case of $n=2$ again,} as for
the higher dimensional case, one can give a similar statement. The
first basic irreducible structure is denoted by $S(m)\otimes S(n)$
which is of weight $m+n$, where $S(m)$ is the $m$-th symmetric
power of $S(1)$; the latter is defined as follows: Take
${\mathbb{C}}^2 = \mathbb{C}e_1 \oplus \mathbb{C}e_2$ with the
usual $\mathfrak{sl}(2,\mathbb{C})$-action so that $n^-e_2 = e_1,
n^-e_1 = 0$ and the Hodge structure
\begin{eqnarray*}
v^+ &=& e_1 + \sqrt{-1}e_2 ~~~~~{\text{of type $(0,1)$}},\\
v^- &=& e_1 - \sqrt{-1}e_2 ~~~~~{\text{of type $(1,0)$}},
\end{eqnarray*}
where $v^+$ and $v^-$ are the eigenvectors of $\sqrt{-1}(n^+ -
n^-)$ with eigenvalues $1$ and $-1$ respectively. One can also
polarize $\mathbb{C}^2$ over $\mathbb{R}$ as follows: $S(v^+, v^-)
= 2\sqrt{-1}$. It is easy to check that the
$\mathfrak{sl}(2,\mathbb{C})$-action on $\mathbb{C}^2$ is
horizontal with respect to the above Hodge structure and
compatible with the polarization; moreover $\{e_1, e_2\}$ is
orthonormal with respect to the hermitian norm corresponding to
this polarization (cf. \S 2). It is well-known that (polarized)
Hodge structures are compatible with the operations of tensor
products and symmetric products. So, $S(m)\otimes S(n)$ inherits
all structures, including Hodge structure,
$(\mathfrak{sl}(2,\mathbb{C}))^2$-action, and polarization;
moreover the action is horizontal and irreducible with respect to
the induced Hodge structure and compatible with the polarization;
in particular, it follows that $\{e^r_1e^{m-r}_2\otimes
e^s_3e^{n-s}_4\}$ is an orthonormal basis of $S(m)\otimes S(n)$
with respect to the corresponding norm, where ${e_1, e_2}$
corresponds to the first factor of $S(m)\otimes S(n)$, while
$\{e_3, e_4\}$ to the second factor. A direct computation shows
that the induced Hodge structure $\{H^{p,q}, p+q= m+n\}$ on
$S(m)\otimes S(n)$ is
\[
H^{p,q}
=\oplus_{r=1}^{p}\oplus_{s=1}^{q}{\mathbb{C}}(v^-_1)^r(v^+_1)^s\otimes
(v^-_2)^{p-r}(v^+_2)^{q-s},
\]
where $\{v^+_1, v^-_1\}$ and $\{v^+_2, v^-_2\}$ correspond to the
first and second factor of $S(m)\otimes S(n)$ respectively, namely
\begin{eqnarray*}
v_1^+ &=& e_1 + \sqrt{-1}e_2 ~~{\text{and}}~~ v_2^+ = e_3 + \sqrt{-1}e_4
~~~~~{\text{of type $(0,1)$}}\\
v_1^- &=& e_1 - \sqrt{-1}e_2 ~~{\text{and}}~~ v_2^- = e_3 -
\sqrt{-1}e_4 ~~~~~{\text{of type $(1,0)$}}.
\end{eqnarray*}

In the following we continue to consider some trivial Hodge
structures on which $(\mathfrak{sl}(2, \mathbb{C}))^2$ acts
horizontally and irreducibly. The one-dimensional complex vector
space $\mathbb{C}$, with the obvious real structure, carries a
unique Hodge structure of weight $2$. Deligne denotes it by $H(1)$
and calls it the "Hodge structure of Tate". Take the standard
basis $h^{1,1}$, which is of type $(1, 1)$, so that
$H(1)=\mathbb{C}h^{1,1}$. For $n\ge 0$, $H(n)$ shall be the $n$-th
symmetric power of $H(1)$, and $H(-n)$ the dual of $H(n)$, the
standard of which will be denoted by $h^{n,n}$.
$(\mathfrak{sl}(2,\mathbb{C}))^2$ trivially acts on $H(n)$ for all
integers $n$. Also, one can polarize $H(1)$ naturally: the
standard basis $h^{1,1}$ enables us to identify $H(1)$ with
$\mathbb{C}$, namely $h^{1,1}$ with $1$; then the nondegenerate
bilinear form $S$ on $H(1)$ is taken as $S(1, 1)=1$, from which
$H(n)$ inherits the standard polarization, still denoted by $S$,
for all integers $n$. It is easy to see that the
$(\mathfrak{sl}(2,\mathbb{C}))^2$-action on $H(n)$ is horizontal
and irreducible and  compatible with the polarization. The second
slight nontrivial Hodge structure is defined on $\mathbb{C}^2$ as
follows: let $e_1, e_2$ be the standard basis vectors of
$\mathbb{C}^2$, which gives the standard real structure on
$\mathbb{C}^2$. For $p\neq q$, we define a Hodge structure
$E(p,q)$ of weight $p+q$ on $\mathbb{C}^2$ by requiring that
\[
e^{p,q} = e_1 - \sqrt{-1}e_2
\]
be of type $(p,q)$, and
\[
e^{q,p} = e_1 + \sqrt{-1}e_2
\]
of type $(q,p)$. Again $(\mathfrak{sl}(2,\mathbb{C}))^2$ trivially
acts on $\mathbb{C}^2$ and the action is horizontal and
irreducible with respect to the Hodge structure $E(p,q)$. The
bilinear form $S$ on $\mathbb{C}^2$, which is described by the
identities
\begin{eqnarray*}
S(e^{p,q}, e^{p,q}) &=& 0, ~~~~~~~~S(e^{q,p}, e^{q,p}) = 0\\
S(e^{p,q}, e^{q,p}) &=& 2\sqrt{-1}^{q-p}, ~~~S(e^{q,p}, e^{p,q}) =
2\sqrt{-1}^{p-q},
\end{eqnarray*}
polarizes $E(p, q)$. It is easy to see that the polarization is
compatible with the action. Now, we are in the position to state
the following generalization of a result of Schmid (cf.
\cite{sch}, Theorem 6.24) to arbitrary dimensions (stated here for
two dimension for simplicity).
\begin{thm}
Let $H_{\mathbb{C}} = H_{\mathbb{R}}\otimes \mathbb{C}$ be a
complex space with a Hodge structure \{$H^{p,q}\}$ of weight $k$
and a horizontal $(\mathfrak{sl}(2,\mathbb{C}))^2$-action. Then
$H_{\mathbb{C}}$ can be decomposed into a direct sum of subspaces
which are invariant and irreducible with respect to the given
structures: the Hodge structure and the horizontal
$(\mathfrak{sl}(2,\mathbb{C}))^2$-action. Every irreducible factor
is isomorphic - relative to the Hodge structure and the horizontal
action - to one of the following types: $H(l)\otimes S(m)\otimes
S(n)$, with $l\in\mathbb{Z}, m, n\ge 0$, and $k = 2l+m+n$; or
$E(p,q)\otimes S(m)\otimes S(n)$, with $p>q, m, n\ge 0$, and $k =
p+q+m+n$. If the Hodge structure of $H_{\mathbb{C}}$ happens to
have a polarization which is compatible with the
$(\mathfrak{sl}(2,\mathbb{C}))^2$-action, then the decomposition
can be chosen to be orthogonal with respect to the polarization,
and the isomorphisms between the irreducible factors and the
irreducible structures of special type can be chosen with the
further restriction that they should preserve the polarizations.
\end{thm}

\noindent {\bf Proof.} Since $\{Z_1, Z_2\}$ commute and act
semisimply on each Hodge subspace $H^{p,q}$, with integral
eigenvalues, one can find a basis of $H^{p,q}$, each element of
which is an eigenvector of both $Z_1$ and $Z_2$ with integral
eigenvalue. Then, we can define a linear $SO(2)$-action on
$H_{\mathbb{C}}$ as follows: let $v\in H^{p,q}$ be an eigenvector
of both $Z_1$ and $Z_2$ with eigenvalues $l_1$ and $l_2$
respectively, then the element $e^{\sqrt{-1}\theta}\in SO(2)$ acts
on $v$ as multiplication by $e^{\sqrt{-1}(l_1+l_2+p-q)\theta}$,
denoted by $e^{\sqrt{-1}\theta}.v$. Obviously the action is
orthogonal with respect to the polarization and trivial on
$S(m)\otimes S(n)$. (Since $v=(v^-_1)^r(v^+_1)^{m-r}\otimes
(v^-_2)^s(v^+_2)^{n-s}\in H^{r+s,m+n-r-s}\subset S(m)\otimes S(n)$
is an eigenvector of $Z^+_1$ and $Z^+_2$ with eigenvalues $m-2r$
and $n-2s$ respectively, so $e^{\sqrt{-1}\theta}$ acts on $v$ by
multiplying
$e^{\sqrt{-1}((m-2r)+(n-2s)+(r+s)-(m+n-r-s))\theta}=1$).
Considering $SO(2)$ as the group of real points of an algebraic
$1$-torus $T$, which is defined and {\it anisotropic} over
$\mathbb{R}$, one can extend the $SO(2)$-action to a
representation of $T$ on $H_{\mathbb{C}}$ over $\mathbb{R}$. In
addition, the $(\mathfrak{sl}(2, {\mathbb{C}}))^2$-action on
$H_{\mathbb{C}}$ also determines a representation on
$H_{\mathbb{C}}$ of the algebraic group $(SL(2,{\mathbb{C}}))^2$
again defined over $\mathbb{R}$. Since the $(\mathfrak{sl}(2,
{\mathbb{C}}))^2$-action is horizontal, the representations of $T$
and $(SL(2,\mathbb{C}))^2$ commute: let $v\in H^{p,q}$ be an
eigenvector of both $Z_1$ and $Z_2$ with the eigenvalues $l_1$ and
$l_2$ respectively, then $X^+_1v$ ($X^+_2v$) is an eigenvector of
both $Z_1$ and $Z_2$ with the eigenvalues $l_1+2$ ($l_1$) and
$l_2$ ($l_2+2$) respectively. So,
\begin{eqnarray*}
&&e^{\sqrt{-1}\theta}.(X^+_1v) (e^{\sqrt{-1}\theta}.(X^+_2v))\\
&&= e^{\sqrt{-1}((l_1+2)+l_2+(p-1)-(q+1))\theta}X^+_1v
(e^{\sqrt{-1}(l_1+(l_2+2)+(p-1)-(q+1))\theta}X^+_2v)\\
&&= e^{\sqrt{-1}(l_1+l_2+p-q)\theta}X^+_1v (e^{\sqrt{-1}(l_1+l_2+p-q)\theta}X^+_2v)\\
&&= X^+_1(e^{\sqrt{-1}\theta}.v) (X^+_2(e^{\sqrt{-1}\theta}.v))
\end{eqnarray*}
and similarly one can show the other commutativity, where we use
the horizontality of the $(\mathfrak{sl}(2,
{\mathbb{C}}))^2$-action. Hence one has a representation of the
product $T\times SL(2, \mathbb{C})$. One can also show that a
$\{Z_1, Z_2\}$-stable subspace of $H_{\mathbb{C}}$ carries a
sub-Hodge structure if and only if it is self-conjugate and
invariant under the $SO(2)$-action: the proof of necessity  is
obvious; we prove the sufficiency. Let $H^1_{\mathbb{C}}$ be a
$\{Z_1,Z_2\}$-stable subspace of $H_{\mathbb{C}}$ which is
self-conjugate and invariant under the $SO(2)$-action, $v\in
H^1_{\mathbb{C}}$. Decompose $v=\sum_{p+q=k, l_1,l_2\in {\bf
Z}}v^{p,q}_{l_1,l_2}$, $v^{p,q}_{l_1,l_2}\in H^{p,q}$ is an
eigenvector of both $Z_1$ and $Z_2$ with eigenvalue $l_1$ and
$l_2$ respectively, since $e^{\sqrt{-1}\theta}.v\in
H^1_{\mathbb{C}}$ for all $\theta\in{\mathbb{R}}$,
$\sum_{l_1+l_2+p-q=s}v^{p,q}_{l_1,l_2}\in H^1_{\mathbb{C}}$. Using
the $\{Z_1, Z_2\}$-stability of $H^1_{\mathbb{R}}$, one can then
show that each $v^{p,q}_{l_1,l_2}$ lies in $H^1_{\mathbb{C}}$, and
hence $H^1_{\mathbb{R}}$ carries a sub-Hodge structure. Thus, one
gets that the subspaces of $H_{\mathbb{C}}$ which are invariant
and irreducible with respect to the Hodge structure and the
$(\mathfrak{sl}(2,{\mathbb{C}}))^2$-action correspond bijectively
to those invariant and irreducible subrepresentations of the
representation of $T\times (SL(2,{\mathbb{C}}))^2$, that are
defined over $\mathbb{R}$. Because the reductivity of the product
$T\times (SL(2,{\mathbb{C}}))^2$, this proves the first part of
the theorem.

An irreducible representation defined over $\mathbb{R}$ of the
group $T\times (SL(2,{\mathbb{C}}))^2$ either remains irreducible
under $(SL(2,{\mathbb{R}}))^2$, in which case $T$ acts trivially,
or splits into two conjugate subspaces, each of which is
$T$-stable and $(SL(2,{\mathbb{C}}))^2$-irreducible, with $T$
acting nontrivially. The first case corresponds to an irreducible
Hodge structure with respect to a horizontal
$(\mathfrak{sl}(2,{\mathbb{C}}))^2$-action of the type
$H(l)\otimes S(m)\times S(n)$ with $k=2l+m+n$, $l\in{\mathbb{Z}},
m,n\ge 0$, which has dimension $(m+1)(n+1)$; the second case
corresponds to an irreducible Hodge structure of the type
$E(p,q)\otimes S(m)\otimes S(n)$, with $p>q, m, n\ge 0$, and $k =
p+q+m+n$, which has dimension $2(m+1)(n+1)$. We omit their proofs
and the proof of the remaining part of the theorem, since they are
standard.$\Box$

We now return to the proof of the theorem at the beginning of this
section. Actually, we only need to prove the case of $r=0$. By the
previous discussions, $e^{\sqrt{-1}N_1}{\tilde F}_{\bf 1}$
(briefly denoted by $F_0$ below) has an $({\mathfrak{sl}}(2,
{\mathbb{C}}))^2$-horizontal action, which is compatible with the
polarization. So, by the above decomposition theorem, we can
assume that $e^{\sqrt{-1}N_1}{\tilde F}_{\bf 1}$ is irreducible.
Since $H(l)$ and $E(p,q)$ are trivial, without loss of generality,
we assume that $e^{\sqrt{-1}N_1}{\tilde F}_{\bf 1} = S(m)\otimes
S(n)$. We choose
\[
{\alpha}_{k,l} = ({\tilde N}^-_1)^{m-k}({\tilde N}^-_2)^{n-l}((v^-_1)^m\otimes(v^-_2)^n)
\]
as a basis of $S(m)\otimes S(n)$. It is easily to check that this basis satisfies
\[
{\alpha}_{k,l}\in F^{k+l}_0\cap W^{\bf 1}_{2k}\cap W^{\bf 2}_{2(k+l)}
\]
and ${\alpha}_{k,l}$ projects nontrivially in $Gr_{F_0}^s =
F^s_0/F^{s+1}_0$, $Gr_{2k}^{W^{\bf 1}} = W^{\bf 1}_{2k}/W^{\bf
1}_{2k-1}$, and $Gr^{W^{\bf 2}}_{2s} = W^{\bf 2}_{2s}/W^{\bf
2}_{2s-1}$ with $s=k+l$, where $\{W^{\bf 1}_k\}$ and $\{W^{\bf
2}_s\}$ are the weight filtrations corresponding to ${\tilde
N}^-_1=N_1$ and ${\tilde N}^-_{\bf 2} = {\tilde N}^-_1 + {\tilde
N}^-_2$ (or $c_1N_1+c_2N_2$) respectively; $W^{\bf 1}_k$ and
$W^{\bf 2}_s$ can actually be described as the subspaces spanned
by eigenvectors of ${\tilde Y}_1$ and ${\tilde Y}_{\bf 2} =
{\tilde Y}_1 + {\tilde Y}_2$ with eigenvalues not greater than
$k-m$ and $s-m-n$ respectively \cite{sch}. (Note that the
notations $W^{\bf j}$ here are just $W^{\bf j}[-k]$ in \S2 and the
same for the remaining part of this section.) Since $F = g(\infty
)F_0$, one has $g(\infty ){\alpha}_{k,l}\in F^{k+l}$ and it
nontrivially projects in $Gr_{F}^s = F^s/F^{s+1}$ with $s=k+l$; by
the means of the construction of $g(\infty )$, one also has that
$g(\infty )({\alpha}_{k,l})$ can be considered as the basis of
$Gr_{2k}^{W^{\bf 1}}\cap Gr^{W^{\bf 2}}_{2(k+l)} = {W^{\bf
1}_{2k}\cap W^{\bf 2}_{2(k+l)}}/ ({W^{\bf 1}_{2k-1}\cap W^{\bf
2}_{2(k+l)} + W^{\bf 1}_{2k}\cap W^{\bf 2}_{2(k+l)-1}})$, since
$g(\infty )$ acts as the identity on $Gr_{2k}^{W^{\bf 1}}\cap
Gr^{W^{\bf 2}}_{2(k+l)}$.

On the other hand, by the nilpotent orbit theorem, using the
notation in the previous sections, one can easily show the
following assertion: {\it Fixing $x\in M$, for $v\in ({\bf
H_{\mathbb{C}})_x}$, $v\in F^p$ if and only if there exist two
holomorphic sections $w_1$ and $w_2$ of $\overline{\bf
H}_{\mathbb{C}}$ such that ${\tilde v} + t_1w_1 + t_2w_2$ is a
section of ${\bf F}^p$.} So, for $g(\infty ){\alpha}_{k,l}$, there
exist two sections ${\beta}^1_{k,l}$ and ${\beta}^2_{k,l}$ of
$\overline{\bf H}_{\mathbb{C}}$ such that
\[
{\sigma}_{k,l} = \widetilde{g(\infty )({\alpha}_{k,l})}
+ t_1{\beta}^1_{k,l} + t_2{\beta}^2_{k,l}
\]
is a section of ${\bf F}^s$ and, under the natural projection,
represents a generator for ${\bf E}^{s}$ for $s=k+l$.

Summing up all the above, one has that $\{{\sigma}_{k,l}| k+l =
s\}$, under the natural projection, represents a set of generators
for ${\bf E}^{s}$, while $\{{\sigma}_{k,l}| k+l \ge s\}$
represents a set of generators for ${\bf F}^{s}$. Furthermore, by
the norm estimates of the previous section, since
$\widetilde{g(\infty )({\alpha}_{k,l})}\in {\overline {\bf
W}}^{\bf 1}_{2k}\cap {\overline {\bf W}}^{\bf 2}_{2(k+l)}$ and
${\Vert t_j{\beta}^j_{k,l}\Vert}^2 =
O(|t_j|^2|\text{log}|t_1||^m|\text{log}|t_2||^n)$, one has
\begin{eqnarray*}
{\Vert{\sigma}_{k,l}\Vert}^2 &\sim& ({\frac{\log|t_1|}
{\log|t_2|}})^{2k-m}\cdot(-\log|t_2|)^{2(k+l-m-n)} \\
&=& (-{\text{log}}|t_1|)^{2k-m}(-{\text{log}}|t_2|)^{2l-n},
\end{eqnarray*}
on the domains of the form $D_{\epsilon} = \{(t_1, t_2)\in
(\triangle^*)^2~|~{\frac{\log|t_1|}{\log|t_2|}}>\epsilon,
-\log|t_2|>\epsilon \}$. (Note again that the notations
${\overline{\bf W}^{\bf j}}$ here are just ${\overline{\bf W}^{\bf
j}}[-k]$ in \S2.) Again since, by the previous section, one can
verify that for $v\in {\bf W}^{\bf 1}_k\cap{\bf W}^{\bf 2}_s$
\begin{eqnarray*}
L({\bf z})({\tilde{v}}) &=& {\tilde{h}}({\bf z})
{\text{exp}}({\frac 1 2}{\text{log}}{\frac{y_1}{y_2}}{\tilde Y}_{\bf 1} +
{\frac 1 2}{\text{log}}y_2{\tilde Y}_{\bf 2})
{\text{exp}}(-\sum_{j=1}^{2}x_jN_j){\tilde v}\\
&=& {\tilde{h}}({\bf z})
[({\frac{y_1}{y_2}})^{\frac k 2}y_2^{\frac{s}{2}}
(\sum_{j=0}^{m+n}{\frac{{\sqrt{-1}}^j}{j!}}
(N_1+N_2)^j)v + (\text{lower order terms})]\\
&=& {\tilde{h}}({\bf z})
[({\frac{y_1}{y_2}})^{\frac k 2}y_2^{\frac{s}{2}}v + \cdots ],
\end{eqnarray*}
where $t_j=e^{2\pi iz_j}$, $z_j=x_j+iy_j$ for $j=1,2$, and
${\tilde{h}}({\bf z})$ is strongly asymptotic to the identity as
${\frac{y_1}{y_2}}, y_2\to\infty$. Thus, $L({\bf
z})(\widetilde{g(\infty )({\alpha}_{k,l})})$ and $L({\bf
z})({\sigma}_{k,l})$ are asymptotically the same. So, applying
Lemma 2, one has that $\{{\sigma}_{k,l}| k+l \ge s\}$ gives an
$L^2$-adapted basis for ${\bf F}^s$. Therefore, by the above norm
estimates of ${\sigma}_{k,l}$, $({\bf F}^s)_{(2)}$ is freely
generated by the sections
$\{t_1^{{\epsilon}_k}t_2^{{\eta}_l}{\sigma}_{k,l} | k+l \ge s\}$
with ${\epsilon}_k = 0$ if $2k\le m$ and ${\epsilon}_k = 1$
otherwise and ${\eta}_l = 0$ if $2l\le n$ and ${\eta}_l = 1$
otherwise. On the other hand, by the above expression of $L({\bf
z})({\tilde{v}})$, one knows that the $(k+l, m+n-k-l)$ component
of $L({\bf z})({\sigma}_{k,l})$ is asymptotic to
$y_1^{\frac{2k-m}{2}}y_2^{\frac{2l-n}{2}}{\alpha}_{k,l}$ (Here,
one needs to use the expansion of $g(\infty )$, see \cite{sch},
the proof of Theorem 6.20), namely ${\sigma}_{k,l}$ carries much
of its norm in this Hodge component. Therefore, the projections in
${\bf E}^s$ of the elements in the set
$\{t_1^{{\epsilon}_k}t_2^{{\eta}_l}{\sigma}_{k,l} | k+l = s\}$
form a set of generators of ${\bf E}^s$. So, one naturally has
$({\bf F}^s)_{(2)}/({\bf F}^{s+1})_{(2)} = {\bf E}_{(2)}$. This
completes the proof of the
theorem.~~~~~~~~~~~~~~~~~~~~~~~~~~~~~~~~~~~~~~~~~~~~~~~~~~~~~$\Box$

\vskip .3cm From the above arguments, we also obtain the following
byproducts.

\begin{thm}
${\nabla}^{1,0}$ (hence also $\theta$) is a bounded operator under
the Poincar\'e-like metric on $(\triangle^*)^2$ and the Hodge
metric on ${\bf H}_{\mathbb{C}}$ (the induced Hodge metric on the
Hodge bundle).
\end{thm}

\noindent {\bf Proof}. From the previous arguments, we know that
$\{{\sigma}_{k,l}| k+l \ge 0\}$ represents a basis for ${\bf
H}_{\mathbb{C}}$ and each ${\sigma}_{k,l}$ can be expressed as
\[
{\sigma }_{k,l} = \widetilde{g(\infty )({\alpha}_{k,l})}
+ t_1{\beta}^1_{k,l} + t_2{\beta}^2_{k,l}
\]
for some holomorphic sections ${\beta}^1_{k,l}$ and
${\beta}^2_{k,l}$ of ${\overline{\bf H}}_{\mathbb{C}}$. By the
definition of ${\overline{\bf H}}_{\mathbb{C}}$ in the Section 2
and 3, we have
\[
{\nabla}^{1,0}({\sigma }_{k,l}) = {\frac{1}{2\pi i}}
({\frac{{\text d}t_1}{t_1}}\otimes N_1 + {\frac{{\text
d}t_2}{t_2}}\otimes N_2)\widetilde{g(\infty )({\alpha}_{k,l})} +
{\nabla}^{1,0}(t_1{\beta}^1_{k,l}) +
{\nabla}^{1,0}(t_2{\beta}^2_{k,l}).
\]
So, if we neglect the higher order terms, we can consider
${\nabla}^{1,0}$ as ${\frac{1}{2\pi i}}({\frac{{\text
d}t_1}{t_1}}\otimes N_1 + {\frac{{\text d}t_2}{t_2}}\otimes N_2)$
acting on ${\bf H}_{\mathbb{C}}$ by multiplication. On the other
hand, by means of the commutativity of $N_1, N_2$ and the
definitions of $\{W^{\bf 1}_{l_1}\}$ and $\{W^{\bf 2}_{l_2}\}$ in
the Sections 2 and 3, one has that $N_1$ lowers the weights of
$\{W^{\bf 1}_{l_1}\}$ and $\{W^{\bf 2}_{l_2}\}$ by $2$, while
$N_2$ lowers the weight of $\{W^{\bf 2}_{l_2}\}$ by $2$ and
preserves $\{W^{\bf 1}_{l_1}\}$. Thus, again by Theorem 5.21 of
\cite{cks1}, one has that, if ${\tilde v}\in {\overline {\bf
W}}^{\bf 1}_{l_1}\cap {\overline {\bf W}}^{\bf 2}_{l_2}$ and the
parts of the highest weights of $v$ in $\{W^{\bf 1}_{l_1}\}$ and
$\{W^{\bf 2}_{l_2}\}$ are nonzero,
\[
{\Vert({\frac{1}{2\pi i}}{\frac{{\text d}t_1}{t_1}}
\otimes N_1){\tilde v}\Vert}^2\sim
{{\text log}^2|t_1|}({\frac{{\text log}|t_1|}
{{\text log}|t_2|}})^{l_1-2}({{\text log}|t_2|})^{l_2-2}
\]
and
\[
{\Vert({\frac{1}{2\pi i}}{\frac{{\text d}t_2}{t_2}}
\otimes N_2){\tilde v}\Vert}^2\sim
{{\text log}^2|t_2|}({\frac{{\text log}|t_1|}
{{\text log}|t_2|}})^{l_1}({{\text log}|t_2|})^{l_2-2},
\]
on $D_{\epsilon}$. Here, we use ${\Vert{\frac{{\text
d}t_1}{t_1}}\Vert}^2\sim {{\text log}^2|t_1|}$ and
${\Vert{\frac{{\text d}t_2}{t_2}}\Vert}^2\sim {{\text
log}^2|t_2|}$, under the Poincar\'e-like metric of $M$ as defined
in the beginning of this section. Hence,
$\Vert{\nabla}^{1,0}\Vert\sim 1$ on $D_{\epsilon}$, namely,
${\nabla}^{1,0}$ is a bounded operator on $D_{\epsilon}$, though
the bound may depend on $\epsilon$. Note that $\nabla^{1,0}$ is
independent of the choice of order $(t_1, t_2)$, so by applying
the same argument to the order $(t_2, t_1)$, we can also show that
$\nabla^{1,0}$ is a bounded operator on the domain $D'_{\epsilon}
= \{(t_1, t_2)~|~{\frac{\log|t_2|}{\log|t_1|}}>\epsilon,
-\log|t_1|>\epsilon\}$. It is clear that for a sufficient small
$\epsilon > 0$, $D_{\epsilon}\cup D'_{\epsilon}$ contains a
neighborhood of the divisor, so $\nabla^{1,0}$ is bounded. The
boundedness of $\theta$ is obtained by using the projection to the
Hodge bundle and the $L^2$-adaptedness theorem.~~~~~~~~$\Box$

\vskip 4mm \noindent {\bf Remarks.} It is easy to see that the
residue $res({\nabla}^{1,0})$ of ${\nabla}^{1,0}$ (according to
the notations of \cite{si1}) at the singularity is $\{N_1, N_2\}$;
again, $\{{\sigma}_{k,l}| k+l \ge s\}$ represent a basis for ${\bf
F}_s$ and the suitable projections of the elements of
$\{{\sigma}_{k,l}| k+l=s\}$ represents a basis of ${\bf E}_s$, so,
by the above formula of the action of ${\nabla}^{1,0}$ on
${\sigma}_{k,l}$, the residue $res(\theta )$ of the induced map
$\theta$ (Higgs field) of ${\nabla}^{1,0}$ on ${\bf E}_s$ (Hodge
bundles, $={\bf F}_s/{\bf F}_{s+1}$) is also $\{N_1, N_2\}$ .
Namely, we have the following corollary, which is essentially due
to Schmid and Simpson \cite{sch, si1}.

\begin{cor}
The three weight filtrations corresponding to
$res({\nabla}^{1,0})$ of ${\nabla}^{1,0}$ (at the singularity,
considered as acting on ${\bf H}_{\mathbb{C}}$), $res(\theta )$ of
$\theta$ (at the singularity, considered as acting on the Hodge
bundles) and $\{N_1, N_2\}$ coincide under some suitable
identifications.
\end{cor}

In the above argument, we use the ordering $(N_1, N_2)$ to
construct an $L^2$-adaptedness basis $\{\sigma_{k,l}|k+l=s\}$ for
the Hodge bundle ${\bf E}^s$; similarly, one can also use the
other ordering $(N_2, N_1)$ to construct the corresponding basis
of ${\bf E}^s$. As a conclusion of this section, we will compare
the two bases. In the following, a monodromy weight filtration is
again considered as some translation (as in the previous) of the
usual one as defined in the \S 2. From the \S 2 and the above
discussion, without loss of generality, we will assume in the
following that $(F, W=W(N_1+N_2))$ is an $\mathbb{R}$-split mixed
Hodge structure. We have two different orderings: $(N_1, N_2)$ and
$(N_1, N_2)$ corresponding to the orderings of the coordinates
$t_1, t_2$ on $(\triangle^*)^2$. In the previous discussion, we
used the first ordering $(N_1, N_2)$ and constructed a basis
$\{{\alpha}_{k,l}\}$ (up to the action of some elements in
$G_{\mathbb{C}}$) of the generic fiber $H_{\mathbb{C}}$ of the
pull-back of ${\bf H}_{\mathbb{C}}$ to the universal covering of
$(\triangle^*)^2$, which flags $F$, $W(N_1)$, and $W$; actually,
under the assumption of $(F, W)$ being $\mathbb{R}$-split, one can
furthermore assume that $\{{\alpha}_{k,l}\}$ corresponds to
simultaneous ${\mathbb{Z}}^2$-gradings of $W(N_1)$ and $W$ over
$\mathbb{R}$, namely ${\alpha}_{k,l}$ being the eigenvector of
${\tilde{Y}}_{\bf 1}$ and ${\tilde{Y}}_{\bf 2}$ with eigenvalues
$2k-m$ and $2(k+l)-m-n$ respectively. Similarly, for the ordering
$(N_2, N_1)$, we can also construct a basis $\{{\alpha}'_{k,l}\}$
of $H_{\mathbb{C}}$, which flags $F$, and grades simultaneously
$W(N_2)$ and $W$ in the same way as above (of course for a certain
other semisimple element ${\tilde{Y}}_{\bf 1}'$ corresponding to
$N_2$, and ${\tilde{Y}}_{\bf 2}$). In the following, we will first
discuss the relation between the both bases. For convenience, we
here state the following
\begin{prop}(\cite{ck}, Theorem 3.3(2); \cite{ca}, Theorem 3)
{\footnote{This is a correction to the original statement in
\cite{ck} by S. Zucker; for this, see {\it Mathematical Reviews},
84a:32046 by S. Zucker.}}) Let $W(N_2, {\text{Gr}}^{W(N_1)}_k)$ be
the weight filtration on ${\text{Gr}}^{W(N_1)}_k$ corresponding to
the induced nilpotent endomorphism of $N_2$ from
${\text{Gr}}^{W(N_1)}_k$ to itself (still denoted by $N_2$). Then
the projection of $W$ to ${\text{Gr}}^{W(N_1)}_k$ just induces
$W(N_2, {\text{Gr}}^{W(N_1)}_k)[-k]$.
\end{prop}

By the previous construction, $\alpha_{k,l}\in (W(N_1))_{2k}\cap
W_{2(k+l )}$ (resp. $\alpha'_{k,l}\in (W(N_2))_{2l}\cap W_{2(k+l
)}$), so $\oplus_l{\mathbb{C}}{\alpha}_{k,l}$ can be identified
with ${\text{Gr}}^{W(N_1)}_{2k}$. Accordingly, $H_{\mathbb{C}}$
can be identified with $\oplus_{k}{\text{Gr}}^{W(N_1)}_{2k}$.
Applying the above proposition to each piece
$({\text{Gr}}^{W(N_1)}_{2k}, N_2)$ and then summing up for $k$, we
have $W_{2l}(N_2)=\oplus_{l'\le l,k}{\mathbb{C}}{\alpha}_{k,l'}$;
on the other hand, the above construction tells us that
$W_{2l}(N_2)=\oplus_{l'\le l,k}{\mathbb{C}}{\alpha}'_{k,l'}$.
Thus, each $\alpha_{k,l}$ can be expressed linearly by
$\{\alpha'_{k',l'}; l'\le l, 0\le k'\le m\}$ and $\alpha'_{k,l}$
also by $\{\alpha_{k',l'}; l'\le l, 0\le k'\le m\}$. The same
discussion also applies to the ordering $(N_2, N_1)$ and one has
that each $\alpha_{k,l}$ can be expressed linearly by
$\{\alpha'_{k',l'}; k'\le k, 0\le l'\le n\}$ and $\alpha'_{k,l}$
also by $\{\alpha_{k',l'}; k'\le k, 0\le l'\le n\}$. Therefore,
each $\alpha_{k,l}$ can be expressed linearly by
$\{\alpha'_{k',l'}; k'\le k, l'\le l\}$ and $\alpha'_{k,l}$ also
by $\{\alpha_{k',l'}; k'\le k, l'\le l\}$. Namely, we have
${\alpha}'_{k,l}=\sum_{k'\le k, l'\le
l}{c_{k',l'}}{\alpha}_{k',l'}$ and ${\alpha}_{k,l}=\sum_{k'\le k,
l'\le l}{c'_{k',l'}}{\alpha}'_{k',l'}$ for some complex number
$c_{k',l'}, c'_{k',l'}$. Denote the flat sections generated by
$\alpha_{k,l}$ and $\alpha'_{k,l}$ by the same symbols. Denote the
monodromized sections of $\alpha_{k,l}$ and $\alpha'_{k,l}$ by
$\sigma_{k,l}$ and $\sigma'_{k,l}$.(Here and in the following, we
always neglect the higher order terms.) It is clear by the
previous construction of an $L^2$-adapted basis for the Hodge
bundles that $\{\sigma_{k,l}\}$ and $\{\sigma'_{k,l}\}$ are two
$L^2$-adapted bases. Since $\alpha_{k,l}\in (W(N_1))_{2k}\cap
W_{2(k+l )}$ (resp. $\alpha'_{k,l}\in (W(N_2))_{2l}\cap W_{2(k+l
)}$), by the estimate of the Hodge norm, we have
\[
\Vert\sigma_{k,l}\Vert^2\sim
({\log}|t_1|)^{2k-m}({{\log}|t_2|})^{2l-n},
~~{\text{on}}~D_{\epsilon}
\]
and
\begin{eqnarray*}
\Vert\sigma'_{k,l}\Vert^2&\sim&({\frac{{\log}|t_2|}
{{\log}|t_1|}})^{2l-n}({{\log}|t_1|})^{2(l+k)-m-n}  \\
&=&({\log}|t_1|)^{2k-m}({{\log}|t_2|})^{2l-n},
~~{\text{on}}~D'_{\epsilon}.
\end{eqnarray*}
On the other hand, since
${\alpha}_{k,l}=\sum_{k'\le k, l'\le l}{c'_{k',l'}}{\alpha}'_{k',l'}$
for some complex number $c'_{k',l'}$,
so $\sigma_{k,l}=\sum_{k'\le k, l'\le l}{c'_{k',l'}}{\sigma}'_{k',l'}$
(since both the parallel translation and the monodromization are linear).
Thus, one has, on $D'_{\epsilon}$,
\begin{eqnarray*}
\Vert\sigma_{k,l}\Vert^2&\sim&
\sum_{k'\le k,l'\le l}(c'_{k',l'})^2\Vert\sigma_{k',l'}\Vert^2 \\
&\sim&(\sum_{k'\le k,l'\le l}(c'_{k',l'})^2{\log}|t_1|)^{2k'-m}
({{\log}|t_2|})^{2l'-n}\\
&\sim&{\log}|t_1|)^{2k-m}({{\log}|t_2|})^{2l-n}.
\end{eqnarray*}
Thus, we obtain the following
\begin{prop}
Near the singularity, there exists a holomorphic basis
$\{{\sigma}_{k,l}~|~k+l=s\}$ for each Hodge bundle ${\bf E}^s$
that satisfies
\[
\Vert\sigma_{k,l}\Vert^2\sim ({\log}|t_1|)^{2k-m}({{\log}|t_2|})^{2l-n}
\]
near the singularity.
\end{prop}

\section{$L^2$-holomorphic Dolbeault complex}

In this section, we still assume that $M=(\triangle^*)^2$ and
$\overline M=(\triangle )^2$, $j:M\to{\overline M}$ is the
inclusion map. Let $(M, {\bf H}_{\mathbb{C}}={\bf
H}_{\mathbb{Z}}\otimes{\mathbb{C}}, \nabla = {\nabla}^{1,0} +
{\nabla}^{0,1} , {\bf F}=\{{\bf F}^p\}, {\bf S})$ be a variation
of polarized Hodge structure with weight $m$ defined over
$\mathbb{Q}$ such that each ${\bf F}^p$ is a holomorphic subbundle
of the local system ${\bf H_{\mathbb{C}}}$ and
\[
{\nabla}^{1,0}{\bf F}^p\subset {\bf F}^{p-1}\otimes{\Omega}^1(M).
\]
Let ${\gamma}_1$ and ${\gamma}_2$ be the monodromy transformations
of the variation, which are always assumed to be unipotent. Let
$N_1$ and $N_2$ be the logarithmic monodromies of ${\gamma}_1$ and
${\gamma}_2$ respectively, which are nilpotent. Note that in the
following, we consider ${\bf H_{\mathbb{C}}}$ as a flat bundle, a
local system, a holomorphic vector bundle relative to
$\nabla^{0,1}$, or a sheaf of local holomorphic sections, this
depends on the context.

In the previous section, using the Poincar\'e-like metric on $M$
and the Hodge metric on ${\bf H_{\mathbb{C}}}$, we define
$\Omega^r({\bf H_{\mathbb{C}}})_{(2)}$ to be the sheaf of local
$L^2$ holomorphic $r$-forms valued in $j_*{\bf H_{\mathbb{C}}}$ on
$\overline M$. One of the purposes of this section will be to show
that the sheaf $\Omega^r({\bf H_{\mathbb{C}}})_{(2)}$ can be
defined algebraically, just using the logarithmic monodromies
$N_1$ and $N_2$, and lies in $j_*{\bf
H}_{\mathbb{C}}\otimes\Omega^r_{\overline M}({\text{log}}D)$. As a
consequence of this fact together with the asymptotic behavior of
$\theta$ and the $L^2$-adaptedness theorem of \S3, we obtain the
$L^2$ holomorphic Dolbeault complex on $\overline M$ in \S1 (for
the precise definition, see this \S)
\[
(*)~~~~~~~~~~~~~~~~{\bf E}_{(2)}\stackrel{\theta}{\to}({\bf
E}\otimes\Omega^1_{\overline M}(\log D))_{(2)}
\stackrel{\theta}{\to}({\bf E}\otimes\Omega^2_{\overline M}(\log
D))_{(2)}\stackrel{\theta}{\to}\cdots,
\]
which does not depend on the two metrics, but is defined
algebraically.

${\bf H_{\mathbb{C}}}$ has a canonical extension $\overline{\bf
H}_{\mathbb{C}}$ on $\overline{M} = {\triangle}^2$ across the
singularity as follows. As a sheaf, the germs of sections of
$\overline{\bf H}_{\mathbb{C}}$ at the singularity are generated
by the elements of $\Gamma (M, {\bf H_{\mathbb{C}}})$ which are of
the form
\[
{\tilde v} =
{\text{exp}}({\frac{1}{2\pi\sqrt{-1}}}(N_1{\text{log}}t_1 +
N_2{\text{log}}t_2))v,
\]
where $v$ is a multivalued flat section of ${\bf H}_{\mathbb{C}}$.
Again fixing an ordering $(N_1, N_2)$ and a point $p\in M$, as in
\S2, one has two monodromy weight filtrations $\{W^{\bf 1}_k\}$
and $\{W^{\bf 2}_s\}$ defined over $\mathbb{Q}$ corresponding to
$N_1$ and $c_1N_1+c_2N_2, c_1>0, c_2>0$ on $({\bf
H_{\mathbb{C}}})_p$ respectively, and they naturally determine two
filtrations of ${\bf H_{\mathbb{C}}}$ by locally constant systems,
denoted by $\{{\bf W}^{\bf 1}_k\}$ and $\{{\bf W}^{\bf 2}_s\}$
respectively; moreover, one also has the canonical extensions
$\{\overline{\bf W}^{\bf 1}_k\}$ and $\{\overline{\bf W}^{\bf
2}_s\}$ of $\{{\bf W}^{\bf 1}_k\}$ and $\{{\bf W}^{\bf 2}_s\}$
across the singularity, which form two filtrations of
$\overline{\bf H}_{\mathbb{C}}$ respectively.

\begin{prop} If restricting everything to $D_{\epsilon}$ for $\epsilon>0$, one
has the following equivalences.
\begin{eqnarray*}
&&\Omega^{0}({\bf H}_{\mathbb{C}})_{(2)} =
t_1t_2\overline{\bf H}_{\mathbb{C}} + t_1\bigcup_{l_2-l_1\le 0}
\overline{\bf W}^{\bf 1}_{l_1}\cap\overline{\bf W}^{\bf 2}_{l_2} +
t_2\overline{\bf W}^{\bf 1}_{0} + \bigcup_{l_1\le 0, l_2\le
l_1}\overline{\bf W}^{\bf 1}_{l_1}\cap
\overline{\bf W}^{\bf 2}_{l_2}; \\
&&\Omega^{1}({\bf H}_{\mathbb{C}})_{(2)} =    \\
&&{\frac{dt_1}{t_1}}\otimes \big(t_1t_2\overline{\bf
H}_{\mathbb{C}} + t_1\bigcup_{l_2-l_1\le 0}\overline{\bf W}^{\bf
1}_{l_1}\cap \overline{\bf W}^{\bf 2}_{l_2} + t_2\overline{\bf
W}^{\bf 1}_{-2}
+ \bigcup_{l_1\le -2, l_2\le l_1}\overline{\bf W}^{\bf 1}_{l_1}
\cap\overline{\bf W}^{\bf 2}_{l_2}\big) +     \\
&&+ {\frac{dt_2}{t_2}}\otimes \big(t_1t_2\overline{\bf
H}_{\mathbb{C}} + t_1\bigcup_{l_2\le l_1-2}\overline{\bf W}^{\bf
1}_{l_1}\cap \overline{\bf W}^{\bf 2}_{l_2} + t_2\overline{\bf
W}^{\bf 1}_{0} +
\bigcup_{l_1\le 0, l_2\le l_1-2}\overline{\bf W}^{\bf 1}_{l_1}
\cap\overline{\bf W}^{\bf 2}_{l_2}\big);  \\
&&\Omega^{2}({\bf H}_{\mathbb{C}})_{(2)} =
{\frac{dt_1}{t_1}}\wedge{\frac{dt_2}{t_2}}\otimes   \\
&&\otimes\big( t_1t_2\overline{\bf H}_{\mathbb{C}} +
t_1\bigcup_{l_2\le l_1-2}\overline{\bf W}^{\bf
1}_{l_1}\cap\overline{\bf W}^{\bf 2}_{l_2} + t_2\overline{\bf
W}^{\bf 1}_{-2} + \bigcup_{l_1\le -2, l_2\le l_1-2}\overline{\bf
W}^{\bf 1}_{l_1}\cap\overline{\bf W}^{\bf 2}_{l_2} \big).
\end{eqnarray*}
Here, ${\tilde v}\in\overline{\bf W}^{\bf 1}_{l_1}\cap
\overline{\bf W}^{\bf 2}_{l_2}$ for ${l_2-l_1\le 0}$ implies that
$v$ has nontrivial projections to both $Gr^{{\bf W}^{\bf
1}_*}_{l_1}$ and $Gr^{{\bf W}^{\bf 2}_*}_{l_2}$; the other terms
have the same explanation. Similarly, one can consider the other
ordering $(N_2, N_1)$ and restrict $\Omega^{i}({\bf
H}_{\mathbb{C}})_{(2)}$, $i=0, 1, 2$, to $D'_{\epsilon}$ to get
the same characterization of the sheaves as above.
\end{prop}
{\bf Remarks:} 1) As $\epsilon$ is sufficiently small,
$D_{\epsilon}\cup D'_{\epsilon}$ covers a neighborhood of the
singularity, so the above proposition actually gives a description
near the singularity of the sheaves $\Omega^{i}({\bf
H}_{\mathbb{C}})_{(2)}$, $i=0, 1, 2$, in terms of the weight
filtrations of $N_1, N_2, N_1+N_2$. 2) From the formulae above, we
also know that the sheaves $\Omega^{i}({\bf
H}_{\mathbb{C}})_{(2)}$ are contained in $j_*{\bf
H}_{\mathbb{C}}\otimes\Omega^r_{\overline M}({\text{log}}D)$,
$i=0, 1, 2$. 3) Although the proposition and the following proof
are presented for the case of dimension $2$, they still work for
general dimension.

\vskip .3cm \noindent {\bf Proof.} We will first show that one can
choose a suitable (flat sections) basis $\{v_i\}$ of ${\bf
H}_{\mathbb{C}}$ near the singularity such that it flags ${\bf
H}_{\mathbb{C}}$ according to the weight filtrations ${\bf
W}_*^{\bf 1}$ and ${\bf W}^{\bf 2}_*$ and $\{\tilde{v}_i\}$ is
$L^2$-adapted in the sense of \S2 on any domain of the form
$D_{\epsilon}$, $\epsilon >0$. We will use freely the notations
and arguments in \S2. By the nilpotent orbit theorem and the
Proposition 2, we know that $(W^{\bf 2}[-m], \psi(0))$ is a
(polarized) mixed Hodge structure; again by Proposition 1, there
exists a unique $\delta\in {\mathfrak{gl}}(H_{\mathbb{R}})$ such
that $(W^{\bf 2}[-m], F=\exp(-\sqrt{-1}\delta)\psi(0))$ is a
(polarized) mixed Hodge structure splitting over $\mathbb{R}$,
where $H_{\mathbb{R}}$ can be considered as the real structure of
the standard fibre $H_{\mathbb{C}}$ of ${\bf H}_{\mathbb{C}}$ as
lifted to the universal covering of $M$. Furthermore, there exits
an element ${\tilde g}(\infty)\in G_{\mathbb{C}}$ such that
$\log{\tilde g}(\infty)\in {\text{Ker}}({\text{ad}}(N_j))$ and
$(W^{\bf 2}[-m], {\tilde F}_{\bf 2} = {\tilde g}(\infty)F)$ is
canonically attached to $(W^{\bf 2}[-m], \psi(0))$. A simple
computation shows that the nilpotent orbit $\theta({\bf z})$
corresponding to the variation can be expressed as
\begin{eqnarray*}
\theta&=&\exp(\sum z_jN_j)\cdot\psi(0)\\
&=&\exp(\sum x_jN_j)\cdot
e^{-1}\cdot\exp(\sqrt{-1}{\text{Ad}}(e)(\sum y_jN_j+\delta))\cdot
eF,
\end{eqnarray*}
where $t_j=e^{2\pi\sqrt{-1}z_j}$, $z_j=x_j+\sqrt{-1}y_j$, and
$e=\exp({\frac 1 2}\log y_1{\tilde Y}_1+{\frac 1 2}\log y_2{\tilde
Y}_2) =\exp({\frac 1 2}\log s_1{\tilde Y}_{\bf 1}+{\frac 1 2}\log
s_2{\tilde Y}_{\bf 2})$, here $s_1={\frac{y_1}{y_2}}$ and
$s_2=y_2$, while ${\text{Ad}}$ is the adjoint representation of
$G_{\mathbb{R}}$. By the $SL_2$-orbit theorem (i),
$\exp(\sqrt{-1}{\tilde N}^-_{\bf 2})\cdot{\tilde F}_{\bf
2}=\exp(\sqrt{-1}N_1)\cdot{\tilde F}_{\bf 1}\in {\bf D}$. So, one
has
$$
F=\exp(-\sqrt{-1}{\tilde N}^-_{\bf 2})\cdot{\tilde
g}^{-1}(\infty)\cdot\exp(\sqrt{-1}N_1)\cdot{\tilde F}_{\bf 1}.
$$
Set $F_0={\tilde
g}^{-1}(\infty)\cdot\exp(\sqrt{-1}N_1)\cdot{\tilde F}_{\bf 1}\in
{\bf D}$ and ${\tilde p}=\exp(\sqrt{-1}{\text{Ad}}(e)(\sum
y_jN_j+\delta))$. By the definitions of $F$ and $e$, it is not
difficult to see that $eF=F$. So, one has
\[
\theta = \exp(\sum x_jN_j)\cdot e^{-1}\cdot{\tilde
p}\cdot\exp(-\sqrt{-1}{\tilde N}^-_{\bf 2})\cdot F_0.
\]
Thus, by the nilpotent orbit theorem again, there exists a
$G_{\mathbb{R}}$-valued function $h({\bf z})$ on ${\mathbb{C}}^2$
such that $h({\bf z})$ goes to the identity as ${\text{Im}}({\bf
z})\to\infty$ and
\[
F_0 = h({\bf z})\cdot\exp(\sqrt{-1}{\tilde N}^-_{\bf
2})\cdot{\tilde p}^{-1}\cdot e\cdot\exp(-\sum
x_jN_j)\cdot{\tilde\phi}({\bf z}).
\]
On the other hand, the argument in Lemma 5.12 of \cite{cks2} also
tells us that
\[
{\tilde p}\to\exp(\sqrt{-1}{\tilde N}^-_{\bf 2}), ~~{\text{as}}~
s_1, s_2\to \infty.
\]
More precisely, $\tilde{p}$ and its inverse are polynomials in
$\{s_1^{-{\frac 1 2}}, s_2^{-{\frac 1 2}}\}$ with constant
coefficients and the constant terms being $\exp(\sqrt{-1}{\tilde
N}^-_{\bf 2})$ and $\exp(-\sqrt{-1}{\tilde N}^-_{\bf 2})$
respectively; moreover, $\exp(\sqrt{-1}{\tilde N}^-_{\bf
2})\cdot{\tilde p}^{-1}$ can be expressed as
$\exp(-\sqrt{-1}\delta(s_2))\cdot\exp(-\sqrt{-1}Z(s_1))$, where
$\delta(s_2)$ and $Z(s_1)$ are some polynomials in $s_2^{-1}$ and
$s_1^{-1}$ respectively with coefficients in $\mathfrak{g}_{\bf
R}$. Denote $h({\bf z})\cdot\exp(\sqrt{-1}{\tilde N}^-_{\bf
2})\cdot{\tilde p}^{-1}$ by ${\tilde h}({\bf z})$, which, by the
behavior of $h(\bf z)$ as ${\text{Im}}({\bf z})\to\infty$, goes to
the identity as $s_1, s_2\to\infty$. (Note that $s_1,
s_2\to\infty$ implies ${\text{Im}}({\bf z})\to\infty$, but not the
converse. That's why we take the domains $D_{\epsilon}$ in the
estimate of the Hodge norm.) Since the image of $\tilde{\phi}$ for
enough large ${\text{Im}}(\bf z)$ and $F_0$ lie in $\bf D$, so
${\tilde h}({\bf z})\in G_{\mathbb{R}}$ for large enough
${\text{Im}}(\bf z)$. Thus, we establish an isometry ${\tilde
h}({\bf z})\cdot e\cdot\exp(-\sum x_jN_j)$ between the fiber over
$\bf z$ and $F_0\in {\bf D}$ for large enough ${\text{Im}}(\bf
z)$, denoted by $L({\bf z})$.

Let $v\in W_{l_1}^{\bf 1}\cap W_{l_2}^{\bf 2}$ and assume that its
projections on $Gr^{W^{\bf 1}_*}_{l_1}$ and $Gr^{W^{\bf
2}_*}_{l_2}$ are nontrivial, i.e., the components of $v$ with
eigenvalues $l_1$ and $l_2$ w.r.t. ${\tilde Y}_{\bf 1}$ and
${\tilde Y}_{\bf 2}$ respectively are nontrivial. (Note that
$W^{\bf j}_{l_j}$ is characterized by the property that $W^{\bf
j}_{l_j}$ is linearly spanned by all eigenvectors of ${\tilde
Y}_{\bf j}$ of eigenvalue not greater than $l_j$. Afterwards,
denote the eigenvector space ${\tilde Y}_{\bf j}$ for the
eigenvalue $l_j$ by $H_{l_j}({\tilde Y}_{\bf j})$) In the
following, we also denote by $v$ the multivalued flat section
generated by $v$, this will be clear from the context. Write
${\tilde v}$ as ${\tilde v}({\bf z})= \exp(\sum z_jN_j)v$ (i.e.
lift $\tilde v$ to the universal covering of $M$) and then
consider
\[
L({\bf z}){\tilde v}({\bf z}) = {\tilde h}({\bf
z})\cdot\exp({\frac 1 2}\log s_1{\tilde Y}_{\bf 1}+{\frac 1 2}\log
s_2{\tilde Y}_{\bf 2}) \cdot\exp(-\sum x_jN_j)\cdot{\tilde v}({\bf
z}).
\]
A direct computation shows
\begin{eqnarray*}
&&\exp({\frac 1 2}\log s_1{\tilde Y}_{\bf 1}+{\frac 1 2}\log
s_2{\tilde Y}_{\bf 2})
\cdot\exp(-\sum x_jN_j)\cdot{\tilde v}({\bf z})   \\
=&&\exp({\frac 1 2}\log s_1{\tilde Y}_{\bf 1}+{\frac 1 2}\log
s_2{\tilde Y}_{\bf 2})
\exp(\sqrt{-1}s_2(N_2+s_1N_1))\cdot v \\
=&&\sum_{k=0}^{m}{\frac{(\sqrt{-1})^k}{k!}}
\sum_{j=0}^kC_k^js_1^{k-j}s_2^k \exp({\frac 1 2}\log s_1{\tilde
Y}_{\bf 1}+{\frac 1 2} \log s_2{\tilde Y}_{\bf 2})\cdot
N_1^{k-j}N_2^jv.
\end{eqnarray*}
On the other hand, by the definitions of $W^{\bf 1}_*$ and $W^{\bf
2}_*$, it is not difficult to see that
\[
N_1^{k-j}N_2^j\cdot v\in W^{\bf 1}_{l_1-2k+2j}\cap W^{\bf
2}_{l_2-2k},
\]
and its $l_1-2k+2j$-eigenvector part and $l_2-2k$-eigenvector part
w.r.t. ${\tilde Y}_{\bf 1}$ and ${\tilde Y}_{\bf 2}$ respectively
are nontrivial. So,
\begin{eqnarray*}
&&\exp({\frac 1 2}\log s_1{\tilde Y}_{\bf 1}+{\frac 1 2}\log
s_2{\tilde Y}_{\bf 2})
\cdot\exp(-\sum x_jN_j)\cdot{\tilde v}({\bf z})   \\
=&&\sum_{k=0}^{m}{\frac{(\sqrt{-1})^k}{k!}}
\sum_{j=0}^kC_k^js_1^{k-j}s_2^k
s_1^{{\frac{l_1}2}-k+j}s_2^{{\frac{l_2}2}-k}N_1^{k-j}N_2^jv
+ ({\text{lower order terms}})  \\
=&&s_1^{\frac{l_1}2}s_2^{\frac{l_2}2}
\sum_{k=0}^m{\frac{(\sqrt{-1})^k}{k!}}(N_1+N_2)^kv
+ ({\text{lower order terms}}).
\end{eqnarray*}
Thus, if $\{v_l\}$ is a basis of $({\bf H}_{\mathbb{C}})_p$, which
flags $({\bf H}_{\mathbb{C}})_p$ according to the filtrations
$W_*^{\bf 1}$ and $W_*^{\bf 2}$ (this can always be done by the
characterization of $W_*^{\bf 1}$ and $W_*^{\bf 2}$ using
${\tilde{Y}}_{\bf j}, j=1, 2$), then
$\{s_1^{-\frac{l_1}2}s_2^{-\frac{l_2}2}L({\bf z}){\tilde v}_l({\bf
z})\}$ goes to a basis of $H_{\mathbb{C}}$ as $s_1, s_2\to
\infty$. So, by the definition of $L^2$-adapted basis and the
lemma 2 in \S2, $\{s_1^{-\frac{l_1}2}s_2^{-\frac{l_2}2}L({\bf
z}){\tilde v}_l({\bf z})\}$, and hence $\{{\tilde v}_l({\bf
z})\}$, is an $L^2$-adapted basis on any domain of the form
$$
{\tilde D}_{\epsilon}=\{{\bf z}=(z_1, z_2)\in {\mathbb{C}}^2~|~
{\frac{{\text{Im}}z_1}{{\text{Im}}z_2}}>\epsilon, {\text{Im}}z_2
>\epsilon\},~~\epsilon >0.
$$
Dropping $\{{\tilde v}_l({\bf z})\}$ and ${\tilde D}_{\epsilon}$
back to $M$, $\{\tilde{v}_l\}$ becomes an $L^2$-adapted basis on
$D_{\epsilon}$. This shows the claim in the beginning of the
proof.

We now turn to prove the proposition. We will prove only the
second identity, the others are similar. We also restrict all the
following discussion to $D_{\epsilon}$ for some fixed
$\epsilon>0$, if not specified. Due to all objects considered
being holomorphic under $\nabla^{0,1}$ and the argument above of
getting an $L^2$-adapted basis, we only need to show when the
elements of the form
$$
{\frac{dt_1}{t_1}}\otimes t_1^{n_1}t_2^{n_2}{\tilde v}_1 +
{\frac{dt_2}{t_2}}\otimes t_1^{n'_1}t_2^{n'_2}{\tilde v}_2
$$
are $L^2$ for $v_1\in {\bf W}^{\bf 1}_{l_1}\cap{\bf W}^{\bf
2}_{l_2}$ and $v_2\in {\bf W}^{\bf 1}_{l'_1}\cap{\bf W}^{\bf
2}_{l'_2}$ satisfying that the projections of both $v_1$ to
$Gr^{{\bf W}^{\bf 1}_*}_{l_1}$ and $Gr^{{\bf W}^{\bf 2}_*}_{l_2}$
and $v_2$ to $Gr^{{\bf W}^{\bf 1}_*}_{l'_1}$ and $Gr^{{\bf W}^{\bf
2}_*}_{l'_2}$ are nontrivial. By the estimate of the Hodge norm of
Theorem 4, we have
\[
|{\frac{dt_1}{t_1}}\otimes t_1^{n_1}t_2^{n_2}{\tilde v}_1|^2\sim
|\log|t_1||^2|t_1|^{2n_1}|t_2|^{2n_2}|\log|t_1||^{l_1}|\log|t_2||^{l_2-l_1}
\]
on $D_{\epsilon}$. So, ${\frac{dt_1}{t_1}}\otimes
t_1^{n_1}t_2^{n_2}{\tilde v}_1$ is $L^2$ iff the following
integral
\[
\int_{D_{\epsilon}}|\log|t_1||^2|t_1|^{2n_1}|t_2|^{2n_2}
|\log|t_1||^{l_1}|\log|t_2||^{l_2-l_1}dV_M
\]
is finite, where $dV_M$ is the volume element of M under the
Poincar\'e-like metric. A simple computation shows that the above
integral is finite iff
\begin{eqnarray*}
&& n_1\ge 1~~{\text{or}}~~ n_1 = 0~~ {\text{and}} ~~ l_1\le -2; \\
&& n_2\ge 1~~{\text{or}}~~ n_2 = 0~~ {\text{and}} ~~ l_2-l_1\le 0.
\end{eqnarray*}
Equivalently, ${\frac{dt_1}{t_1}}\otimes t_1^{n_1}t_2^{n_2}{\tilde
v}_1$ belongs to
\[
{\frac{dt_1}{t_1}}\otimes \big(t_1t_2\overline{\bf H}_{\mathbb{C}}
+ t_1\bigcup_{l_2-l_1\le 0}\overline{\bf W}^{\bf 1}_{l_1}\cap
\overline{\bf W}^{\bf 2}_{l_2} + t_2\overline{\bf W}^{\bf 1}_{-2}
+ \bigcup_{l_1\le -2, l_2\le l_1}\overline{\bf W}^{\bf
1}_{l_1}\cap\overline{\bf W}^{\bf 2}_{l_2}\big).
\]
From the discussion above, the above formula should be understood
as follows. ${\tilde v}_1\in\overline{\bf W}^{\bf 1}_{l_1}\cap
\overline{\bf W}^{\bf 2}_{l_2}$ for ${l_2-l_1\le 0}$ implies that
$v_1$ has nontrivial projections to both $Gr^{{\bf W}^{\bf
1}_*}_{l_1}$ and $Gr^{{\bf W}^{\bf 2}_*}_{l_2}$; the other terms
have a similar explanation.

\vskip .3cm The same computation shows that
${\frac{dt_2}{t_2}}\otimes t_1^{n'_1}t_2^{n'_2}{\tilde v}_2$ is
$L^2$ iff
\begin{eqnarray*}
&&n'_1\ge 1~~{\text{or}}~~n'_1 = 0~~{\text{and}}~~l'_1\le 0; \\
&&n'_2\ge 1~~{\text{or}}~~n'_2 = 0~~{\text{and}}~~l'_2-l'_1\le -2.
\end{eqnarray*}
Namely, ${\frac{dt_2}{t_2}}\otimes t_1^{n'_1}t_2^{n'_2}{\tilde
v}_2$ belongs to
\[
{\frac{dt_2}{t_2}}\otimes \big(t_1t_2\overline{\bf H}_{\mathbb{C}}
+ t_1\bigcup_{l'_2\le l'_1-2}\overline{\bf W}^{\bf 1}_{l'_1}\cap
\overline{\bf W}^{\bf 2}_{l'_2} + t_2\overline{\bf W}^{\bf 1}_{0}
+ \bigcup_{l'_1\le 0, l'_2\le l'_1-2}\overline{\bf W}^{\bf
1}_{l'_1}\cap\overline{\bf W}^{\bf 2}_{l'_2}\big).
\]
Clearly, the terms of the formula above should have the same
explanation as above. Thus, we finish the proof of the
proposition.~~~~~~~~~~~~~~~~~~~~~~~~~~~~~~~~~$\Box$

\vskip .3cm In the following, we turn to the case of general
dimension. Let
\[
{\bf E} = \oplus_{p=0}^{m}{\bf E}^p, ~~\theta=\sum_p\theta^p
\]
be the Higgs bundle corresponding to the variation as defined in
the introduction. From \S3, we know $\theta$ is an $L^2$-bounded
holomorphic $1$-form operator and
$$
F^{p+1}\Omega^r({\bf H}_{\mathbb{C}})_{(2)}/F^{p}\Omega^r({\bf
H}_{\mathbb{C}})_{(2)}\cong(\Omega^r\otimes{\bf E}^{p-r})_{(2)},
$$
where $F^{p}\Omega^r({\bf H}_{\mathbb{C}})_{(2)}$ is the
filtration of $\Omega^r({\bf H}_{\mathbb{C}})_{(2)}$ according to
$\{{\bf F}^p\}$ (cf. \S3). On the other hand, by the previous
proposition, $\Omega^r({\bf H}_{\mathbb{C}})_{(2)}\subset j_*{\bf
H}_{\mathbb{C}}\otimes\Omega^r_{\overline M}({\text{log}}D)$ and
can be described by just using the weight filtrations of $N_1,
N_2, N_1+N_2$, so $(\Omega^r\otimes{\bf E}^{p-r})_{(2)}$ is
actually the set of local $L^2$-sections of $j_*{\bf
E}^{p-r}\otimes\Omega^r_{\overline M}({\text{log}}D)$ and can
essentially be described in terms of the weight filtrations of
$N_1, N_2, N_1+N_2$. Omitting $j_*$, we from now on denote
$(\Omega^r\otimes{\bf E}^{p-r})_{(2)}$ by $({\bf
E}^{p-r}\otimes\Omega^r_{\overline M}({\text{log}}D))_{(2)}$,
$\oplus_p({\bf E}^{p-r}\otimes\Omega^r_{\overline
M}({\text{log}}D))_{(2)}$ by $({\bf E}\otimes\Omega^r_{\overline
M}(\log D))_{(2)}$. The $L^2$ boundedness of $\theta$ and
$\theta\wedge\theta=0$ tell us that the following sequence
\[
(*)~~~~~~~~~~~~~{\bf E}_{(2)}\stackrel{\theta}{\to}({\bf
E}\otimes\Omega^1_{\overline M}(\log D))_{(2)}
\stackrel{\theta}{\to}({\bf E}\otimes\Omega^2_{\overline M}(\log
D))_{(2)}\stackrel{\theta}{\to}\cdots
\]
is a complex, called the $L^2$ holomorphic Dolbeault complex.

\vskip .3cm\noindent {\bf Remark:} By the above argument and
Proposition 6, one can actually write down $({\bf
E}\otimes\Omega^r_{\overline M}(\log D))_{(2)}$ explicitly in
terms of the basis constructed in $\S 3$. Let
$\{\sigma_{k,l}|k+l=s\}$ be the basis of ${\bf E}^s$ constructed
in $\S3$. By the construction, we know that $\sigma_{k,l}\in
{\overline{\bf W}}^{\bf 1}_{2k}\cap{\overline{\bf W}}^{\bf
2}_{2(k+l)}$ (under a suitable identification, cf. Corollary 1) up
to some higher order terms. So, using the estimate of the Hodge
norm, we can show which holomorphic $r$-forms valued in the line
bundle generated by $\sigma_{k,l}$ are $L^2$.

\section{$L^2$ $\overline\partial$-Poincar\'e lemma}

Let $(M, {\bf H}_{\mathbb{C}}={\bf
H}_{\mathbb{Z}}\otimes{\mathbb{C}}, \nabla = {\nabla}^{1,0} +
{\nabla}^{0,1} , {\bf F}=\{{\bf F}^p\}, {\bf S})$ be a variation
of polarized Hodge structure with weight $m$ defined over
$\mathbb{Q}$. In the previous section, we defined the $L^2$
holomorphic Dolbeault complex $\{({\bf
E}\otimes\Omega^r_{\overline M}(\log D))_{(2)}, \theta\}$ of the
variation on $\overline M$, which is independent of both metrics
and essentially determined by the logarithmic monodromies. For
uniformity of notations of this section, we from now on denote
$({\bf E}^{p-r}\otimes\Omega^r_{\overline M}(\log D))_{(2)}$ by
${\text Gr}_{F}^p{\Omega}^{r}({\bf H}_{\mathbb{C}})_{(2)}$, which
is a piece of $\{({\bf E}\otimes\Omega^r_{\overline M}(\log
D))_{(2)}$ as seen at the end of \S4.

By the infinitesimal period relation of ${\nabla}^{1,0}$ and the
definition and the boundedness of $\theta$, we have
\[
\theta(\{{\text Gr}_{F}^p{\Omega}^{r}({\bf
H}_{\mathbb{C}})_{(2)})\subset \{{\text
Gr}_{F}^p{\Omega}^{r+1}({\bf H}_{\mathbb{C}})_{(2)}.
\]
(More precisely, we should restrict to a piece of $\theta$ above.)
Thus, we obtain a holomorphic Dolbeault subcomplex of $\{({\bf
E}\otimes\Omega^{\cdot}_{\overline M}(\log D))_{(2)}, \theta\}$ on
$\overline M$
\[
\{{\text Gr}_{F}^p{\Omega}^{\cdot}({\bf H}_{\mathbb{C}})_{(2)},
\theta\}.
\]
The main purpose of this section is then to show that the
hypercohomology ${\mathbb{H}}^*({\overline M}, \{{\text
Gr}_{F}^p{\Omega}^{.}({\bf H}_{\mathbb{C}})_{(2)}, \theta\})$ is
equal to a certain $L^2$-cohomology on $\overline M$. Thus, it is
computable from a certain complex of $L^2$-differential forms.

To this end, we first define some fine sheaves on $\overline M$.
Let $F^pA^k({\bf H}_{\mathbb{C}})_{(2)} =
\oplus_{r+s=k}(A^{r,s}\otimes {\bf F}^{p-r})_{(2)}$ and $D''_1 =
{\overline{\partial}} + {\nabla}^{1,0}$ for $p\ge 0$. Here,
$A^{r,s}$ is the sheaf of germs of local forms of type $(r, s)$
(not necessarily smooth) on $\overline M$ and $(A^{r,s}\otimes
{\bf F}^{p-r})_{(2)}$ is the sheaf of germs of local $L^2$ ${\bf
F}^{p-r}$-valued forms $\phi$ of type $(r,s)$ on $\overline M$ for
which ${\overline
\partial}\phi$ are $L^2$ in the weak sense.
Here the action of $\overline{\partial}$ is defined as follows:
Let $\phi$ be a form of type $(r, s)$ and $v$ a holomorphic
section of ${\bf F}^p$, then
$$
\overline\partial(\phi\otimes v)=\overline\partial\phi\otimes v.
$$
Note that $\{F^{p}A^k({\bf H}_{\mathbb{C}})_{(2)}, p\ge 0\}$ is in
general not a filtration of $A^k({\bf H}_{\mathbb{C}})_{(2)}$ (for
the definition, cf. \S1). It is clear that $D''_1(F^pA^k({\bf
H}_{\mathbb{C}})_{(2)})\subset F^pA^{k+1}({\bf
H}_{\mathbb{C}})_{(2)}$ by the boundedness of $\nabla^{1,0}$
(\S3). By the Hodge filtration $\{{\bf F}^p\}$, take the
successive quotients
\[
[{\text Gr}_F^pA^k({\bf H}_{\mathbb{C}})]_{(2)} := F^pA^k({\bf
H}_{\mathbb{C}})_{(2)}/F^{p+1}A^k({\bf H}_{\mathbb{C}})_{(2)},
\]
which, by the $L^2$-adaptedness theorem (which is clearly true for
the differentiable case by the proof), can be identified with
$\oplus_{r+s=k}(A^{r,s}\otimes {\bf E}^{p-r})_{(2)}$. Here,
$(A^{r,s}\otimes {\bf E}^{p-r})_{(2)}$ is the sheaf on $\overline
M$ of germs of local $L^2$ ${\bf E}^{p-r}$-valued forms $\phi$ of
type $(r,s)$ for which ${\overline \partial}\phi$ are $L^2$ in the
weak sense. Denote the induced map of $D_1''$ by $D''$, which is
actually ${\overline{\partial}} + \theta$ and clearly satisfies
$(D'')=0$. We now obtain a complex of fine sheaves on $\overline
M$
\[
\{[{\text Gr}_F^pA^{\cdot}({\bf H}_{\mathbb{C}})]_{(2)}, D''\},
\]
for $p\ge 0$ and the holomorphic Dolbeault subcomplex $\{{\text
Gr}_{F}^p{\Omega}^{\cdot}({\bf H}_{\mathbb{C}})_{(2)}, \theta\}$
is obviously its subcomplex. Similar to the holomorphic Dolbeault
subcomplex above, $([{\text Gr}_F^pA^{\cdot}({\bf
H}_{\mathbb{C}})]_{(2)}, D'')$ can also be considered as a piece
of a larger complex of fine sheaves: Let ${\bf E}$ be the Higgs
bundle of the variation with the Higgs field $\theta$, denote by
$[{\text Gr}_F^{*}A^{k}({\bf H}_{\mathbb{C}})]_{(2)}$ the sheaf of
germs of local $L^2$ $k$-forms $\phi$ (not necessarily smooth)
with values in $\bf E$ on $\overline M$, for which
$\overline\partial\phi$ are also $L^2$ in the weak sense. It is
clear that $\{[{\text Gr}_F^{*}A^{\cdot}({\bf
H}_{\mathbb{C}})]_{(2)}, D''\}$ is a complex of fine sheaves and
$\{[{\text Gr}_F^pA^{\cdot}({\bf H}_{\mathbb{C}})]_{(2)}, D''\}$
is a piece of it.

The purpose of the remaining part in this section will be to show
the following
\begin{thm}
The holomorphic Dolbeault complex $\{{\text
Gr}_{F}^p{\Omega}^{\cdot}({\bf H}_{\mathbb{C}})_{(2)}, \theta\}$
is quasi-isomorphic to the complex $\{[{\text
Gr}_F^pA^{\cdot}({\bf H}_{\mathbb{C}})]_{(2)}, D''\}$ under the
inclusion map for $p\ge 0$.
\end{thm}

Concerning the definition of {\it quasi-isomorphism} and relative
notations, one can refer to \cite{gh} or the appendix in this
paper. Summing up for $p\ge 0$, the above theorem implies

\begin{cor}
The holomorphic Dolbeault complex (*) $\{({\bf
E}\otimes\Omega^{\cdot}_{\overline M}(\log D))_{(2)}, \theta\}$ is
quasi-isomorphic to the complex (**) $\{[{\text
Gr}_F^*A^{\cdot}({\bf H}_{\mathbb{C}})]_{(2)}, D''\}$ under the
inclusion map.
\end{cor}

The standard result (cf. \cite{gri}) then tells us that the
hypercohomologies of (*) and (**) are isomorphic:
\[
{\mathbb{H}}^*({\overline M}, \{({\bf
E}\otimes\Omega^{\cdot}_{\overline M}(\log D))_{(2)}, \theta\})
\simeq {\mathbb{H}}^*({\overline M}, \{[{\text
Gr}_F^*A^{\cdot}({\bf H}_{\mathbb{C}})]_{(2)}, D''\});
\]
on the other hand, since $\{[{\text Gr}_F^*A^{\cdot}({\bf
H}_{\mathbb{C}})]_{(2)}, D''\}$ is a complex of fine sheaves, so
its hypercohomology is just the cohomology $H^*(\{\Gamma
({\overline M}, [{\text Gr}_F^*A^{\cdot}({\bf
H}_{\mathbb{C}})]_{(2)}), D''\})$ of the corresponding complex of
global sections of the sheaves (cf. \cite{gri} or the appendix).
Thus we obtain Theorem A
\[
{\mathbb{H}}^*({\overline M}, \{({\bf
E}\otimes\Omega^{\cdot}_{\overline M}(\log D))_{(2)}, \theta\})
\simeq H^*(\{\Gamma ({\overline M}, [{\text Gr}_F^*A^{\cdot}({\bf
H}_{\mathbb{C}})]_{(2)}), D''\}).
\]

In order to prove that the inclusion map is a quasi-isomorphism,
we need to prove the $\overline{\partial}$-Poincar\'e lemma under
the present situation. More concretely, we need to prove:

\vskip .3cm \noindent {\bf $L^2$ $\overline{\partial}$-Poincar\'e
Lemma} valued in a hermitian vector bundle: { \it Let
$\tilde{\psi}^{r,k-r}\in (A^{r,k-r}\otimes {\bf E}^{p-r})_{(2)}$
($\subset [{\text Gr}_F^pA^k({\bf H}_{\mathbb{C}})]_{(2)}$) be a
local section near the divisor with
${\overline\partial}\tilde{\psi}^{r,k-r} = 0$. Then there exists a
local section ${\psi}^{r, k-r-1}\in (A^{r,k-r-1}\otimes {\bf
E}^{p-r})_{(2)}$ near the divisor (possibly on a smaller defining
domain) with}
\[
\overline{\partial}{\psi}^{r, k-r-1} = \tilde{\psi}^{r,k-r}.
\]
Note that for a neighborhood not containing the divisor, the $L^2$
$\overline{\partial}$-Poincar\'e lemma is just the classical
$\overline{\partial}$-Poincar\'e lemma. As seen in \S3, near the
divisor, ${\bf E}^s$ is generated by some holomorphic sections
$\{{\sigma}_{k,l}|k+l=s\}$, which are linearly independent
everywhere and satisfy the norm estimates of the following form
(cf. Prop. 4)
\[
{\Vert{\sigma}_{k,l}\Vert}^2\sim
(-{\text{log}}|t_1|)^{2k-m}(-{\text{log}}|t_2|)^{2l-n}.
\]
So, we can reduce the above $L^2$ $\overline{\partial}$-Poincar\'e
lemma to the following

\vskip .3cm \noindent {\bf $L^2$ $\overline{\partial}$-Poincar\'e
Lemma} valued in a line bundle: {\it Let $L$ be a holomorphic line
bundle on $(\triangle^*)^n$ with a generating section $\sigma$,
and with a Hermitian metric satisfying, for some fixed integers
$k, l$,
\[
{\Vert\sigma\Vert}^2\sim (-{\text{log}}|t_1|)^{k}(-{\text{log}}|t_2|)^{l}
\]
near the divisor. Then for each $L$-valued $L^2$ $(r,s)$-form
$\phi\otimes\sigma$ near the divisor with
${\overline\partial}(\phi\otimes\sigma) = 0$, there exists an
$L$-valued $L^2$ $(r,s-1)$-form $\psi\otimes\sigma$ (possibly on a
smaller defining domain) with
${\overline\partial}(\psi\otimes\sigma) = \phi\otimes\sigma$,
i.e., ${\overline\partial}\psi = \phi$.}

\vskip .3cm Very unfortunately, we can prove the above $L^2$
$\overline{\partial}$-Poincar\'e lemma only for $k\neq k(r)$ and
$l\neq l(r)$, where $k(r)$ and $l(r)$ are some integers depending
on $r$ (for details, see the following remark). It is however
lucky enough that the proof of quasi-isomorphism can be obtained
just from the above $L^2$ $\overline\partial$-Poincar\'e lemma for
$k\neq k(r)$ and $l\neq l(r)$ and the properties of $\theta$
established in \S3. We now begin to prove the $L^2$
$\overline{\partial}$-Poincar\'e lemma above for $k\neq k(r)$ and
$l\neq l(r)$. Actually we only need to prove it for $r=0$, while
at the time $k(0)=l(0)=1$ (for details, see the following remark).
For emphasis, we restate it as
\begin{prop}
Let $L$ be a holomorphic line bundle on $(\triangle^*)^n$ with  a
generating section $\sigma$, and with a Hermitian metric
satisfying
\[
{\Vert\sigma\Vert}^2\sim (-{\text{log}}|t_1|)^{k}
(-{\text{log}}|t_2|)^{l} ~~{\text{for}} ~k, ~l\neq 1
\]
near the divisor. Then for each $L$-valued $L^2$ $(0,s)$-form
($s\ge 1$) $\phi\otimes\sigma$ with
${\overline\partial}(\phi\otimes\sigma) = 0$ near the divisor,
there exists an $L$-valued $L^2$ $(0,s-1)$-form
$\psi\otimes\sigma$ (possibly on a smaller defining domain) with
${\overline\partial}(\psi\otimes\sigma) = \phi\otimes\sigma$,
i.e., ${\overline\partial}\psi = \phi$.
\end{prop}

\noindent {\bf Proof.} For the sake of simplicity, we again assume
$n=2$. In order to prove the proposition, we actually need to
prove only the following assertion: {\it Let $\phi$ be a form of
type $(0, s)$ on $({\triangle^*})^2$ with
\[
\Vert\phi\otimes\sigma\Vert^2_{(2)}=:
\iint_{(\triangle^*)^2}\Vert\phi\Vert^2 \Vert\sigma\Vert^2
{\frac{dt_1\wedge d{\overline
t}_1}{|t_1|^2(-{\text{log}}|t_1|)^2}} {\frac{dt_2\wedge
d{\overline t}_2}{|t_2|^2(-{\text{log}}|t_2|)^2}} < \infty
\]
and $\overline{\partial}{\phi}=0$. Then there exists a form $\psi$
of type $(0, s-1)$ on $(\triangle^*)^2$ {\footnote {In general the
defining domain of $\psi$ will be smaller than that of $\phi$.
Afterwards, we will not specify this, since it is clear from the
context.}} with
\[
\Vert\psi\otimes\sigma\Vert^2_{(2)}=:
{\iint_{(\triangle^*)^2}}\Vert\psi\Vert^2 \Vert\sigma\Vert^2
{\frac{dt_1\wedge d{\overline
t}_1}{|t_1|^2(-{\text{log}}|t_1|)^2}} {\frac{dt_2\wedge
d{\overline t}_2}{|t_2|^2(-{\text{log}}|t_2|)^2}} < \infty
\]
and $\overline{\partial}\psi = \phi$. Here, the norms
$\Vert\phi\Vert$ and $\Vert\psi\Vert$ are measured under the
Poincar\'e-like metric.} We will use Fourier series to construct a
formal solution $\psi$ to $\overline{\partial}\psi = \phi$ for the
given $\phi$ satisfying, for some constant $C>0$,
\[
\Vert\psi\otimes\sigma\Vert^2_{(2)}\le C\Vert\phi\otimes\sigma\Vert^2_{(2)},
\]
which shows that $\psi$ is indeed a real solution to the
$\overline\partial$-problem. Using the technique of approximation,
one can furthermore assume that $\phi$ has compact support in
$(\triangle^*)^2$.

\noindent {\bf The proof of the assertion:} Using polar
coordinates ${(r_i, {\theta}_i)}$ for coordinate components $t_i$
with $|t_i|=r_i$, $i=1, 2$, the Poincar\'e-like metric can then be
rewritten as (neglecting the higher term which essentially plays
no role)
\[
{\frac{dr_1^2 + r_1^2d{\theta}_1^2}{r_1^2{\text{log}}^2r_1}}
+ {\frac{dr_2^2 + r_2^2d{\theta}_2^2}{r_2^2{\text{log}}^2r_2}}.
\]
Since we assume that the dimension of the base manifold is $2$, we
have only two cases to consider: 1) $\phi$ is a $(0, 1)$-form; 2)
$\phi$ is a $(0, 2)$-form.

\noindent {\bf Case 1:} Write $\phi$ as $f^1d{\overline t}_1 +
f^2d{\overline t}_2$. The conditions on $\phi$ are then equivalent
to the following integrability condition
\[
{\frac{{\partial}{f^1}}{\partial{\overline{t}_2}}} =
{\frac{{\partial}{f^2}}{\partial{\overline{t}_1}}}
\]
and
\begin{eqnarray*}
&&{\iint_{(\triangle^*)^2}} (|f^1|^2|t_2|^{-2}|\log|t_2||^{-2}
\Vert\sigma\Vert^2
{dt_1\wedge d{\overline t}_1}\wedge{dt_2\wedge d{\overline t}_2} +\\
&&{\iint_{(\triangle^*)^2}} (|f^2|^2|t_1|^{-2}|\log|t_1||^{-2}
\Vert\sigma\Vert^2 {dt_1\wedge d{\overline t}_1}\wedge{dt_2\wedge
d{\overline t}_2} < \infty;
\end{eqnarray*}
and we need to prove that there exists a function $u$ on
$(\triangle^*)^2$ satisfying
\[
{\iint_{(\triangle^*)^2}}|u|^2|t_1|^{-2}|t_2|^{-2}
(-{\text{log}}|t_1|)^{-2}(-{\text{log}}|t_2|)^{-2}
\Vert\sigma\Vert^2 {dt_1\wedge d{\overline t}_1}\wedge{dt_2\wedge
d{\overline t}_2} < \infty
\]
and
\[
{\frac{{\partial}{u}}{\partial{\overline{t}_1}}} = f^1;~~
{\frac{{\partial}{u}}{\partial{\overline{t}_2}}} = f^2.
\]
In order to do this, write $f^1$, $f^2$, and $u$ as $(r_1,
r_2)$-dependent Fourier series
\begin{eqnarray*}
&&f^1=\sum f^1_{m,n}(r_1, r_2)\exp(\sqrt{-1}m\theta_1+\sqrt{-1}n\theta_2), \\
&&f^2=\sum f^2_{m,n}(r_1, r_2)\exp(\sqrt{-1}m\theta_1+\sqrt{-1}n\theta_2), \\
&&u=\sum u_{m,n}(r_1,
r_2)\exp(\sqrt{-1}m\theta_1+\sqrt{-1}n\theta_2).
\end{eqnarray*}
Since ${\frac{\partial}{{\partial\overline t}_i}}= {\frac 1
2}e^{\sqrt{-1}\theta_i}[{\frac{\partial}{{\partial r}_i}}+
{\frac{\sqrt{-1}}{r_i}}{\frac{\partial}{{\partial\theta}_i}}],
i=1, 2$, the integrability condition
${\frac{{\partial}{f^1}}{\partial{\overline{t}_2}}} =
{\frac{{\partial}{f^2}}{\partial{\overline{t}_1}}}$ becomes
\[
{\frac{\partial{f^1_{m+1,n}}}{{\partial r}_2}}-{\frac{n}{r_2}}f^1_{m+1,n} =
{\frac{\partial{f^2_{m,n+1}}}{{\partial r}_1}}-{\frac{m}{r_1}}f^1_{m,n+1};
\]
while the equations
${\frac{{\partial}{u}}{\partial{\overline{t}_1}}} = f^1$,
${\frac{{\partial}{u}}{\partial{\overline{t}_2}}} = f^2$ become
\[
{\frac 1 2}[{\frac{\partial{u_{m,n}}}{{\partial
r}_1}}-{\frac{m}{r_1}}u_{m,n}] = f^1_{m+1,n},~~ {\frac
12}[{\frac{\partial{u_{m,n}}}{{\partial
r}_2}}-{\frac{n}{r_2}}u_{m,n}] = f^2_{m,n+1},
\]
for all $m, n\in{\mathbb{Z}}$; or
\begin{eqnarray*}
&&{\frac{\partial}{\partial{r_2}}}(r^{-m}_1r^{-n}_2f^1_{m+1,n})
= {\frac{\partial}{\partial{r_1}}}(r^{-m}_1r^{-n}_2f^1_{m,n+1}); \\
&&{\frac{\partial}{\partial{r_1}}}(r^{-m}_1r^{-n}_2u_{m,n}) =
2r_1^{-m}r_2^{-n}f^1_{m+1,n},~~
{\frac{\partial}{\partial{r_2}}}(r^{-m}_1r^{-n}_2u_{m,n}) =
2r_1^{-m}r_2^{-n}f^2_{m,n+1}.
\end{eqnarray*}
Put, for all $m, n\in{\mathbb{Z}}$,
\[
u_{m,n}(r_1, r_2) =
2r_1^mr_2^n\int_{L}\rho_1^{-m}\rho_2^{-n}f^1_{m+1,n}(\rho_1,
\rho_2){\text{d}}\rho_1 +
\rho_1^{-m}\rho_2^{-n}f^2_{m,n+1}(\rho_1, \rho_2){\text{d}}\rho_2,
\]
where the integral is a curve integral along the following
oriented curve ( for some positive constant $A<1$)
\begin{eqnarray*}
L =
\begin{cases}
(0, 0)\to (r_1, r_2)~~{\text{if}}~(m<0; n<0), (m<0; n=0, l>1), \\
~~~~~~~~~~~~~~~~~~~~~~~~~(m=0, k>1; n<0)~{\text{or}}~(m=0, k>1; n=0, l>1);\\
(0, A)\to (r_1, r_2)~~{\text{if}}~(m<0; n>0), (m<0; n=0, l>1), \\
~~~~~~~~~~~~~~~~~~~~~~~~~(m=0, k>1; n>0)~{\text{or}}~(m=0, k>1; n=0, l<1);\\
(A, 0)\to (r_1, r_2)~~{\text{if}}~(m>0; n<0), (m>0; n=0, l>1), \\
~~~~~~~~~~~~~~~~~~~~~~~~~(m=0, k<1; n<0)~{\text{or}}~(m=0, k<1; n=0, l>1);\\
(A, A)\to (r_1, r_2)~~{\text{if}}~(m>0; n>0), (m>0; n=0, l<1), \\
~~~~~~~~~~~~~~~~~~~~~~~~~(m=0, k<1; n>0)~{\text{or}}~(m=0, k<1; n=0, l<1),
\end{cases}
\end{eqnarray*}
and formally
\[
u=\sum u_{m,n}(r_1,
r_2)\exp(\sqrt{-1}m\theta_1+\sqrt{-1}n\theta_2).
\]
We will next show that
\begin{eqnarray*}
&&{\iint_{(\triangle^*)^2}}|u|^2 \Vert\sigma\Vert^2
{\frac{dt_1\wedge d{\overline
t}_1}{|t_1|^2(-{\text{log}}|t_1|)^2}}
{\frac{dt_2\wedge d{\overline t}_2}{|t_2|^2(-{\text{log}}|t_2|)^2}}\\
&\le& C{\iint_{(\triangle^*)^2}}\Vert\phi\Vert^2
\Vert\sigma\Vert^2 {\frac{dt_1\wedge d{\overline
t}_1}{|t_1|^2(-{\text{log}}|t_1|)^2}} {\frac{dt_2\wedge
d{\overline t}_2}{|t_2|^2(-{\text{log}}|t_2|)^2}}
\end{eqnarray*}
for some constant $C$, i.e. $\Vert u\sigma\Vert^2_{(2)}\le C\Vert
\phi\otimes\sigma\Vert_{(2)}^2$; in terms of polar coordinates and
the Fourier series of $u, f^1, f^2$ and using the asymptotic
behavior of $\Vert\sigma\Vert$, it is just
\begin{eqnarray*}
&&\sum_{m,n}\iint_{[0, A]^2}|u_{m,n}(r_1, r_2)|^2
r_1^{-1}r_2^{-1}|\log r_1|^{k-2}|\log r_2|^{l-2}{\text d}r_1{\text d}r_2\\
&\le&C\big\{\sum_{m,n}\iint_{[0, A]^2}|f^1_{m+1,n}(r_1, r_2)|^2
r_1r_2^{-1}|\log r_1|^{k}|\log r_2|^{l-2}{\text d}r_1{\text d}r_2\\
&&+\sum_{m,n}\iint_{[0, A]^2}|f^2_{m,n+1}(r_1, r_2)|^2
r_1^{-1}r_2|\log r_1|^{k-2}|\log r_2|^{l}{\text d}r_1{\text
d}r_2\big\}.
\end{eqnarray*}
If this is done, then $u$ is the required solution, i.e.
$\overline\partial u=\phi$.

We first do this for the oriented curves $L$ from $(0, 0)$ to
$(r_1, r_2)$. In this case, we choose a particular path which is
from $(0, 0)$ to $(0, r_2)$ (resp. $(r_1, 0)$), then to $(r_1,
r_2)$. (More precisely, the segment from $(0, 0)$ to $(0, r_2)$
(resp. $(r_1, 0)$) should be substituted by a path sufficiently
closed to the $r_2$-axis (resp. $r_1$-axis) from $(0, 0)$ to $(s,
r_2)$ (resp. $(r_1, s)$) for some sufficiently small $s>0$, but
since we assume that $\phi$ have compact support, this path plays
essentially no role.) $u_{m,n}$ can then be rewritten as
\begin{eqnarray*}
u_{m,n}(r_1, r_2) &=& 2r_1^m\int_0^{r_1}{\rho_1^{-m}}
f_{m+1,n}^1(\rho_1, r_2){\text d}{\rho_1}\\
({\text{resp.}}&=& 2r_2^n\int_0^{r_2}{\rho_2^{-n}}
f_{m,n+1}^2(r_1, {\rho}_2){\text d}{\rho_2}).
\end{eqnarray*}
In the following, without loss of generality, we always assume
that $\Vert\sigma\Vert^2=|\log|t_1||^{k}|\log|t_2||^{l}$. Letting
$m<0$ and $n<0$, one then has
\begin{eqnarray*}
&&\iint_{[0,A]^2} |u_{m,n}(r_1,r_2)|^2r_1^{-1}r_2^{-1}
|\log r_1|^{k-2}|\log r_2|^{l-2}dr_1dr_2\\
&=&\iint_{[0,A]^2}\big(2r_1^m\int_0^{r_1}{\rho_1^{-m}}f_{m+1,n}^1(\rho_1,
r_2){\text d}{\rho_1}\big)^2
r_1^{-1}r_2^{-1}|\log r_1|^{k-2}|\log r_2|^{l-2}dr_1dr_2\\
\end{eqnarray*}
\begin{eqnarray*}
&=&\iint_{[0,A]^2}4r_1^{2m}\big(\int_0^{r_1}\rho_1^{-2m}|f_{m+1,n}^1(\rho_1,
r_2)|^2d\rho_1\big)
\big(\int_0^{r_1}d\rho_1\big)r_1^{-1}r_2^{-1}|\log r_1|^{k-2}|\log r_2|^{l-2}dr_1dr_2\\
&=&\iint_{[0,A]^2}4r_1^{2m}\big(\int_0^{r_1}\rho_1^{-2m}|f_{m+1,n}^1(\rho_1,
r_2)|^2d\rho_1\big)
r_2^{-1}|\log r_1|^{k-2}|\log r_2|^{l-2}dr_1dr_2\\
&=&\int_0^A\big\{\int_0^A4\big(\int_0^{r_1}\rho_1^{-2m}
|f_{m+1,n}^1(\rho_1, r_2)|^2d\rho_1\big)
r_1^{2m}|\log r_1|^{k-2}dr_1\big\}r_2^{-1}|\log r_2|^{l-2}dr_2\\
&=&\int_0^A\big\{\lim_{\epsilon\to
0}\big(-4\int_r^A\rho_1^{2m}|\log\rho_1|^{k-2}d\rho_1\cdot
\int_0^{r_1}\rho_1^{-2m}|f_{m+1,n}^1(\rho_1, r_2)|^2d\rho_1|^A_{\epsilon}\big)\\
&&+4\int_0^Ar_1^{-2m}|f_{m+1,n}^1(r_1,
r_2)|^2\big(\int_{r_1}^A\rho_1^{2m}|\log\rho_1|^{k-2}d\rho_1\big)dr_1
\big\}r_2^{-1}|\log r_2|^{l-2}dr_2\\
&=&4\int_0^A\big\{\int_0^Ar_1^{-2m}|f_{m+1,n}^1(r_1,
r_2)|^2\big(\int_{r_1}^A\rho_1^{2m}|\log\rho_1|^{k-2}d\rho_1\big)dr_1
\big\}r_2^{-1}|\log r_2|^{l-2}dr_2\\
&\le&4\int_0^A\big\{\int_0^Ar_1^{-2m}|f_{m+1,n}^1(r_1,
r_2)|^2r_1^{2m+1}|\log r_1|^{k-2}\log{\frac{A}{r_1}}dr_1
\big\}r_2^{-1}|\log r_2|^{l-2}dr_2\\
&&({\text{Here, we used that}}~\rho^n|\log\rho|^k~
{\text{is a decreasing function in}}~\rho~ {\text{for}}~n<0.)  \\
&\sim&\int_0^A\big\{\int_0^Ar_1|f_{m+1,n}^1(r_1, r_2)|^2
|\log r_1|^{k-1}dr_1\big\}r_2^{-1}|\log r_2|^{l-2}dr_2\\
&\le&\int_0^A\big\{\int_0^A|f_{m+1,n}^1(r_1, r_2)|^2
|\log r_1|^{k}dr_1\big\}r_2^{-1}|\log r_2|^{l-2}dr_2\\
&=&\iint_{[0,A]^2}|f_{m+1,n}^1(r_1, r_2)|^2|\log r_1|^{k}|\log
r_2|^{l-2}r_2^{-1}dr_1dr_2;
\end{eqnarray*}
similarly, by using the other path mentioned above, one also has
\begin{eqnarray*}
&&\iint_{[0,A]^2} |u_{m,n}(r_1,r_2)|^2r_1^{-1}r_2^{-1}
|\log r_1|^{k-2}|\log r_2|^{l-2}dr_1dr_2\\
&\le&\iint_{[0,A]^2}|f_{m,n+1}^2(r_1, r_2)|^2|\log r_1|^{k-2}|\log
r_2|^{l}r_1^{-1}dr_1dr_2.
\end{eqnarray*}
Letting $m=0$, $k>1$ and $n=0$, $l>1$, one then has
\begin{eqnarray*}
u_{0,0}(r_1, r_2) &=& 2\int_0^{r_1}f_{1,0}^1(\rho_1, r_2){\text d}{\rho_1}\\
({\text{resp.}}&=& 2\int_0^{r_2}f_{0,1}^2(r_1, {\rho}_2){\text d}{\rho_2}),
\end{eqnarray*}
and
\begin{eqnarray*}
&&\iint_{[0,A]^2} |u_{0,0}(r_1,r_2)|^2r_1^{-1}r_2^{-1}
|\log r_1|^{k-2}|\log r_2|^{l-2}dr_1dr_2\\
&=&\iint_{[0,A]^2}\big(2\int_0^{r_1} f_{1,0}^1(\rho_1, r_2){\text
d}{\rho_1}\big)^2
r_1^{-1}r_2^{-1}|\log r_1|^{k-2}|\log r_2|^{l-2}dr_1dr_2\\
\end{eqnarray*}
\begin{eqnarray*}
&\le&4\int_0^A\big\{\int_0^A\big(\int_0^{r_1}|
f_{1,0}^1(\rho_1, r_2)|^2\rho_1|\log\rho_1|^{1+\eta}d\rho_1\big)\cdot\\
&&\big(\int_0^{r_1}{\rho_1}^{-1}|\log\rho_1|^{-1-\eta}d\rho_1\big)
r_1^{-1}|\log r_1|^{k-2}dr_1\big\}r_2^{-1}|\log r_2|^{l-2}dr_2\\
&&({\text{Here, $\eta$ is a positive constant $<k-1$}}.)\\
&\le&{\frac{4}{\eta}}\int_0^A\big\{\int_0^A
\big(\int_0^{r_1}|f_{1,0}^1(\rho_1, r_2)|^2\rho_1|\log\rho_1|^{1+\eta}d\rho_1\big)
r_1^{-1}|\log r_1|^{k-2-\eta}dr_1\big\}r_2^{-1}|\log r_2|^{l-2}dr_2\\
&=&{\frac{4}{\eta}}\int_0^A\big\{\lim_{\epsilon\to
0}\big(-\big(\int_0^{r_1}|f_{1,0}^1(\rho_1, r_2)|^2\rho_1|
\log\rho_1|^{1+\eta}d\rho_1\big)\big(\int_{r_1}^A
|\log\rho_1|^{k-2-\eta}{\rho_1}^{-1}d\rho_1\big)\big)|^A_{\epsilon}\\
&&+\int_0^A|f_{1,0}^1(r_1, r_2)|^2r_1|\log
r_1|^{1+\eta}\big(\int_{r_1}^A|\log\rho_1|^{k-2-\eta}
{\rho_1}^{-1}d\rho_1\big)dr_1\big\}r_2^{-1}|\log r_2|^{l-2}dr_2\\
&=&{\frac{4}{\eta(k-1-\eta)}}\int_0^A\big\{\int_0^A
|f_{1,0}^1(r_1, r_2)|^2r_1|\log r_1|^{k}dr_1
\big\}r_2^{-1}|\log r_2|^{l-2}dr_2\\
&=&{\frac{4}{\eta(k-1-\eta)}} \iint_{[0,A]^2}|f_{1,0}^1(r_1,
r_2)|^2|\log r_1|^{k}|\log r_2|^{l-2}r_1r_2^{-1}dr_1dr_2;
\end{eqnarray*}
similarly, using the other path, one can also get
\begin{eqnarray*}
&&\iint_{[0,A]^2} |u_{0,0}(r_1,r_2)|^2r_1^{-1}r_2^{-1}
|\log r_1|^{k-2}|\log r_2|^{l-2}dr_1dr_2\\
&\le&{\frac{4}{\eta(k-1-\eta)}} \iint_{[0,A]^2}|f_{0,1}^2(r_1,
r_2)|^2|\log r_1|^{k-2}|\log r_2|^{l}r_1^{-1}r_2dr_1dr_2.
\end{eqnarray*}
For the other two subcases ($m<0$; $n=0$, $l>1$ and $m=0$, $k>1$;
$n<0$), we can get the same estimates as above.

For the case of the oriented curves from $(A, A)$ to $(r_1, r_2)$,
we have a similar estimate. In this case, the curve integral
above, by two different paths, can be rewritten as
\begin{eqnarray*}
u_{m,n}(r_1, r_2)&=&
- 2r_2^n\int_{r_2}^A\rho_2^{-n}f_{m,n+1}^2(r_1, \rho_2){\text d}\rho_2 \\
({\text{resp}}. &=&
-2r_1^m\int_{r_1}^A\rho_1^{-m}f^1_{m+1,n}(\rho_1, r_2){\text
d}\rho_1).
\end{eqnarray*}
Supposing $m>0$ and $n>0$, one has the following estimate
\begin{eqnarray*}
&&\iint_{[0,A]^2} |u_{m,n}(r_1,r_2)|^2r_1^{-1}r_2^{-1}
|\log r_1|^{k-2}|\log r_2|^{l-2}dr_1dr_2\\
&\le&
4\iint_{[0,A]^2}r_2^{2n}\big(\int_{r_2}^A\rho_2^{-n}f_{m,n+1}^2(r_1,
\rho_2){\text d}\rho_2
\big)^2r_1^{-1}r_2^{-1}|\log r_1|^{k-2}|\log r_2|^{l-2}dr_1dr_2\\
\end{eqnarray*}
\begin{eqnarray*}
&=&4\int_0^A\big\{\int_0^Ar_2^{2n}\big(\int_{r_2}^A\rho_2^{-n}f^2_{m,n+1}(r_1,
\rho_2)d\rho_2\big)^2
r_2^{-1}|\log r_2|^{l-2}dr_2\big\}r_1^{-1}|\log r_1|^{k-2}dr_1\\
&=&4\int_0^A\big\{\int_0^Ar_2^{2n}\big(\int_{r_2}^A\rho_2^{-2n+1}f^2_{m,n+1}(r_1,
\rho_2)d\rho_2\big)
\big(\int_{r_2}^A{\frac{1}{\rho_2}}d\rho_2\big) r_2^{-1}|\log
r_2|^{l-2}dr_2 \big\}   \cdot   \\
&&r_1^{-1}|\log r_1|^{k-2}dr_1\\
&\sim&4\int_0^A\big\{\int_0^A r_2^{2n-1}|\log r_2|^{l-1}
\big(\int_{r_2}^A\rho_2^{-2n+1}f^2_{m,n+1}(r_1, \rho_2)d\rho_2\big)dr_2
\big\}r_1^{-1}|\log r_1|^{k-2}dr_1\\
&=&4\int_0^A\big\{\lim_{\epsilon\to 0}\big(\int_0^{r_2}\rho_2^{2n-1}
|\log\rho_2|^{l-1}d\rho_2
\cdot\int_{r_2}^A\rho_2^{-2n+1}f^2_{m,n+1}(r_1, \rho_2)d\rho_2\big)|_{\epsilon}^A \\
&&+\int_0^A\big(\int_0^{r_2}\rho_2^{2n-1}|\log\rho_2|^{l-1}d\rho_2\big)
r_2^{-2n+1}f^2_{m,n+1}(r_1, r_2)dr_2
\big\}r_1^{-1}|\log r_1|^{k-2}dr_1\\
&=&4\int_0^A\big\{\int_0^A\big(\int_0^{r_2}\rho_2^{2n-1}
|\log\rho_2|^{l-1}d\rho_2\big)r_2^{-2n+1}f^2_{m,n+1}(r_1, r_2)dr_2
\big\}r_1^{-1}|\log r_1|^{k-2}dr_1\\
&\le&4\int_0^A\big\{\int_0^Ar_2^{2n}|\log r_2|^{l-1}\cdot
r_2^{-2n+1}|f^2_{m,n+1}(r_1, r_2)|^2dr_2
\big\}r_1^{-1}|\log r_1|^{k-2}dr_1\\
&&({\text{$\rho^{n}|\log\rho|^{l}$ is an increasing function as $n>0$.}})\\
&\le&4\int_0^A\big\{\int_0^A|f^2_{m,n+1}(r_1, r_2)|^2r_2|\log r_2|^{l}dr_2
\big\}r_1^{-1}|\log r_1|^{k-2}dr_1\\
&=&4\iint_{[0,A]^2}|f^2_{m,n+1}(r_1, r_2)|^2 r_1^{-1}r_2|\log
r_1|^{k-2}|\log r_2|^{l}dr_1dr_2;
\end{eqnarray*}
similarly, using the other path, one can also show
\begin{eqnarray*}
&&\iint_{[0,A]^2} |u_{m,n}(r_1,r_2)|^2r_1^{-1}r_2^{-1}
|\log r_1|^{k-2}|\log r_2|^{l-2}dr_1dr_2\\
&=&4\iint_{[0,A]^2}|f^1_{m+1,n}(r_1, r_2)|^2 r_1r_2^{-1}|\log
r_1|^{k}|\log r_2|^{l-2}dr_1dr_2.
\end{eqnarray*}
Supposing $m=0$, $k<1$ and $n=0$, $l<1$, one then has
\begin{eqnarray*}
&&\iint_{[0,A]^2}|u_{0,0}(r_1,r_2)|^2r_1^{-1}r_2^{-1}|\log r_1|^{k-2}|\log r_2|^{l-2}dr_1dr_2\\
&\le&4\iint_{[0,A]^2}\big(\int_{r_2}^Af^2_{0,1}(r_1, \rho_2){\text
d}\rho_2
\big)^2r_1^{-1}r_2^{-1}|\log r_1|^{k-2}|\log r_2|^{l-2}dr_1dr_2\\
&=&4\int_0^A\big\{\int_0^A\big(\int_{r_2}^Af^2_{0,1}(r_1, \rho_2){\text d}\rho_2\big)^2r_2^{-1}|\log r_2|^{l-2}dr_2
\big\}r_1^{-1}|\log r_1|^{k-2}dr_1\\
&\le&4\int_0^A\big\{\int_0^A\big(\int_{r_2}^A
|f^2_{0,1}(r_1, \rho_2)|^2\rho_2|\log\rho_2|^{1-\eta}d\rho_2\big)\cdot\\
\end{eqnarray*}
\begin{eqnarray*}
&&\big(\int_{r_2}^A{\rho_2}^{-1}|\log\rho_2|^{-1+\eta}d\rho\big)r_2^{-1}|\log
r_2|^{l-2}dr_2
\big\}r_1^{-1}|\log r_1|^{k-2}dr_1\\
&&({\text{$\eta$ is some positive constant $<1$.}})\\
&\sim&{\frac{4}{\eta}}\int_0^A
\big\{\int_0^A\big(\int_{r_2}^A|f^2_{0,1}(r_1, \rho_2)|^2\rho_2
|\log\rho_2|^{1-\eta}d\rho_2\big)\cdot
r_2^{-1}|\log r_2|^{l-2+\eta}dr_2\big\}r_1^{-1}|\log r_1|^{k-2}dr_1\\
&=&{\frac{4}{\eta}}\int_0^A\big\{\big(\lim_{\eta\to 0}
\big(\int_{r_2}^A|f^2_{0,1}(r_1, \rho_2)|^2\rho_2|\log\rho_2|^{1-\eta}d\rho_2\big)
\big(\int_0^{r_2}\rho_2^{-1}|\log\rho_2|^{l-2+\eta}d\rho_2\big)\big)|_\eta^A\\
&&+\int_0^A|f^2_{0,1}(r_1, r_2)|^2r_2|\log r_2|^{1-\eta}
\big(\int_0^{r_2}\rho_2^{-1}|\log\rho_2|^{l-2+\eta}d\rho_2\big)dr_2
\big\}r_1^{-1}|\log r_1|^{k-2}dr_1\\
&=&{\frac{4}{\eta(1-l-\eta)}}\int_0^A\big\{
\int_0^A|f^2_{0,1}(r_1, r_2)|^2r_2|\log r_2|^{1-\eta}|\log r_2|^{l-1+\eta}dr_2
\big\}r_1^{-1}|\log r_1|^{k-2}dr_1\\
&=&{\frac{4}{\eta(1-l-\eta)}}\int_0^A\big\{
\int_0^A|f^2_{0,1}(r_1, r_2)|^2r_2|\log r_2|^ldr_2
\big\}r_1^{-1}|\log r_1|^{k-2}dr_1\\
&=&{\frac{4}{\eta(1-l-\eta)}}\iint_{[0,A]^2} |f^2_{0,1}(r_1,
r_2)|^2r_1^{-1}r_2|\log r_1|^{k-2}|\log r_2|^ldr_1dr_2;
\end{eqnarray*}
similarly, using the other path, one also has
\begin{eqnarray*}
&&\iint_{[0,A]^2}|u_{0,0}(r_1,r_2)|^2r_1^{-1}r_2^{-1}
|\log r_1|^{k-2}|\log r_2|^{l-2}dr_1dr_2\\
&=&{\frac{4}{\eta(1-l-\eta)}}\iint_{[0,A]^2} |f^1_{1,0}(r_1,
r_2)|^2r_1r_2^{-1}|\log r_1|^{k}|\log r_2|^{l-2}dr_1dr_2.
\end{eqnarray*}
One can deal similarly with the two other subcases of this case
($m>0$; $n=0$, $l<1$ and $m=0$, $k<1$; $n>0$).

For the two other cases for the curve integral, one can obtain the
same estimates as above, we omit these.

\noindent {\bf Case 2:} Writing $\phi$ as $f d{\overline
t_1}\wedge d{\overline t}_2$, the condition of $\phi$ can then be
restated as
\[
{\iint_{(\triangle^*)^2}}|f|^2 \Vert\sigma\Vert^2 {dt_1\wedge
d{\overline t}_1}\wedge{dt_2\wedge d{\overline t}_2} < \infty.
\]
We now need to prove that there exists a form $\psi =
u^1d{\overline{t}}_1 + u^2d{\overline{t}}_2$ of type $(0, 1)$
satisfying
\[
{\frac{{\partial}u^2}{\partial{\overline{t}_1}}} -
{\frac{{\partial}u^1}{\partial{\overline{t}_2}}} = f
\]
and
\begin{eqnarray*}
&&{\iint_{(\triangle^*)^2}}
(|u^1|^2|t_2|^{-2}(-{\text{log}}|t_2|)^{-2}) \Vert\sigma\Vert^2
{dt_1\wedge d{\overline t}_1}\wedge{dt_2\wedge d{\overline t}_2} +\\
&&{\iint_{(\triangle^*)^2}}
(|u^2|^2|t_1|^{-2}(-{\text{log}}|t_1|)^{-2}) \Vert\sigma\Vert^2
{dt_1\wedge d{\overline t}_1}\wedge{dt_2\wedge d{\overline t}_2} <
\infty.
\end{eqnarray*}
As in the case $1$, we still use the method of Fourier series.
Write $f$, $u^1$, and $u^2$ as $(r_1, r_2)$-dependent Fourier
series as follows
\begin{eqnarray*}
&&f=\sum f_{m,n}(r_1, r_2)\exp(\sqrt{-1}m\theta_1+\sqrt{-1}n\theta_2), \\
&&u^1=\sum u^1_{m,n}(r_1, r_2)\exp(\sqrt{-1}m\theta_1+\sqrt{-1}n\theta_2), \\
&&u^2=\sum u^2_{m,n}(r_1, r_2)\exp(\sqrt{-1}m\theta_1+\sqrt{-1}n\theta_2).
\end{eqnarray*}
The above equation of $u^1$ and $u^2$ can then be rewritten as
\[
{\frac{\partial}{\partial r_1}}(r_1^{-m}r_2^{-n}u^2_{m,n+1})
-{\frac{\partial}{\partial r_2}}(r_1^{-m}r_2^{-n}u^1_{m+1,n})
=2r_1^{-m}r_2^{-n}f_{m+1,n+1}.
\]
Put
\begin{eqnarray*}
u^1_{m,n}(r_1, r_2) =
\begin{cases}
-r_2^n\int_0^{r_2}\rho_2^{-n}f_{m,n+1}(r_1, \rho_2)
d\rho_2,~n<0~{\text{or}}~n=0, l>1;\\
 r_2^n\int_{r_2}^A\rho_2^{-n}f_{m,n+1}(r_1, \rho_2)
 d\rho_2,~n>0~{\text{or}}~n=0, l<1
\end{cases}
\end{eqnarray*}
and
\begin{eqnarray*}
u^2_{m,n}(r_1, r_2) =
\begin{cases}
 r_1^m\int_0^{r_1}\rho_1^{-m}f_{m+1,n}(\rho_1, r_2)
 d\rho_1,~m<0~{\text{or}}~m=0, k>1;\\
-r_1^m\int_{r_1}^A\rho_1^{-m}f_{m+1,n}(\rho_1,
r_2)d\rho_1,~m>0~{\text{or}}~m=0, k<1.
\end{cases}
\end{eqnarray*}
We now want to show that
\begin{eqnarray*}
&&\sum_{m,n}\iint_{[0,A]^2}|u^1_{m,n}(r_1, r_2)|^2
|\log r_1|^k|\log r_2|^{l-2}r_1r_2^{-1}dr_1dr_2\\
&&+\sum_{m,n}\iint_{[0,A]^2}|u^2_{m,n}(r_1, r_2)|^2
|\log r_1|^{k-2}|\log r_2|^{l}r_1^{-1}r_2dr_1dr_2\\
&\le&C\sum_{m,n}\iint_{[0,A]^2}|f_{m,n}(r_1, r_2)|^2 |\log
r_1|^k|\log r_2|^{l}r_1r_2dr_1dr_2
\end{eqnarray*}
for some positive constant $C$. We perform the estimate only for
$u^1_{m,n}$ with $n<0$, the other cases are similar,
\begin{eqnarray*}
&&\iint_{[0,A]^2}|u^1_{m,n}(r_1, r_2)|^2
|\log r_1|^k|\log r_2|^{l-2}r_1r_2^{-1}dr_1dr_2\\
&=&\iint_{[0,A]^2}\big(-r_2^n\int_0^{r_2}
\rho_2^{-n}f_{m,n+1}(r_1, \rho_2)d\rho_2
\big)^2|\log r_1|^k|\log r_2|^{l-2}r_1r_2^{-1}dr_1dr_2\\
&=&\int_0^A\big\{\int_0^Ar_2^{2n}
\big(\int_0^{r_2}\rho_2^{-2n}|f_{m,n+1}(r_1, \rho_2)|^2d\rho_2\big)
\big(\int_0^{r_2}d\rho_2\big)
|\log r_2|^{l-2}r_2^{-1}dr_2\big\}|\log r_1|^kr_1dr_1\\
\end{eqnarray*}
\begin{eqnarray*}
&=&\int_0^A\big\{\int_0^Ar_2^{2n}
\big(\int_0^{r_2}\rho_2^{-2n}|f_{m,n+1}(r_1,
\rho_2)|^2d\rho_2\big)
|\log r_2|^{l-2}dr_2\big\}|\log r_1|^kr_1dr_1\\
&\sim&\int_0^A\big\{\int_0^A|f_{m,n+1}(r_1, r_2)|^2
|\log r_2|^{l-1}r_2dr_2\big\}|\log r_1|^kr_1dr_1\\
&\le&\iint_{[0,A]^2}|f_{m,n+1}(r_1, r_2)|^2 |\log r_1|^k|\log
r_2|^lr_1r_2dr_1dr_2.
\end{eqnarray*}
Thus, we complete the proof of the assertion, and hence the
proposition. ~~~$\Box$

\vskip 4mm \noindent {\bf Remarks.} 1)~~  Near the divisor, denote
by $L_1$ (resp. $L_2$) the line bundle generated by
${\frac{dt_1}{t_1}}$ (resp. ${\frac{dt_2}{t_2}}$); under the
Poincar\'e-like metric, $\Vert{\frac{dt_1}{t_1}}\Vert^2\sim
|\log|t_1||^2$ and $\Vert{\frac{dt_2}{t_2}}\Vert^2\sim
|\log|t_2||^2$. So, by the previous proposition, if $k\neq -1$ and
$l\neq 1$ (resp. $k\neq 1$ and $l\neq -1$; $k\neq -1$ and $l\neq
-1$), the $L^2$ $\overline\partial$-Poincar\'e lemma is true for
$L^2$ $(0, s)$-forms ($s\ge 1$) with values in $L_1\otimes L$
(resp. $L_2\otimes L$; $L_1\otimes L_2\otimes L$). On the other
hand, by the proposition in \S4, any $L^2$ $(1, s)$-form valued in
${\bf E}^p$ can be written as the sum of some forms valued in
$L_1\otimes {\bf E}^p$ or $L_2\otimes {\bf E}^p$ and any $L^2$
$(2, s)$-forms valued in ${\bf E}^p$ can be written as a form
valued in $L_1\wedge L_2\otimes {\bf E}^p$, that is why we proved
the $L^2$ $\overline\partial$-Poincar\'e lemma only for $r=0$. All
these will be used in the following proof of quasi-isomorphism.
2)~~ The classical $L^2$ existence theorem for the
$\overline{\partial}$-problem \cite{av, ho} (refer to \cite{dem})
works only for partial cases. This can be seen from the following
discussion. Let $L$ be a holomorphic line bundle on
$(\triangle^*)^2$ with a generating section $\sigma$ and a
Hermitian metric satisfying, near the divisor,
\[
{\Vert\sigma\Vert}^2 =
(-{\text{log}}|t_1|)^{k}(-{\text{log}}|t_2|)^{l}~~ {\text{for some
real numbers $k$ and $l$}}.
\]
Take the Poincare-like metric on $(\triangle^*)^2$ as before
\[
\omega = {\frac{\sqrt{-1}}{2}}\big({\frac{dt_1\wedge d{\overline t}_1}
{|t_1|^2(-{\text{log}}|t_1|)^2}} +
{\frac{dt_2\wedge d{\overline t}_2}{|t_2|^2(-{\text{log}}|t_2|)^2}}\big).
\]
A standard computation tells us that the curvature of the metric
on the line bundle $L$ is
\[
{\sqrt{-1}}\Theta(L) =
{\frac{\sqrt{-1}}{2}}\big({\frac{kdt_1\wedge d{\overline t}_1}
{|t_1|^2(-{\text{log}}|t_1|)^2}} +
{\frac{ldt_2\wedge d{\overline t}_2}{|t_2|^2(-{\text{log}}|t_2|)^2}}\big).
\]
Using the notations of page 30 of \cite{dem}, the curvature
eigenvalues, relative to $\omega$, are $\gamma_1=k$ and
$\gamma_2=l$.(Here, without loss of generality, we have assumed
that $k\le l$.) Then for any $L$-valued $(p, q)$-form $u=\sum
u_{JK}dz_{J}\wedge d{\overline z}_{K}\otimes\sigma$,
\begin{eqnarray*}
\langle[{\sqrt{-1}}\Theta(L), \Lambda ]u, u\rangle &=&
\sum_{|J|=p, |K|=q} \big(\sum_{j\in J}\gamma_j + \sum_{j\in
K}\gamma_j -
\sum_{1\le j\le n}\gamma_j\big)|u_{JK}|^2          \\
&\ge& (\gamma_1 +\cdots + \gamma_q - \gamma_{p+1} - \cdots - \gamma_2)|u|^2.
\end{eqnarray*}
So, by the $L^2$ existence theorem for the
$\overline{\partial}$-problem (refer to Theorem 5.1 of
\cite{dem}), one has that if $\gamma_1 +\cdots + \gamma_q -
\gamma_{p+1} - \cdots - \gamma_2 > 0$, $q\ge 1$, then for each
$L$-valued $L^2$ $(p,q)$ form  $\psi$ on $(\triangle^*)^2$ with
$\overline{\partial}\psi = 0$, there exists an $L$-valued $L^2$
$(p,q-1)$ form $\phi$ on $(\triangle^*)^2$ such that
$\overline{\partial}\phi = \psi$ and
\[
\Vert\phi\Vert_{(2)}\le C\Vert\psi\Vert_{(2)}.
\]
Here, $C=C(k, l, p, q)$ is a constant depending only on $k, l, p,
q$. In particular, if $u$ is a $(0, 1)$-form, the condition above
is $k-k-l=-l > 0$; if $u$ is a $(0, 2)$-form, the condition above
is $k+l-k-l=0$, which is empty; if $u$ is a $(1, 1)$-form, the
condition above is $k-l \le 0$, which is also empty; if $u$ is a
$(1, 2)$-form, the condition above is $k+l-l=k > 0$; if $u$ is a
$(2, 1)$-form, the condition above is $k > 0$; if $u$ is a $(2,
2)$-form, the condition above is $k+l > 0$. For the higher
dimensional case, we have a similar argument.

\vskip 4mm \noindent {\bf The proof of quasi-isomorphism:} Fix
$p\ge 0$, we now prove that the complex
$\{{\text{Gr}}^p_F\Omega^{\cdot}({\bf H}_{\mathbb{C}})_{(2)}
=(\Omega^{\cdot}\otimes{\bf E}^{p-\cdot})_{(2)}, \theta\}$ is
quasi-isomorphic to the complex $\{[{\text{Gr}}^p_FA^{\cdot}({\bf
H}_{\mathbb{C}})]_{(2)} =\oplus_{r+s=\cdot}(A^{r, s}\otimes{\bf
E}^{p-r})_{(2)}, D''={\overline\partial}+\theta\}$. It is
sufficient to show that the corresponding cohomological sheaves of
the two complexes are the same near the singularity. So, in the
following, we always restrict our discussion to some small
neighborhood near the singularity. We first review some facts
which are established (though which may not be stated explicitly,
but are clear from the context) in \S3 . We know that, if
neglecting the higher order terms, $\theta\sim
N_1\otimes{\frac{dt_1}{t_1}} + N_2\otimes{\frac{dt_1}{t_1}}$; on
the other hand, the higher order terms essentially play no role,
so, we assume that $\theta = N_1\otimes{\frac{dt_1}{t_1}} +
N_2\otimes{\frac{dt_1}{t_1}}$ in the following discussion. From
the construction of $\sigma_{k, l}$, as an element of some
$L^2$-adapted basis of $({\bf E}^{k+l})_{(2)}$, one has that, if
both ${\frac{dt_1}{t_1}}\otimes N_1\sigma_{k, l}$ and
${\frac{dt_2}{t_2}}\otimes N_2\sigma_{k, l}$ are nonzero,
\[
\Vert{\frac{dt_1}{t_1}}\otimes N_1\sigma_{k, l}\Vert^2\sim
\Vert\sigma_{k, l}\Vert^2~~{\text{and}}~~
\Vert{\frac{dt_2}{t_2}}\otimes N_2\sigma_{k, l}\Vert^2\sim
\Vert\sigma_{k, l}\Vert^2,
\]
and they can be considered as two elements in an $L^2$ adapted
basis $(\Omega^1_M\otimes({\bf E}^{k+l-1}))_{(2)}$; if
$\Vert\sigma_{k, l}\Vert^2\sim|\log|t_1|||\log|t_2||^{2l-n}$
(resp. $|\log|t_1||^{2k-m}|\log|t_2||$, $|\log|t_1|||\log|t_2||$),
then $N_1\sigma_{k, l}$ (resp. $N_2\sigma_{k, l}$,
$N_1N_2\sigma_{k, l}$) is not equal zero; if
$\Vert{\frac{dt_1}{t_1}}\otimes \sigma_{k, l}\Vert^2\sim
\log|t_1|||\log|t_2||^{2l-n}$ (resp.
$\Vert{\frac{dt_2}{t_2}}\otimes \sigma_{k, l}\Vert^2\sim
\log|t_1||^{2k-m}|\log|t_2||$), then there exists an element
$\sigma_{k+1, l}$ (resp. $\sigma_{k, l+1}$) in some $L^2$-adapted
basis $\{\sigma_{s,t}\}$ {\footnote {which, very probably, is
slightly different from, but is equal to, up to a higher order
term and hence in the following is viewed as, the one constructed
in the previous section.}} of $({\bf E}^{k+l+1})_{(2)}$ with
$N_1\sigma_{k+1, l}=\sigma_{k, l}$ (resp. $N_2\sigma_{k,
l+1}=\sigma_{k, l}$) and
\[
\Vert\sigma_{k+1, l}\Vert^2\sim\Vert{\frac{dt_1}{t_1}}\otimes\sigma_{k, l}
\Vert^2~~({\text{resp.}}~~
\Vert\sigma_{k, l+1}\Vert^2\sim\Vert{\frac{dt_2}{t_2}}
\otimes\sigma_{k, l}\Vert^2).
\]

Let $U$ be an open subset containing the origin $(0, 0)\in
(\triangle)^2$. In the following argument, we sometimes need to
shrink $U$ to a smaller open subset containing the origin, this
will not be pointed out explicitly since it will be obvious from
the context. In order to prove that the two cohomological sheaves
are the same near the singularity, it is sufficient to show that
{\it for any $D''$-closed form $\phi\in[{\text{Gr}}_F^pA^k({\bf
H}_{\mathbb{C}})]_{(2)}$ on $U$, there are a $\theta$-closed form
$\eta\in {\text{Gr}}_F^p{\Omega}^k({\bf H}_{\mathbb{C}})_{(2)}$
and a form $\psi\in[{\text{Gr}}_F^pA^{k-1}({\bf
H}_{\mathbb{C}})]_{(2)}$ on $U$ satisfying $\phi = \eta +
D''\psi$, $k=0, 1, 2, 3, 4$}. In the following, we only consider
the case of $k=2$, the other cases are similar, even much easier.
It is easy to see that the complex
$\{[{\text{Gr}}_F^pA^{\cdot}({\bf H}_{\mathbb{C}})]_{(2)}, D''\}$
is actually a double complex with the differentials
$\{\overline\partial, \theta\}$. Putting the two complexes
together, we have the following diagram
\begin{equation*}
\begin{CD}
{{\bf E}^p_{(2)}}@>>>{(A^{0,0}\otimes{\bf E}^p)_{(2)}}@>{\overline\partial}>>
{(A^{0,1}\otimes{\bf E}^p)_{(2)}}@>{\overline\partial}>>
{(A^{0,2}\otimes{\bf E}^p)_{(2)}}@>>>{0}\\
@VV{\theta}V @VV{\theta}V @VV{\theta}V @VV{\theta}V \\
{(\Omega^1\otimes{\bf E}^{p-1})_{(2)}}@>>>{(A^{1,0}
\otimes{\bf E}^{p-1})_{(2)}}@>{\overline\partial}>>
{(A^{1,1}\otimes{\bf E}^{p-1})_{(2)}}@>{\overline\partial}>>
{(A^{1,2}\otimes{\bf E}^{p-1})_{(2)}}@>>>{0}\\
@VV{\theta}V @VV{\theta}V @VV{\theta}V @VV{\theta}V \\
{(\Omega^2\otimes{\bf E}^{p-2})_{(2)}}@>>>{(A^{2,0}
\otimes{\bf E}^{p-2})_{(2)}}@>{\overline\partial}>>
{(A^{2,1}\otimes{\bf E}^{p-2})_{(2)}}@>{\overline\partial}>>
{(A^{2,2}\otimes{\bf E}^{p-2})_{(2)}}@>>>{0}.
\end{CD}
\end{equation*}
Let $\phi\in [{\text{Gr}}_F^pA^{2}({\bf H}_{\mathbb{C}})]_{(2)}$
be a $D''$-closed form on $U$. By the consideration of type, it
can be decomposed into $\phi=\phi^{2,0}+\phi^{1,1}+\phi^{0,2}$
with $\phi^{r,2-r}\in (A^{r,2-r}\otimes{\bf E}^{p-r})_{(2)}$
satisfying $\overline\partial\phi^{r,2-r}+\theta\phi^{r-1,3-r}=0,
r=0, 1, 2, 3$ ($\phi^{3,-1}, \phi^{-1,3}$=0 automatically). We
first consider the parts of types $(0, 2)$ and $(1, 1)$:
$\phi^{0,2}$ and $\phi^{1, 1}$. {\it We want to eliminate the part
$\phi^{0,2}$ from $\phi$ by the solutions for the
$\overline\partial$-problem and the properties of $\theta$}, more
precisely, we want to find some $\psi\in [{\text{Gr}}_F^pA^1({\bf
H}_{\mathbb{C}})]_{(2)}$ satisfying $\phi-D''\psi$ not containing
the part of type $(0, 2)$. By the construction in the previous
section, near the singularity, ${\bf E}^p$ can be decomposed into
the holomorphic direct sum of some trivial line bundles according
to an $L^2$-adapted basis of ${\bf E}^p_{(2)}$
$\{\sigma_{k,l}^p\}$, satisfying
\[
\Vert\sigma_{k,l}^p\Vert^2\sim|\log|t_1||^{k'}|\log|t_2||^{l'}
\]
for some integers $k', l'$ depending on $k, l$ respectively. So,
$\phi^{0,2}=\sum\phi^{0,2}_{p,k,l}\otimes\sigma^p_{k,l}$, here
$\{\phi^{0,2}_{p,k,l}\}$ are some $(0, 2)$-forms and each summand
is $L^2$. If for a $\sigma_{k,l}$, $k', l'\neq 1$, then by the
previous proposition, there exists a local section
$\psi^{0,1}_{p,k,l}\otimes\sigma^p_{k,l}\in(A^{0,1}\otimes{\bf
E}^p)_{(2)}(U)$ satisfying
$\overline\partial(\psi^{0,1}_{p,k,l}\otimes\sigma^p_{k,l})
=\phi^{0,2}_{p,k,l}\otimes\sigma^p_{k,l}$.
Thus the part of type $(0,2)$ of
$\phi-D''(\psi^{0,1}_{p,k,l}\otimes\sigma^p_{k,l})$ does not
contain a summand with value in the line bundle generated by
$\sigma^p_{k,l}$. Therefore we can assume that $\phi^{0,2}$ has
only summands of form $\phi^{0,2}_{p,k,l}\otimes\sigma^p_{k,l}$
satisfying
\[
\Vert\sigma^p_{k,l}\Vert^2\sim|\log|t_1|||\log|t_2||^{l'}~~{\text{or}}~~
|\log|t_1||^{k'}|\log|t_2||.
\]
Let $\phi^{0,2}_{p,k,l}\otimes\sigma^p_{k,l}$ be such a summand
with $\Vert\sigma^p_{k,l}\Vert^2\sim|\log|t_1|||\log|t_2||^{l'}$
(the other case is similar).  From the above review, we know that
${\frac{dt_1}{t_1}}\otimes N_1\sigma^p_{k,l}$ and
${\frac{dt_2}{t_2}}\otimes N_2\sigma^p_{k,l}$ (if $\neq 0$) are
two elements of some $L^2$-adapted basis of $(\Omega^1\otimes{\bf
E}^{p-1})_{(2)}$. So, by
$\theta\phi^{0,2}+\overline\partial\phi^{1,1}=0$, we have that
there exist two summands
$\phi^{0,1}_1\wedge{\frac{dt_1}{t_1}}\otimes N_1\sigma^p_{k,l}$
and $\phi^{0,1}_2\wedge{\frac{dt_2}{t_2}}\otimes
N_2\sigma^p_{k,l}$ of $\phi^{1,1}$ (when decomposed according to
some $L^2$-adapted basis containing ${\frac{dt_1}{t_1}}\otimes
N_1\sigma^p_{k,l}$ and ${\frac{dt_2}{t_2}}\otimes
N_2\sigma^p_{k,l}$ of $(\Omega^1\otimes{\bf E}^{p-1})_{(2)}$) with
$\overline\partial\phi^{0,1}_1=\overline\partial\phi^{0,1}_2=\phi^{0,2}_{p,k,l}$.
Set $\psi=\phi^{0,1}_1\otimes\sigma^p_{k,l}$. Since
$\Vert\sigma^p_{k,l}\Vert\sim\Vert{\frac{dt_1}{t_1}}\otimes
N_1\sigma^p_{k,l}\Vert$ and
$\overline\partial\psi=\phi^{0,2}_{p,k,l}\otimes\sigma^p_{k,l}$,
so $\psi\in(A^{0,1}\otimes{\bf E}^p)_{(2)}$. Thus again the part
of type $(0, 2)$ of $\phi-D''\psi$ does not contain
$\phi^{0,2}_{p,k,l}\otimes\sigma^p_{k,l}$ as a summand when
decomposed according to some $L^2$-adapted basis of ${\bf
E}^p_{(2)}$ containing $\sigma^p_{k,l}$. Similarly, we can
eliminate other such terms of $\phi$.

Summing the above argument up, we now assume that $\phi$ does not
contain the part of type $(0, 2)$. Namely,
$\phi=\phi^{2,0}+\phi^{1,1}$. We now want to eliminate the term
$\phi^{1,1}$ in $\phi$ by using the same method as above. But as
seen in the following, the situation is slightly different. Choose
an $L^2$-adapted basis $\{\sigma_{k,l}^{p-1}\}$ of ${\bf
E}^{p-1}_{(2)}$ with
$\Vert\sigma^{p-1}_{k,l}\Vert^2\sim|\log|t_1||^{k'}|\log|t_2||^{l'}$,
where $k', l'$ depends on $k, l$ respectively. It is easy to see
that $\{{\frac{dt_1}{t_1}}\otimes\sigma^{p-1}_{k,l},
{\frac{dt_2}{t_2}}\otimes\sigma^{p-1}_{k,l}\}$ can be considered
as an $L^2$-adapted basis of $(\Omega^1\otimes{\bf
E}^{p-1})_{(2)}$. According to the basis, $\phi$ can then be
decomposed as follows
\[
\phi^{1,1}=\sum\phi^{0,1}_{p-1,k,l,1}\wedge{\frac{dt_1}{t_1}}
\otimes\sigma^{p-1}_{k,l}
+\sum\phi^{0,1}_{p-1,k,l,2}\wedge{\frac{dt_2}{t_2}}
\otimes\sigma^{p-1}_{k,l},
\]
for some $(0, 1)$-forms $\phi^{0,1}_{p-1,k,l,1}$ and
$\phi^{0,1}_{p-1,k,l,2}$. Since $\phi^{0,2}=0$, each term above is
zero under $\overline\partial$. So, if
$\Vert{\frac{dt_1}{t_1}}\otimes\sigma^{p-1}_{k,l}\Vert^2
\sim|\log|t_1||^{k_1}|\log|t_2||^{l_1}$
with $k_1, l_1\neq 1$, then the $\overline\partial$-problem tells
us that there exists a local section
$\psi^{0,0}_{p-1,k,l,1}{\frac{dt_1}{t_1}}\otimes\sigma^{p-1}_{k,l}\in
(A^{1,0}\otimes{\bf E}^{p-1})_{(2)}(U)$ satisfying
$\overline\partial(\psi^{0,0}_{p-1,k,l,1}{\frac{dt_1}{t_1}}
\otimes\sigma^{p-1}_{k,l})
=\phi^{0,1}_{p-1,k,l,1}{\frac{dt_1}{t_1}}\otimes\sigma^{p-1}_{k,l}$;
thus the part of type $(1, 1)$ of
$\phi-D''(\psi^{0,0}_{p-1,k,l,1}{\frac{dt_1}{t_1}}\otimes\sigma^{p-1}_{k,l})$
does not contain a $(0, 1)$-form with value in the line bundle
generated by ${\frac{dt_1}{t_1}}\otimes\sigma^{p-1}_{k,l}$. So, we
can assume that the above decomposition of $\phi^{1,1}$ contains
only the terms satisfying
\[
\Vert{\frac{dt_1}{t_1}}\otimes\sigma^{p-1}_{k,l}\Vert^2
\sim|\log|t_1|||\log|t_2||^{l_1}
~~{\text{or}}~~|\log|t_1||^{k_1}|\log|t_2||
\]
and
\[
\Vert{\frac{dt_2}{t_2}}\otimes\sigma^{p-1}_{k,l}\Vert^2
\sim|\log|t_1|||\log|t_2||^{l_1}
~~{\text{or}}~~|\log|t_1||^{k_1}|\log|t_2||.
\]
We can do a further reduction as follows: Let the decomposition of
$\phi^{1,1}$ have a term satisfying
$\Vert{\frac{dt_1}{t_1}}\otimes\sigma^{p-1}_{k,l}\Vert^2
\sim|\log|t_1|||\log|t_2||^{l_1}$. By the review before, there
exists an element $\sigma^{p}_{k+1,l}$, lying in some
$L^2$-adapted basis of ${\bf E}^p_{(2)}$, satisfying
$N_1\sigma^{p}_{k+1,l}=\sigma^{p-1}_{k,l}$ and
$\Vert\sigma^{p}_{k+1,l}\Vert^2\sim|\log|t_1|||\log|t_2||^{l_1}$.
Set $\psi=\phi^{0,1}_{p-1,k,l,1}\otimes\sigma^{p}_{k+1,l}$, which
is obviously a local section in $(A^{0,1}\otimes{\bf E}^p)_{(2)}$.
Since $\phi^{0,2}=0$, one has
$\overline\partial(\phi^{0,1}_{p-1,k,l,1})=0$. So
$\overline\partial(\psi)=0$. Consider
$\phi-D''\psi=\phi-\theta\psi=\phi-\phi^{0,1}_{p-1,k,l,1}
\wedge{\frac{dt_1}{t_1}}\otimes\sigma^{p-1}_{k,l}
-\phi^{0,1}_{p-1,k,l,1}\wedge{\frac{dt_2}{t_2}}\otimes
N_2\sigma^{p}_{k+1,l}$. We know that , if
${\frac{dt_2}{t_2}}\otimes N_2\sigma^{p}_{k+1,l}\neq 0$,
$\Vert{\frac{dt_2}{t_2}}\otimes
N_2\sigma^{p}_{k+1,l}\Vert^2\sim|\log|t_1|||\log|t_2||^{l_1}$. So,
if doing the same argument for the terms containing
${\frac{dt_2}{t_2}}$, we can always assume that in the
decomposition of $\phi^{1,1}$ there exist only the terms
satisfying
\[
\Vert{\frac{dt_1}{t_1}}\otimes\sigma^{p-1}_{k,l}\Vert^2
\sim|\log|t_1||^{k_1}|\log|t_2||
\]
and
\[
\Vert{\frac{dt_2}{t_2}}\otimes\sigma^{p-1}_{k,l}\Vert^2
\sim|\log|t_1|||\log|t_2||^{l_1}.
\]
In the following, we first discuss the terms containing
${\frac{dt_1}{t_1}}$:
$\phi^{0,1}_{p-1,k,l,1}\wedge{\frac{dt_1}{t_1}}\otimes\sigma^{p-1}_{k,l}$
with $\Vert{\frac{dt_1}{t_1}}\otimes\sigma^{p-1}_{k,l}\Vert^2
\sim|\log|t_1||^{k_1}|\log|t_2||$ and $k_1\neq 1$. The terms
containing ${\frac{dt_2}{t_2}}$ with $l_1\neq 1$ can be discussed
similarly. As pointed out in the review before,
$N_2\sigma^{p-1}_{k,l}\neq 0$. So
${\frac{dt_1}{t_1}}\wedge{\frac{dt_2}{t_2}}\otimes
N_2\sigma^{p-1}_{k,l}$ can be considered as an element of some
$L^2$-adapted basis of $(\Omega^2\otimes{\bf E}^{p-2})_{(2)}$.
Again since $\overline\partial\phi^{2,0}+\theta\phi^{1,1}=0$, so
$\phi^{2,0}$ contains a part of form
$\phi^{0,0}{\frac{dt_1}{t_1}}\wedge{\frac{dt_2}{t_2}}\otimes
N_2\sigma^{p-1}_{k,l} \in(A^{2,0}\otimes{\bf E}^{p-2})_{(2)}$ with
$\overline\partial\phi^{0,0}=\phi^{0,1}_{p-1,k,l,1}$. Set
$\psi=\phi^{0,0}{\frac{dt_1}{t_1}}\otimes\sigma^{p-1}_{k,l}$,
which is clearly a local section in $(A^{1,0}\otimes{\bf
E}^{p-1})_{(2)}$ and satisfies
$\overline\partial\psi=\phi^{0,1}_{p-1,k,l,1}\wedge{\frac{dt_1}{t_1}}
\otimes\sigma^{p-1}_{k,l}$. It is clear that $\phi-D''\psi$ does
not contain the term discussed. So we can assume that $\phi^{1,1}$
contains only two terms {\footnote{Here and in the sequel,
w.l.o.g., we always assume that the $L^2$-adapted basis
$\{\sigma^{p-1}_{k,l}\}$ is constructed using an irreducible
representation of $(\mathfrak{sl}_2({\mathbb{C}}))^2$ by the
method of \S 3.}}
$\phi^{0,1}_{p-1,k,l,1}\wedge{\frac{dt_1}{t_1}}\otimes\sigma^{p-1}_{k,l,1}$
and
$\phi^{0,1}_{p-1,k,l,2}\wedge{\frac{dt_2}{t_2}}\otimes\sigma^{p-1}_{k,l,2}$
satisfying
\[
\Vert{\frac{dt_1}{t_1}}\otimes\sigma^{p-1}_{k,l,1}\Vert^2
\sim|\log|t_1|||\log|t_2||
\]
and
\[
\Vert{\frac{dt_2}{t_2}}\otimes\sigma^{p-1}_{k,l,2}\Vert^2
\sim|\log|t_1|||\log|t_2||,
\]
where $\sigma^{p-1}_{k,l,1}$ and $\sigma^{p-1}_{k,l,2}$ are two
elements of the $L^2$-adapted basis of ${\bf E}^{p-1}_{(2)}$
chosen before. By the above asymptotic behavior of
$\sigma^{p-1}_{k,l,1}$ and $\sigma^{p-1}_{k,l,2}$ and the
construction of an $L^2$-adapted basis of ${\bf E}^{p-1}_{(2)}$ in
the previous section, one can actually choose the above
$L^2$-adapted basis $\{\sigma_{k,l}^{p-1}\}$ in the way that there
exist an element $\sigma^p_{k,l}$ in some $L^2$-adapted basis of
${\bf E}^p_{(2)}$ satisfying
$\Vert\sigma^p_{k,l}\Vert^2\sim|\log|t_1|||\log|t_2||$ and
$N_1\sigma^p_{k,l}=\sigma^{p-1}_{k,l,1}$ and
$N_2\sigma^p_{k,l}=\sigma^{p-1}_{k,l,2}$. We now go on doing a
reduction for $\phi$. As above, we still discuss the term
containing ${\frac{dt_1}{t_1}}$:
$\phi^{0,1}_{p-1,k,l,1}\wedge{\frac{dt_1}{t_1}}\otimes\sigma^{p-1}_{k,l,1}$.
As before, one still has
$\overline\partial\phi^{0,1}_{p-1,k,l,1}=0$. By the above choice
of the $L^2$-adapted basis, one has
\begin{eqnarray*}
&&~~~\phi^{0,1}_{p-1,k,l,1}\wedge{\frac{dt_1}{t_1}}\otimes\sigma^{p-1}_{k,l,1}\\
&&=\phi^{0,1}_{p-1,k,l,1}\wedge{\frac{dt_1}{t_1}}\otimes N_1\sigma^{p}_{k,l}\\
&&=-\theta(\phi^{0,1}_{p-1,k,l,1}\otimes\sigma^{p}_{k,l})
-\phi^{0,1}_{p-1,k,l,1}\wedge{\frac{dt_2}{t_2}}\otimes N_2\sigma^{p}_{k,l}\\
&&=-\theta(\phi^{0,1}_{p-1,k,l,1}\otimes\sigma^{p}_{k,l})
-\phi^{0,1}_{p-1,k,l,1}\wedge{\frac{dt_2}{t_2}}\otimes\sigma^{p-1}_{k,l,2}\\
&&=-D''(\phi^{0,1}_{p-1,k,l,1}\otimes\sigma^{p}_{k,l})
-\phi^{0,1}_{p-1,k,l,1}\wedge{\frac{dt_2}{t_2}}\otimes\sigma^{p-1}_{k,l,2}.
\end{eqnarray*}
Obviously, $\phi^{0,1}_{p-1,k,l,1}\otimes\sigma^{p}_{k,l}$ is a
local section of $(A^{0,1}\otimes{\bf E}^p)_{(2)}$. So, one can
furthermore assume that $\phi^{1,1}$ has only the term
$\phi^{0,1}_{p-1,k,l,2}\wedge{\frac{dt_2}{t_2}}\otimes\sigma^{p-1}_{k,l}$
with
\[
\Vert{\frac{dt_2}{t_2}}\otimes\sigma^{p-1}_{k,l}\Vert^2\sim|\log|t_1|||\log|t_2||.
\]
Again using the relation
$\overline\partial\phi^{2,0}+\theta\phi^{1,1}=0$, we can finally
assume that $\phi^{1,1}=0$, namely, $\phi$ has only the part of
type $(2, 0)$ up to an exact $D''$-form. Obviously, a $D''$-closed
form of type $(2, 0)$ is holomorphic and $\theta$-closed. This
completes the proof of quasi-isomorphism, and hence Theorem A.
~~~~~~~~~~~~~~~~~~~~~~~~~~~$\Box$


\vskip 1cm \noindent {\bf APPENDIX:} ~In this appendix, the aim is
two-fold. One is to prove that the complex $\{{\text
Gr}_{F}^p{\Omega}^{.}({\bf H}_{\mathbb{C}})_{(2)}, \theta\}$ is
quasi-isomorphic to the complex $([{\text Gr}_F^pA^.({\bf
H}_{\mathbb{C}})]_{(2)}, D'')$ under the inclusion map for $p\ge
0$, if the $L^2$ $\overline{\partial}$-Poincar\'e lemma is always
true. This is actually an application of a general result to this
special case; from the following proof, we can easily abstract
this. The proof is standard. We first establish some notations.

On a topological space $X$, a {\it complex of sheaves} $(K^*, d)$
is given by sheaves of Abelian groups $K^p$ together with sheaf
maps $d$
\[
K^0\stackrel{d}{\to}K^1
\stackrel{d}{\to}\cdots\stackrel{d}{\to}K^p
\stackrel{d}{\to}K^{p+1}\stackrel{d}{\to}\cdots
\]
satisfying $d^2 = 0$. Associated to a complex of sheaves $(K^*,
d)$ are the cohomology sheaves $H^q(K^*)$: Setting, for an open
subset $U\subset X$, $K^q(U) = H^0(U, K^q)$, the presheaf
\[
U\to {\frac{{\text Ker}\{d: K^q(U)\to K^{q+1}(U)\}}{dK^{q-1}(U)}}
\]
gives rise to a sheaf $H^q(K^*)$, whose stalk is
\[
{(H^q(K^*))}_x = \lim_{x\in U}{\frac{{\text Ker}
\{d: K^q(U)\to K^{q+1}(U)\}}{dK^{q-1}(U)}}.
\]
A section $\sigma$ of $H^q(K^*)$ over an open set $U\subset X$ is
given by a covering $\{U_{\alpha}\}$ of $U$ and
${\sigma}_{\alpha}\in K^q(U_{\alpha})$ such that
\begin{eqnarray*}
d\sigma &=& 0,  \\
{\sigma}_{\alpha}-{\sigma}_{\beta} &=& d{\eta}_{\alpha\beta},
~~~{\eta}_{\alpha\beta}\in K^{q-1}(U_{\alpha}\cap U_{\beta});
\end{eqnarray*}
the section is zero in case
\[
{\sigma}_{\alpha} = d{\eta}_{\alpha}, ~~~{\eta}_{\alpha}\in K^{q-1}(U_{\alpha}),
\]
after perhaps refining the given covering. We note that, by the
above definitions, {\it the cohomology sheaves $H^q(K^*) = 0$ for
$q>0$ iff the Poincar\'e lemma holds for the complex of sheaves
$(K^*, d)$}.

\vskip .4cm \noindent {\bf Definition} {\it A map
\[
j: L^*\to K^*
\]
between complexes of sheaves is a quasi-isomorphism if it induces
an isomorphism on cohomology sheaves:
\[
j_*: H^q(L^*)\to H^q(K^*),~~~q\ge 0.
\]}

We now turn to the proof of quasi-isomorphisim, provided that the
$L^2$ $\overline{\partial}$-Poincar\'e lemma be true. By the
definition, we need to show that the two sheaves, $H^k({\text
Gr}_{F}^p{\Omega}^{.}({\bf H}_{\mathbb{C}})_{(2)})$ and
$H^k([{\text Gr}_F^pA^.({\bf H}_{\mathbb{C}})]_{(2)})$, have the
same germs for all $k\ge 0$ and $p\ge 0$. Let $U$ be an open
subset of $\overline M$, then the definition tells us
\begin{eqnarray*}
H^k({\text Gr}_{F}^p{\Omega}^{.}({\bf H}_{\mathbb{C}})_{(2)})(U)
&=& {\frac{{\text Ker}(\Omega^{k}\otimes{\bf
E}^{p-k}\to\Omega^{k+1}\otimes{\bf E}^{p-k-1})(U)}
{{\text Im}(\Omega^{k-1}\otimes{\bf E}^{p-k+1}\to\Omega^{k}
\otimes{\bf E}^{p-k})(U)}}, \\
H^k([{\text Gr}_F^pA^.({\bf H}_{\mathbb{C}})]_{(2)})(U) &=&
{\frac{{\text Ker}([{\text Gr}_F^pA^k({\bf
H}_{\mathbb{C}})]_{(2)}\to[{\text Gr}_F^pA^{k+1}({\bf
H}_{\mathbb{C}})]_{(2)})(U)} {{\text Im}([{\text
Gr}_F^pA^{k-1}({\bf H}_{\mathbb{C}})]_{(2)}\to[{\text
Gr}_F^pA^k({\bf H}_{\mathbb{C}})]_{(2)})(U)}}.
\end{eqnarray*}
Here, "$\to$" in the first formula is under $\theta$, while $\to$
in the second formula is under $D''$. Let $\phi =
\sum_{r+s=k}\phi^{r,s}\in{{\text Ker}([{\text Gr}_F^pA^k({\bf
H}_{\mathbb{C}})]_{(2)} \to[{\text Gr}_F^pA^{k+1}({\bf
H}_{\mathbb{C}})]_{(2)})(U)}$ with $\phi^{r,s}\in (A^{r,s}\otimes
{\bf E}^{p-r})_{(2)}$, which represents the germ of an element in
$H^k([{\text Gr}_F^pA^.({\bf H}_{\mathbb{C}})]_{(2)})$.
Equivalently, $D''\phi = 0$, i.e.,
\begin{eqnarray*}
&&\theta{\phi}^{k,0} = 0,  \\
&&{\overline\partial}{\phi}^{k,0} + \theta{\phi}^{k-1,1} = 0,  \\
&&{\overline\partial}{\phi}^{k-1,1} + \theta{\phi}^{k-2,2} = 0,  \\
&&\cdots  \\
&&{\overline\partial}{\phi}^{k-r,r} + \theta{\phi}^{k-r-1,r+1} = 0,  \\
&&\cdots  \\
&&{\overline\partial}{\phi}^{1,k-1} + \theta{\phi}^{0,k} = 0,  \\
&&{\overline\partial}{\phi}^{0,k} = 0.
\end{eqnarray*}
By means of the last formula and the Poincar\'e lemma for
$\overline\partial$, one has that on $U$ (possibly a smaller open
set), there exists a local section ${\psi}^{0,k-1}\in
(A^{0,k-1}\otimes {\bf E}^{p})_{(2)}$ satisfying
\[
{\phi}^{0,k} = {\overline\partial}{\psi}^{0,k-1}.
\]
Substituting the above formula into
${\overline\partial}{\phi}^{1,k-1} + \theta{\phi}^{0,k} = 0$ and
using $\theta{\overline\partial} + {\overline\partial}\theta = 0$,
one has
\[
{\overline\partial}({\phi}^{1,k-1} - \theta{\psi}^{0,k-1}) = 0.
\]
The same reasoning derives that there exists a local section
${\psi}^{1,k-2}\in (A^{1,k-2}\otimes {\bf E}^{p-1})_{(2)}$ on a
(possibly smaller than $U$) open subset, satisfying
\[
{\phi}^{1,k-1} = {\overline\partial}{\psi}^{1,k-2} + \theta{\psi}^{0,k-1}.
\]
Substituting the above formula into
${\overline\partial}{\phi}^{2,k-2} + \theta{\phi}^{1,k-1} = 0$ and
using $\theta{\overline\partial} + {\overline\partial}\theta = 0$
and $\theta\wedge\theta = 0$, one has
\[
{\overline\partial}({\phi}^{2,k-2} - \theta{\psi}^{1,k-2}) = 0.
\]
Inductively, one has that there exists a sequence of germs of
local sections ${\psi}^{r,k-r-1}\in (A^{r,k-r-1}\otimes {\bf
E}^{p-r})_{(2)}$, $2\le r\le k-1$, satisfying
\[
{\phi}^{r,k-r} = {\overline\partial}{\psi}^{r,k-r-1} + \theta{\psi}^{r-1,k-r}.
\]
Setting $r=k-1$ and substituting the above formula into
${\overline\partial}{\phi}^{k,0} + \theta{\phi}^{k-1,1} = 0$, one
has
\[
{\overline\partial}({\phi}^{k,0} - \theta{\psi}^{k-1,0}) = 0.
\]
So, ${\phi}^{k,0} - \theta{\psi}^{k-1,0}$ is a local holomorphic
$k$-form with values in ${\bf E}^{p-k}$ (if $p-k\ge 0$), denoted
by $h^k$, i.e.,
\[
{\phi}^{k,0} = \theta{\psi}^{k-1,0} + h^k,
\]
and hence $\theta h^k = 0$, which implies that $h^k$ represents
the germ of some element in $H^k({\text Gr}_{F}^p{\Omega}^{.}({\bf
H}_{\mathbb{C}})_{(2)})$ . Summing up
$\sum_{r=1}^k{\phi}^{r,k-r}$, one has
\begin{eqnarray*}
\phi
&=& {\overline\partial}{\psi}^{0,k-1}+\theta{\psi}^{0,k-1}+
{\overline\partial}{\psi}^{1,k-2}+\cdots +
\theta{\psi}^{k-2,1}+{\overline\partial}{\psi}^{k-1,0}+
\theta{\psi}^{k-1,0}+h^k   \\
&=& D''({\psi}^{0,k-1}+{\psi}^{1,k-2}+\cdots + {\psi}^{k-2}+
{\psi}^{k-1,0})+h^k,
\end{eqnarray*}
which implies that the germ represented by $\phi$ in the sheaf
$H^k([{\text Gr}_F^pA^.({\bf H}_{\mathbb{C}})]_{(2)})$ can be
identified with the germ represented by $h^k$ in the sheaf
$H^k({\text Gr}_{F}^p{\Omega}^{.}({\bf H}_{\mathbb{C}})_{(2)})$.
We need to show that this identification is well-defined, namely,
if the germ in $H^k([{\text Gr}_F^pA^.({\bf
H}_{\mathbb{C}})]_{(2)})$ represented by $\phi$ can also be
identified with the germ in $H^k({\text Gr}_{F}^p{\Omega}^{.}({\bf
H}_{\mathbb{C}})_{(2)})$ represented by another local holomorphic
$k$-form with values in ${\bf E}^{p-k}$ $h^k_1$ (that is to say,
there exists a local section $\psi_1 \in [{\text
Gr}_F^pA^{k-1}({\bf H}_{\mathbb{C}})]_{(2)}$ with $\phi =
D''\psi_1 + h^k_1$), then $h^k$ and $h^k_1$ represent the same
germ in $H^k({\text Gr}_{F}^p{\Omega}^{.}({\bf
H}_{\mathbb{C}})_{(2)})$. The proof of this assertion is
completely similar to the previous arguments-----the same steps
let us obtain a local holomorphic $k-1$-form $h^{k-1}$ with values
in ${\bf E}^{p-k+1}$ satisfying
\[
h^k_1 - h^k = \theta{h^{k-1}}.
\]
Thus, we get a natural map $j^*$ from $H^k([{\text Gr}_F^pA^.({\bf
H}_{\mathbb{C}})]_{(2)})$ to $H^k({\text
Gr}_{F}^p{\Omega}^{.}({\bf H}_{\mathbb{C}})_{(2)})$. The assertion
of uniqueness also shows that the map $j^*$ is surjective. The
proof of injectivity is easy: Suppose that
$j^*([\phi_1]) = j^*([\phi_2])$ for $\phi_i\in    \\
{{\text Ker}([{\text Gr}_F^pA^k({\bf H}_{\mathbb{C}})]_{(2)}
\to[{\text Gr}_F^pA^{k+1}({\bf H}_{\mathbb{C}})]_{(2)})(U)}$,
$i=1,2$. By the above arguments, there exist some local sections
$\psi_i\in [{\text Gr}_F^pA^{k-1}({\bf H}_{\mathbb{C}})]_{(2)}$
and some local holomorphic $k$-forms $h^{k}_i$ with values in
${\bf E}^{p-k}$ satisfying $\phi_i = D''\psi_i + h^k_i$ for
$i=1,2$. On the other hand, since $j^*([\phi_1]) = j^*([\phi_2])$,
so there exists a local holomorphic $k-1$-form $h^{k-1}$ with
values in ${\bf E}^{p-k+1}$ satisfying $h^k_1-h^k_2 =
\theta{h^{k-1}}$. So, one has
\[
 \phi_1-\phi_2 = D''(\psi_1-\psi_2) + (h^k_1-h^k_2) = D''(\psi_1-\psi_2+h^{k-1}).
\]
Therefore, as germs, $[\phi_1]=[\phi_2]$. This completes the proof
of the quasi-isomorphism theorem.

The second aim of the appendix is to state that the
hypercohomology of a complex of fine sheaves is isomorphic to the
cohomology of the corresponding complex of global sections. First,
we need to show what the hypercohomology of a sheaves complex is.
As before, let $(K^*, d)$ be a complex of sheaves with
differentials $d$ on a topological space $X$. Take a covering
${\underline U}=\{U_{\alpha}\}$ of $X$ and let $C^p({\underline
U}, K^q)$ be the $\breve{\text C}$ech cochains of degree $p$ with
values in $K^q$. The two operators
\begin{eqnarray*}
&&\delta: C^p({\underline U}, K^q)\to C^{p+1}({\underline U},
K^q)\\
&&d: C^p({\underline U}, K^q)C^p({\underline U}, K^{q+1}),
\end{eqnarray*}
satisfy $\delta^2=d^2=0$ and $d\delta+\delta d=0$; and one has a
double complex
\[
\{C^{p,q}:=C^p({\underline U}, K^q); \delta, d\}.
\]
Put $C^n({\underline U})=\bigoplus_{p+q=n}C^p({\underline U},
K^q)$ and $D=\delta+d$, then $\{C^*, D\}$ is a complex. A
refinement ${\underline U}'$ of $\underline U$ induces mappings
\begin{eqnarray*}
&&C^p({\underline U}, K^q)\to C^p({\underline U}', K^q),\\
&&H^*(C^*(\underline U), D)\to H^*(C^*({\underline U}'), D),
\end{eqnarray*}
and then we define the hypercohomology of $(K^*, d)$ as
\[
{\mathbb{H}}^*{X, \{K^*, d\}}=\lim_{\underline
U}H^*(C^*(\underline U), D).
\]
For the complex of sheaves $(K^*, d)$, one can consider the
complex of its global sections. Putting $\Gamma({K^p})=H^0(X,
K^p)$, one then has a complex $\{{\Gamma({K^p})}, d\}$ and its
usual cohomology $H^*(\{{\Gamma({K^p})}, d\})$. Then one has

\vskip .3cm \noindent {\bf Theorem.} {\it If the complex of
sheaves $(K^*, d)$ is a complex of fine sheaves, then
\[
{\mathbb{H}}^*({X, \{K^*, d\}})\cong H^*(\{{\Gamma({K^p})}, d\}).
\]}

\bigskip

\noindent
Juergen Jost,\\
Max-Planck Institute for Mathematics in the Sciences, Inselstr. 22-26, 04103 Leipzig, Germany;\\
Yi-Hu Yang,\\
Department of Applied Mathematics, Tongji University Shanghai 200092, China;\\
Kang Zuo,\\
Department of Mathematics, the Chinese University of Hong Kong, Hong Kong, China.


\begin{thebibliography}{99}


\bibitem{av} A. Andreotti and E. Vesentini, Carleman estimates
for the Laplace-Beltrami equation in complex manifolds. Publ.
Math. I.H.E.S., {\bf 25} (1965), 81-130.

\bibitem{ca} E. Cattani, Mixed Hodge structures, compactifications
and monodromy weight filtration, in {\it Topics in Transcendental
Algebraic Geometry}, ed P. Griffiths, Annals of Mathematics
Staudies, Vol. {\bf 106}, 75-100.

\bibitem{cg}Cornalba and Griffiths, Analytic cycles and vector
bundles on non-compact algebraic varieties, Invent. Math. {\bf 28}
(1975), 1--106.

\bibitem{ck} E. Cattani, A. Kaplan, Polarized mixed Hodge
structures and the local monodromy of a variation of Hodge
structures, Inventiones Math. {\bf 67} (1982), 101-115.

\bibitem{cks1} E. Cattani, A. Kaplan and W. Schmid, Degeneration
of Hodge structures, Annals of Mathematics, {\bf 123} (1986),
457-535.

\bibitem{cks2} E. Cattani, A. Kaplan and W. Schmid, $L^2$ and
intersection cohomologies for a polarizable variation of Hodge
structure, Inventiones Math., {\bf 87} (1987), 217-252.

\bibitem{d} P. Deligne, Th\'eorie de Hodge. II, Publ. Math. IHES {\bf 40} (1972), 5-57

\bibitem{dem} J.-P. Demailly, $L^2$ Vanishing Theorems for Positive
Line Bundles and Adjunction Theory, LNM, {\bf 1646}, 1994.

\bibitem{gri}P. Griffiths, Periods of integrals on algebraic manifolds. III.
Some global differential-geometric properties of the period
mapping, Inst. Hautes Etudes Sci. Publ. Math. {\bf 38} (1970)
125--180.

\bibitem{gh} P. Griffiths and Harris, {\it Principles of algebraic
geometry},
Pure and Applied Mathematics. Wiley-Interscience [John Wiley \&
Sons], New York, 1978.

\bibitem{ho} L. H\"ormander $L^2$ estimates and existence theorems for
the $\overline{\partial}$ operator, Acta Mathematica, {\bf 113}
(1965), 89-152.

\bibitem{sch} W. Schmid, Variation of Hodge structure: The singularities
of period mapping, Inventiones Math., {\bf 22}(1973),211-319.

\bibitem{si1} C. Simpson, Harmonic bundles on noncompact curves,
J. Amer. Math. Soc., {\bf 3}(1990), 713-770.

\bibitem{si2} C. Simpson, Higgs bundles and local systems,
IHES Publ. {\bf 75} (1992), 5-95.

\bibitem{zuc1} S. Zucker, Generalized intermediate Jacobians and the
theorem on normal functions, Inventiones Math., {\bf 33} (1976),
185-222.

\bibitem{zuc2} S. Zucker, Hodge theory with degenerating coefficients:
$L^2$-cohomology in the Poincar\'e metric, Annals of Mathematics,
{\bf 109}(1979), 415-476.

\bibitem{zuo} K. Zuo, On the negativity of kernels of Kodaira-Spencer maps
on Hodge bundles and applications, {\it Kodaira's issue}, Asian J.
Math. {\bf 4} (2000), no. 1, 279--301.


\end{thebibliography}
\end{document}